\documentclass[leqno]{article}
\usepackage{amssymb,amsmath,amscd,amsthm,amsfonts,graphicx}
\usepackage[dvips]{pstcol}
\usepackage[all]{xy}

\newtheorem{theorem}{Theorem}[section]
\newtheorem{lemma}[theorem]{Lemma}
\newtheorem{proposition}[theorem]{Proposition}
\newtheorem{corollary}[theorem]{Corollary}

\newtheorem{definition}[theorem]{Definition}

\def\bs{\bigskip}
\def\ms{\medskip}
\def\semi{\hbox{ $\times $ \kern-.972em \raise.12719em\hbox{ $_{^|}$}  }}

\def\be{\begin{enumerate}}
\def\ee{\end{enumerate}}
\def\bi{\begin{itemize}}
\def\ei{\end{itemize}}

\title{Conjugacy in Garside Groups I: \\ Cyclings, Powers, and Rigidity \\
   }\author{Joan S. Birman\footnote{Partially supported by
the U.S. National Science Foundation under Grants DMS-9973232 and
0405586.} \and Volker Gebhardt \and Juan
Gonz\'alez-Meneses\footnote{Partially supported by
MTM2004-07203-C02-01 and FEDER.} }
\date{May 8, 2006}

\begin{document}

\maketitle


\begin{abstract}

\noindent  In this paper a relation between iterated cyclings and
iterated powers of elements in a Garside group is shown. This yields
a characterization of elements in a Garside group having a rigid
power, where `rigid' means that the left normal form changes only in
the obvious way under cycling and decycling. It is also shown that,
given $X$ in a Garside group, if some power $X^m$ is conjugate to a
rigid element, then $m$ can be bounded above by $||\Delta||^3$.  In
the particular case of braid groups $\{B_n, \ n\in {\mathbb N}\}$,
this implies that a pseudo-Anosov braid has a small power whose
ultra summit set consists of rigid elements. This solves one of the
problems in the way of a polynomial solution to the conjugacy
decision problem (CDP) and the conjugacy search problem (CSP) in
braid groups. In addition to proving the rigidity theorem, it will
be shown how this paper fits into the authors' program for finding a
polynomial algorithm to the CDP/CSP, and what remains to be done.
\end{abstract}

\setlength{\parskip}{0pt}

\tableofcontents

\setlength{\parskip}{1ex plus 0.5ex minus 0.2ex}

\bs

\section{Introduction}

Braid groups $B_n, \ n =1,2,3,\dots$, were introduced in a
foundational paper by Emil Artin \cite{Artin1925} in 1925.  In it
Artin gave the well-known presentation:
\begin{equation}
\label{equation:classical presentation}
B_n= \left< \sigma_1,\ldots, \sigma_{n-1} \left|
\begin{array}{ll}
 \sigma_i \sigma_j = \sigma_j \sigma_i  & \mbox{if } \; |i-j|>1,  \\
 \sigma_i \sigma_j \sigma_i = \sigma_j \sigma_i \sigma_j &  \mbox{if }\; |i-j|=1.\ \end{array}\right.
 \right>.
\end{equation}
\begin{figure}[htpb!]
$\centerline{\xymatrix@C=.3cm@R=1cm{
 1 \ar[dd] & \ar@{}[dd]|\cdots  &  i-1 \ar[dd] & i \ar@(d,u)[ddr]|\hole & i+1 \ar@(d,u)[ddl]& i+2 \ar[dd] & \ar@{}[dd]|\cdots & n \ar[dd]\\ \\
  &  &  &  &  & & &
} }$
 \caption{ The elementary braid $\sigma_i$. }
\label{figure:braids-elem}
\end{figure}
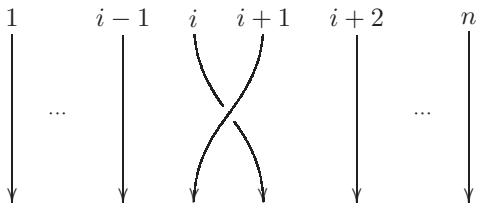

The elementary braid $\sigma_i$ is depicted in Figure
\ref{figure:braids-elem}. To study $B_n$ Artin used the fact that
there is a canonical homomorphism $\pi:B_n\to \Sigma_n$, where the
image is the symmetric group, defined by sending a braid to the
associated permutation of its endpoints.  He went on to uncover the
structure of the kernel of $\pi$, and used what he learned to solve
the {\it word problem} in $B_n$: to decide, for arbitrary words
$X,Y$ in the generators and their inverses, whether they represent
the same element of $B_n$.  Artin also posed the {\it conjugacy
decision problem} (CDP): to decide whether, for arbitrary $X,Y \in
B_n$, there exists $Z\in B_n$ such that $Y = Z^{-1}XZ$.  A different
but related problem, the {\it conjugacy search problem} (CSP) asks
to find $Z$, provided that one knows that it exists.

During the period 1925-1969 various efforts were made to solve the
conjugacy problem, building on techniques which had been introduced
in \cite{Artin1925}, but there was no significant progress. Then, in
1969 F. Garside \cite{Garside} brought completely new techniques to
bear, looking at $B_n$ in a very new way which stressed the
similarity of its combinatorics to those of $\Sigma_n$,  rather than
focussing on ker($\pi:B_n\to\Sigma_n$).  Garside succeeded in
solving both the word and conjugacy search problems simultaneously,
and in a unified way.  His methods were soon shown to apply to other
groups too \cite{B-S},\cite{Deligne}, and over the years broadened
to an entire class of groups which subsequently became known as {\it
Garside groups}.   The ideas that Garside introduced, and their
subsequent improvements, are the subject of this paper, which  is
the first in a series with the unifying title `Conjugacy in Garside
groups I, II, III,$\dots$  They have a common goal: to improve
Garside's algorithm for the CDP/CSP in a Garside group to obtain, in
the particular case of the braid group $B_n$, an algorithm which is
polynomial both in $n$ and an appropriate measure $||X||,||Y||$ of
the complexity of $X$ and $Y$. This would have implications as
regards the security of certain codes in public key cryptography
\cite{AAG,KLC}.

The existence of such a {polynomial algorithm} for the word problem
in all Garside groups is now known via the work of
\cite{El-M,Epstein} for the braid groups. In \cite{D-P} and
\cite{Dehornoy-2002} the class of Garside groups is defined in a
more general setting, and shown to be biautomatic.  It is a
consequence of the way the definitions were chosen in \cite{D-P,
Dehornoy-2002} that in fact, Garside's algorithm  solves the word
and CDP/CSP's in all Garside groups.  Our work in this paper is a
step in a program that we have developed to prove that the CDP/CSP
in $B_n$ is polynomial in both $n$ and $||X||$. But all results in
this paper are valid in every Garside group, except the results in
$\S$\ref{subsection:Consequences for pseudo-Anosov braids} and
Theorem~\ref{T:bound for pseudo-Anosov}, where we consider
applications of these results to the special case of braid groups,
in particular to pseudo-Anosov braids.

Before we can state exactly what we do in this paper, and describe
it in context, we need to set up necessary notation and review the
known results and techniques.  The combinatorial structure that we
will use, and the new structure that we have uncovered, is quite
complicated and, we think, interesting. In order to make this paper
accessible to non-experts we give details and examples which those
who are acquainted with the literature will probably wish to bypass
quickly, moving on to $\S$\ref{subsection:summary of results}, where
we describe the essential content of this paper and its context in
our larger goal, and thence to $\S$\ref{section:cyclings and
powers}, where our new contributions begin.

\ms

{\bf Acknowledgements:} J.Birman and J. Gonzalez-Meneses, who were
working together, and V. Gebhardt, became acquainted with each
other's partial results at a conference in the Banff International
Research Station for Mathematical Innovation and Discovery, in
October 2004. There was some overlap, and also some recognition that
differing viewpoints could lead to progress, so they decided to pool
forces at that time. The three authors thank the PIMS, MSRI, MITACS
and IM-UNAM for their wisdom in sponsoring international conferences
which foster exactly this kind of fruitful interchange and
collaboration.

J. Birman thanks the Project MTM2004-07203-C02-01 of the Spanish
Ministerio de Ciencia y Tecnolog\'\i a for hosting her visit to
Seville in November 2004, so that she and J. Gonzalez-Meneses could
work together on this project.

J. Gonz\'alez-Meneses thanks the project MTM2004-07203-C02-01 and
the Columbia University Department of Mathematics for hosting his
two visits to New York, in July 2004 and March-April 2006.

V. Gebhardt thanks the department of Algebra of the University of
Seville, and the Junta de Andaluc\'\i a, for funding his visit to
Seville in January 2006.

The work in this paper and \cite{BGGM-II, BGGM-III}, which are being
completed at this writing,  was done simultaneously and
independently from the work of S.J. Lee and E.K.Lee in
\cite{Lee-Lee,Lee-Lee2,Lee-Lee3}. We first became aware of that work
when we were in the process of writing up this one.  Some of this
work, notably Proposition~\ref{P:SU is not emptyset} below and the
results in \cite{BGGM-II}, were reported on in talks at conferences
in Banff in October 2004 and Luminy in June 2005.

\subsection{Garside groups}
\label{subsection:Garside groups}
Among the known equivalent definitions of Garside groups, we  use the one which was suggested to us by John Crisp \cite{Crisp}, because it seems the most natural of the many possible definitions.  A group $G$ is said to be a {\em Garside group} if it satisfies
properties (A), (B) and (C) below:
\be
\item [(A)] $G$ admits a lattice order $(G,\preceq,\vee,\wedge)$, invariant under
left-multiplication.

This means that there is a partial order $\preceq$ on the elements
of $G$ such that $a\preceq b$ implies $ca\preceq cb$ for every $c\in
G$. Also, every pair of elements $s,t\in G$  admits a unique lcm
$s\vee t$ and a unique gcd $s\wedge t$ with respect to $\preceq$.
This partial order $\preceq$ defines a submonoid $P\subset G$,
called the {\em positive cone} of $G$, defined by $P=\{p\in G;\;
1\preceq p\}$.  Notice that the invariance of $\preceq$ under
left-multiplication implies that $P\cap P^{-1}=\{1\}$, and also that
$  a \preceq b \quad \Longleftrightarrow \quad a^{-1}b \in P. $
Hence the submonoid $P$ determines the partial order $\preceq$, so
we shall equally talk about {\em the lattice} $(G,P)$. We remark
that if $a,b\in P$ then $a\preceq b$ if and only if $a$ is a prefix
of $b$, that is, there exists $c\in P$ such that $ac=b$. This is why
$\preceq$ is sometimes called the {\em prefix} order.

There is also a related {\it suffix} order, defined by $b\succeq a$ if $ba^{-1}\in P$. It is important that $a\preceq b$ does not imply that $b\succeq a$. Sometimes we will get genuinely new information by using both orderings, even when the proofs are little more than copies of one-another.

\item [(B)] There exists an element $\Delta \in P$, called the {\it Garside element}, satisfying:
\begin{enumerate}
\item [(a)]
The interval $[1,\Delta]=\{s\in G;\; 1\preceq s \preceq \Delta\}$
generates $G$. Its elements are called the {\it simple elements} of
$G$.  We shall always assume that $[1,\Delta]$ is finite, that is,
that $G$ has {\it finite type}.

 \item [(b)] Conjugation by $\Delta$ preserves the positive cone $P$: $\Delta^{-1}
 P \Delta = P$.
\end{enumerate}

We remark that if $\Delta$ satisfies both (a) and (b), then
$[1,\Delta]$ also generates $P$ as a monoid, which is one of the
properties usually required in the definition of a Garside element.

\item  [(C)] The monoid $P$ is {\it atomic}.

This means that for every $x\in P$ there exists an upper bound on
the length of a (strict) chain $1\prec x_1 \prec \cdots \prec x_r
=x$.  In other words, if we define the {\it atoms} of $G$ as the
elements $a\in P$ which cannot be decomposed in $P$ (there are no
nontrivial $b,c\in P$ such that $a=bc$), then for every $x\in P$
there exists an upper bound on the number of atoms in a product $x =
a_1a_2\cdots a_r$ with each $a_i$ an atom. In particular,  if $P$ is
atomic, one can define the {\em length} of an element $x\in P$ as
the maximal length of such a chain, that is,
$$
||x||=max\{n:\; x=a_1a_2\cdots a_n, \; \mbox{ where } \;a_i\in
P\backslash \{1\}\; \}.
$$
Notice that the atoms generate $G$.

\ee

These data determine a {\em Garside structure} on $G$, which may be
defined as follows: Let $G$ be a countable group, $P$ be a
submonoid, and $\Delta\in P$. The triple $(G,P,\Delta)$ is said to
be a {\em (finite type) Garside structure} on $G$ if  $(G,P)$ is a
lattice, $\Delta$ is a Garside element (with $[1,\Delta]$ finite),
and $P$ is atomic. We remark that a given group $G$ may admit more
than one Garside structure.

{\bf Example 1:} Our first example is very simple: We consider the braid group $B_3$ and its two known Garside structures:

1A. The classical Garside structure is associated to the
presentation  (\ref{equation:classical presentation}) of $B_3$. The
Garside element is $\sigma_1\sigma_2\sigma_1 =
\sigma_2\sigma_1\sigma_2$. The elements in $P$ correspond to the
braids in which all crossings are positive.   The atoms are
$\sigma_1$ and $\sigma_2$.

1B.  If we set $x=\sigma_1, y=\sigma_2, z =
\sigma_2\sigma_1\sigma_2^{-1}$ we get the presentation $\left< x,y,z
\; | \; xy = yz = zx \right>$.    The Garside element is now
$\Delta= xy$.  See \cite{Xu-1992a} for the way in which this
structure was used to solve the shortest word problem in $B_3$ and
to give an algorithm for determining the genus of knots and links
which are closed 3-braids.  This Garside structure was generalized
to all $n$ in \cite{BKL-1998}.

{\bf Example 2: Free abelian groups of finite rank.} This is another very simple example of a Garside group:
$$
 \mathbb Z^n= \left< \; x_1,\ldots, x_n \quad | \quad x_ix_j=x_jx_i, \quad i<j
 \right>.
$$
The positive cone is
$$
 \mathbb  N^n=\{x_1^{e_1}\cdots x_n^{e_n} \; ; \quad e_i\geq 0, \; \forall
 i\}.
$$
The Garside element is $\Delta=x_1\cdots x_n$, and the simple
elements have the form $x_1^{e_1}\cdots x_n^{e_n}$ where
$e_i\in\{0,1\}$ for every $i=1,\ldots,n$. Hence there are $2^n$
simple elements.\ms

\begin{figure} [ht]
\centerline{\includegraphics{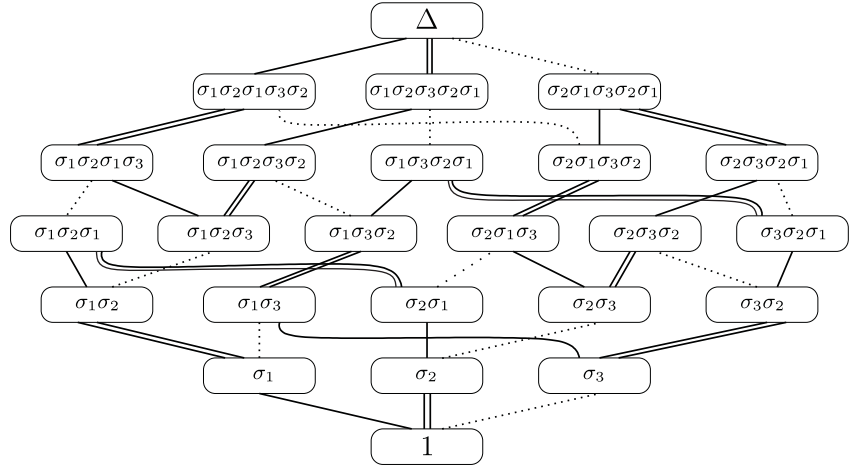}}
 \caption{The lattice of simple elements in $B_4^+$. They are $4! =24$
elements.}\label{figure:simpleB4}
\end{figure}

{\bf Example 3: The braid group $B_n$, with the classical Garside structure:}
Garside used the presentation (\ref{equation:classical presentation}).  The usual Garside
structure in this group is determined by $(B_n, B_n^+,\Delta)$,
where $B_n^+$ is the monoid of {\em positive braids}, consisting of
the elements in $B_n$ that can be written as a product of
$\sigma_i$'s with no $\sigma_i^{-1}$, and
\begin{equation}
\label{Garside's Delta in Bn}
\Delta= (\sigma_1)
(\sigma_2 \sigma_1) (\sigma_3 \sigma_2 \sigma_1)\cdots (\sigma_{n-1}
\cdots \sigma_1)
\end{equation}
is a half-twist on all of the strands.  The atoms are $\sigma_1,\dots,\sigma_{n-1}$.
The elements in $P$ correspond to the braids in which all
crossings are positive. The Garside element $\Delta$ can be characterized as the only positive braid in
which every pair of strands cross exactly once. The simple elements,
in the case of $B_n$, are the positive braids in which every pair of strands
cross {\em at most} once. It follows that every simple element
corresponds to a permutation on the set of $n$ elements (the
strands). Hence there are $n!$ simple elements in $B_n^+$, with this
Garside structure.  Figure~\ref{figure:simpleB4} shows the Hasse
diagram representing the lattice of simple elements in $B_4^+$. In
this diagram, an element $a$ is joined by a line to an element $b$
in the upper row if and only if $a\preceq b$. Moreover, each line
type corresponds to a right multiplication by an atom: a single line
corresponds to $\sigma_1$, a double line to $\sigma_2$ and a dotted
line to $\sigma_3$. This lattice of simple elements determines the
whole Garside structure of the group.

\ms

{\bf Example 4: Spherical type Artin-Tits groups:}  \cite{Bourbaki}
The previous three examples were particular cases of Artin-Tits
groups. All Artin-Tits groups of spherical type are known to be
Garside groups~\cite{B-S}. Given a finite set $S$, a {\em Coxeter
matrix} over $S$ is a symmetric matrix $M=(m_{st})_{s,t\in S}$,
where $m_{ss}=1$ for all $s\in S$ and $m_{s,t}\in
\{2,3,\cdots,\infty\}$. Every Coxeter matrix $M$ defines a group
$A_M$ given by the following presentation:
$$
   A_M = \left< \;S \; \left| \quad \underbrace{stst\cdots}_{m_{st}\;\mbox{\scriptsize
terms}} = \underbrace{tsts\cdots}_{m_{st}\;\mbox{\scriptsize
terms}}, \quad \mbox{ for all } \; s,t\in S
   \right>\right.,
$$
where $m_{st}=\infty$ means that there is no relation involving $s$
and $t$. The group $A_M$ is called the {\em Artin-Tits group}
associated to $M$, also called {\em Artin group} or {\em generalized
braid group}.

If one adds to the above presentation the relations $s^2=1$ for all
$s\in S$, one obtains the group $W_M$, called the {\em Coxeter
group} associated to $M$. An Artin-Tits group is said to be of {\em
spherical type} if its corresponding Coxeter group is finite.

The usual Garside structure in these groups is given by $(A_M,
A_M^+,\Delta)$, where $A_M^+$ is the monoid of {\em positive
elements}, consisting of products of elements of $S$ (the above
presentation of $A_M$, considered as a monoid presentation, gives
$A_M^+$), and the Garside element $\Delta$ is defined as follows.
The set of generators $S$ can be decomposed into two sets $S=S_1\cup
S_2$, where elements contained in the same set $S_i$ commute. This
decomposition can be easily obtained from the {\em Coxeter graph} of
the group. The reader is referred to \cite{Bourbaki} for the
definition of a Coxeter graph and its associated Coxeter matrix, and
also for a list of the Coxeter graphs associated to the finite
Coxeter groups. If the Coxeter graph $\Gamma$ is connected, there is
only one possible decomposition of $S$ in the above sense. Now
define
$$\displaystyle
\Delta_1=\prod_{s\in S_1}{s} \quad  {\rm and} \quad  \displaystyle
\Delta_2=\prod_{s\in S_2}{s}.$$

Then one has
$$ \Delta= \underbrace{\Delta_1 \Delta_2 \Delta_1 \Delta_2 \cdots }_{h \mbox{
 terms}},$$
where $h$ is the {\em Coxeter number} of the corresponding Coxeter
group. The Coxeter numbers corresponding to the spherical type
Artin-Tits groups are the following.

$$
\begin{tabular}{|c||c|c|c|c|c|c|c|c|c|c|c|}
\hline
 Type & $A_l$ & $B_l$ & $D_l$ & $E_6$ & $E_7$ & $E_8$ & $F_4$ & $G_2$
 & $H_3$ & $H_4$ & $I_2(p)$ \\ \hline
 $h$ & $l+1$ & $2l$ & $2l-2$ & 12 & 18 & 30 & 12 & 6 & 10 & 30 & $p$
 \\
\hline
\end{tabular}
$$

\bigskip

As an example, the Garside element of the the spherical Artin-Tits
group of type $B_l$ is $\Delta=((s_1s_3\cdots s_l)(s_2s_4\cdots
s_{l-1}))^{l}$ when $l$ is odd, and $\Delta=((s_1s_3\cdots
s_{l-1})(s_2s_4\cdots s_{l}))^{l}$ when $l$ is even. Notice that the
Artin-Tits monoid of type $A_l$ is precisely the Artin braid monoid
on $l+1$ strands. Notice also that the Garside structure given by this
construction coincides with the original Artin structure for braid
groups described above.

We remark that every spherical type Artin-Tits group admits another
Garside structure, discovered in~\cite{bessis}, called the {\it
dual} Garside structure. In the case of braid groups, the dual
Garside structure is precisely the one discovered
in~\cite{BKL-1998}.

\ms
 {\bf Example 5: Torus knot groups.}  The fundamental group of the
complement of a $(p,q)$-torus knot, where $p,q>1$ are coprime, is
given by the following presentation:
$$
\left< x,y \; | \; x^p=y^q  \right>.
$$
If we consider the monoid given by this presentation, it is a
Garside monoid with Garside element $\Delta= x^p$.

\ms

 {\bf Example 6:}  The following two examples of Garside
groups do not belong to a class of known groups, but they have
interesting properties which are not satisfied by the groups in the
previous examples. They were discovered and studied by Picantin
in~\cite{Picantin-PhD}. In both cases, we give presentations of the
groups which, considered as monoid presentations, yield the
corresponding Garside monoids. Hence we shall only define the
Garside element, in each case.
\begin{enumerate}

\item $G = \left< x,y,z \; | \; xzxy=yzx^2, \; yzx^2z = zxyzx = xzxyz
\right>$. The Garside element is $\Delta= xzxyzx$. In most examples
of Garside groups, the Garside element $\Delta$ is the least common
multiple (with respect to $\preceq$) of the atoms. In this example,
since the relations are homogeneous, the atoms are just the letters
$x$, $y$ and $z$, and one has ${\rm lcm}_\preceq (x,y,z)=xzxyz$. But
$xzxyz\neq \Delta$. Indeed, since conjugation by $\Delta$ must
preserve the set of atoms, all atoms must be left and right divisors
of $\Delta$, but we have $xzxyz\not\succeq y$, hence ${\rm
lcm}_\preceq (x,y,z)=xzxyz \neq {\rm lcm}_\succeq (x,y,z)$. This is
an example of a Garside monoid in which the ${\rm lcm}_{\preceq}$ of
the atoms is not a Garside element.

\item $G=\left< x,y \; | \; xyxyx = y^2 \right>$. Garside element
$\Delta = y^3$.  In this case, $\mbox{lcm}_{\preceq}(x,y)=
\mbox{lcm}_{\succeq}(x,y)=y^2$, but $y^2$ is not a Garside element.
Indeed, since conjugation by $\Delta$ must preserve the set of
simple elements, the set of positive left-divisors and the set of
positive right-divisors of $\Delta$ must coincide. But this does not
happen for $y^2$. For instance $xyxy\preceq y^2$ but $y^2
\not\succeq xyxy$.  This is also an example of a Garside monoid in
which the relations are not homogeneous, hence the length of a
positive element is not given by the letter length of any
representative.
\end{enumerate}

\ms
 {\bf Construction of new Garside monoids.} We already
provided several examples of Garside monoids and groups. Using these
monoids as building blocks, one can construct new Garside monoids
and groups thanks to the following result. In~\cite{Picantin-PhD}
there is a definition of the so called {\em crossed product} of
monoids, which also allows to construct new Garside monoids.

\begin{theorem}{\em \cite{Picantin-PhD}}
The crossed product of Garside monoids is a Garside monoid.
\end{theorem}

An example of crossed product, given in~\cite{Lee}, is the
semidirect product $\mathbf Z\ltimes G^n$, where the action of
$\mathbb Z$ on the free product $G^n$ ($G$ is a Garside group) is
given by cyclic permutations of coordinates.

\subsection{Solving the word and conjugacy problems in Garside groups}
\label{subsection:Solving the word and conjugacy problem in Garside
groups}

From now on, we will fix a Garside group $G$ with a finite type
Garside structure $(G,P,\Delta)$.  We will show how to solve the
word problem, giving a well known normal form for elements in a
Garside group.  The basic reference is \cite{El-M}. While everything
in that paper relates to the braid groups, most of it generalizes
easily to arbitrary Garside groups,  which were singled out as a
class several years later \cite{D-P, Dehornoy-2002}. \ms

\begin{definition}{\bf (Left normal form):} {\rm Given $X\in G$, we will say that
a decomposition $X=\Delta^p x_1\cdots x_r$ ($r\geq 0)$ is the {\em
left normal form} of $X$ is it satisfies
\begin{enumerate}
 \item $p\in \mathbb Z$ is maximal such that $\Delta^p\preceq X$.  That is,  $x_1\cdots x_r\in P$
and $\Delta\not\preceq x_1\cdots x_r$
 \item $x_i=(x_i\cdots x_r)\wedge \Delta$, for $i=1,\ldots,r$. That  is,  $x_i$ is the biggest {\em simple} prefix of $x_i\cdots
x_r$.
\end{enumerate}
}
\end{definition}

One can also show by induction that $\Delta^p x_1\cdots x_i =
X \wedge \Delta^{p+i}$, for $i=1,\ldots,r$.

It is known that one can check whether a given decomposition
$X=\Delta^p x_1\cdots x_r$ is a left normal form by looking at each
pair of consecutive factors $x_ix_{i+1}$. We say that a pair of
simple elements $a,b\in [1,\Delta]$ is {\em left weighted}, if the
product $ab$ is in left normal form as written, that is, if
$a=(ab)\wedge \Delta$. Then $\Delta^p x_1\cdots x_r$ is a left
normal form if and only if $x_1\neq \Delta$ and each pair
$x_ix_{i+1}$ is left weighted.

Notice that, if we consider the set of simple elements as a set of
generators for $G$, then the decomposition defined above is a normal
form in the usual sense, that is, a unique way to write any element
of $G$ as a product of the generators and their inverses.  If one
wishes to obtain a normal form with respect to any other set of
generators (for example the set of atoms), one just needs to choose
a unique way to write each simple element in terms of the desired
generators, and replace this in the left normal form. \ms

We now give several standard terms that will be  needed to work with
Garside groups.\ms

 If $X=\Delta^p x_1\cdots x_r$ is in left normal form, the {\em
infimum, supremum} and {\em canonical length} of $X$, are defined by
$\inf(X)=p$, $\sup(X)=p+r$ and $\ell(X)=r$, respectively.

The {\it shift map}  $\tau$ is the inner automorphism $\tau:\: G \rightarrow G$
given by $\tau(x)=\Delta^{-1} x \Delta$.

Given a simple element $x$, we define $x^*=x^{-1}\Delta$. That is,
$x^*$ is the only simple element such that $x\:x^*=\Delta$, and is
called the {\em right complement} of $x$. The element $x^*$ is the
maximal element $s$ (with respect to $\preceq$) such that $xs$ is
simple.  A product $ab$ of simple elements $a$ and $b$ is left
weighted if and only if $a^*\wedge b=1$.  It will be convenient to
define the {\it right complement map} $\partial :\: [1,\Delta]
\rightarrow [1,\Delta]$ by $\partial(x)=x^*$:

\begin{lemma}
The map $\partial:\: [1,\Delta] \rightarrow [1,\Delta]$ is a
bijection, and $\partial^2=\tau$.
\end{lemma}

\begin{proof}
We show that $\partial$ is a bijection by defining its inverse
$\partial^{-1}:\: [1,\Delta]\rightarrow [1,\Delta]$ as
$\partial^{-1}(y)= \Delta y^{-1}$. The element $\Delta y^{-1}$,
sometimes denoted $^*y$, is called the {\em left complement} of $y$.
It is the only simple element such that $^* y \: y = \Delta$.

On the other hand, $\partial^2(x)= \partial(x^{-1}\Delta) =
(\Delta^{-1} x)\Delta = \tau(x)$, as we wanted to show.
\end{proof}

\begin{corollary}
There exists a positive integer $e$ such that $\Delta^e$ belongs to
the center of $G$. More precisely, one has
$\tau([1,\Delta])=[1,\Delta]$ and $\tau(A)=A$, where $A$ is the set
of atoms in $G$, and  $\tau^e=\mbox{id}_G$ for some positive integer
$e$, so that $\Delta^e$ is central.
\end{corollary}

\begin{proof}
Since $\partial([1,\Delta])=[1,\Delta]$, it follows that
$\tau([1,\Delta])=\partial^2([1,\Delta])=[1,\Delta]$. This also
implies that $\tau(A)=A$. Indeed, suppose that there is some atom
$a$ such that $\tau(a)$ is not an atom. Then $\tau(a)$ is a simple
element that can be decomposed into a product of two simple elements
$\tau(a)=st$. But then $\tau^{-1}(s)$ and $\tau^{-1}(t)$ are simple
elements such that $a=\tau^{-1}(\tau(a))=\tau^{-1}(s) \tau^{-1}(t)$.
A contradiction, since  $a$ is an atom. Hence $\tau(A)\subset A$.
Since $A$ is a finite set, and $\tau: G\rightarrow G$ is a
bijection, it follows that $\tau(A)=A$.

Finally, since $\tau$ induces a permutation in $A$, there exists a
positive integer $e$ such that $\tau^e$ induces the trivial
permutation on $A$. Since the atoms generate $G$, it follows that
$\tau^e$ is the trivial automorphism of $G$. That is to say,
$\Delta^e$ is central.
\end{proof}

{\bf Remark:} In the braid group $B_n$ one has $e=2$, so $\Delta^2$
is central. Furthermore, the center of $B_n$ is the cyclic group
generated by $\Delta^2$.

The right complement plays an important role when comparing the left
normal forms of $X$ and $X^{-1}$.

\begin{theorem}\label{T:left normal form of the inverse}  {\em \cite{El-M}} If
$$\
X = \Delta^p x_1\cdots x_r,
$$
in left normal form, then the left normal form of $X^{-1}$ is equal to
$$
 X^{-1}=\Delta^{-p-r}x_r'\cdots x_1',
$$
where $x_i'=\tau^{-p-i}(\partial(x_i))$ for $i=1,\ldots,r$.
\end{theorem}

{\bf Remark:} Notice that
$x_i'=\tau^{-p-i}(\partial(x_i))=\partial^{-2p-2i+1}(x_i)$, so the
left normal form of $X^{-1}$ is equal to
$$
X^{-1}= \Delta^{-p-r} \;\partial^{-2p-2r+1}(x_r)\;
\partial^{-2p-2r+3}(x_{r-1})\cdots
\partial^{-2p-1}(x_1).
$$

\begin{corollary}\label{C:inf, sup and length of the inverse}
For every $X\in G$, one has $\inf(X^{-1})=-\sup(X)$, $\;
\sup(X^{-1})=-\inf(X)$ and $\ell(X^{-1})=\ell(X)$.
\end{corollary}

See Section 9.5 of~\cite{Epstein} for a proof  that an $n$-braid of
length $m$ can be put in left normal form in running time $O(m^2n
\log n)$, with the usual Garside structure of $B_n$, and
see~\cite{BKL-1998} to find how one can compute the normal form in
time $O(m^2n)$, using the {\it dual} Garside structure of $B_n$,
usually known as Birman-Ko-Lee structure. In general, using the
normal form algorithm, the complexity of computing the left normal
form of a given element in a Garside group $G$ is $O(m^2 p)$, where
$p$ is the complexity of computing the gcd of two simple elments in
$G$. The number $p$ usually depends on the length of $\Delta$
(simple elements are smaller than $\Delta$) and on the number of
atoms in $G$, since one usually computes the gcd of two elements by
iteratively testing if there is some atom which is a common prefix.

We now explain the algorithms for solving the conjugacy decision and
search problems (CDP/CSP) in Garside groups that were given in
\cite{Garside, El-M, F-GM, Gebhardt}:  given two elements $X,Y\in
G$, determine if $X$ and $Y$ are conjugate and, if this is the case,
compute a conjugating element $Z$ such that $X^Z=Z^{-1}XZ=Y$. Each
algorithm in~\cite{Garside,El-M,F-GM, Gebhardt} is an improvement of
the previous one, but the basic idea is the same in all of them:
Given an element $X\in G$, the algorithm computes a finite subset
$I_X$ of the conjugacy class of $X$ which has the following
properties:
\begin{enumerate}
\item[(1)] For every $X\in G$, the set $I_X$ is finite, non-empty and
only depends on the conjugacy class of $X$.  In particular, two
elements $X,Y\in G$ are conjugate if and only if $I_X=I_Y$ or,
equivalently, $I_X\cap I_Y\neq \emptyset$.
\item[(2)] Given $X\in G$, a representative $\widetilde X\in I_X$ and an
element $a\in G$ such that $X^a=\widetilde X$ can be computed
effectively.
\item[(3)] Given a non-empty subset $I\subset I_X$, there is a finite
process which either proves that $I=I_X$ or produces an element
$Z\in I$ and an element $b\in G$ such that $Z^b\in I_X\backslash I$.
In particular, $I_X$ can be constructed from any representative as
the closure under this process.
\end{enumerate}

Given $X,Y\in G$, solving the CDP/CSP then involves the following
steps.
\begin{itemize}
\item[(a)] Find representatives $\widetilde X\in I_X$ and
$\widetilde  Y\in I_Y$.
\item[(b)]  Repeatedly use the process from (3), keeping track of
the conjugating elements, to compute further elements of $I_X$ until
either
\begin{enumerate}
 \item[(i)] $\widetilde Y$ is found as an element of $I_X$, proving $X$
 and $Y$ to be conjugate and providing a conjugating element, or
 \item[(ii)] the entire set $I_X$ has been constructed without
 encountering $\widetilde Y$, proving that $X$ and $Y$ are not
 conjugate.
\end{enumerate}
\end{itemize}

We now discuss, briefly, each particular algorithm in~\cite{Garside,
El-M,F-GM,Gebhardt}. In Garside's original algorithm~\cite{Garside},
the set $I_X$ is the {\em Summit Set} of $X$, denoted $SS(X)$, which
is the set of conjugates of $X$ having {\em maximal infimum}. This
was improved by Elrifai and Morton~\cite{El-M} who considered
$I_X=SSS(X)$, the {\em super summit set} of $X$, consisting of the
conjugates of $X$ having {\em minimal canonical length}. They also
show that $SSS(X)$ is the set of conjugates of $X$ having maximal
infimum and minimal supremum, at the same time.

For instance, in the braid group $B_n$ with the usual Garside
structure, one has $SS(\sigma_1)=SSS(\sigma_1)=\{\sigma_1,\ldots,
\sigma_{n-1}\}$. A small example in which $SSS(X)$ is strictly
smaller than $SS(X)$ is given by $X=\Delta \sigma_1\sigma_1\in B_3$,
for which $SSS(X)=\{\Delta\cdot\sigma_1\sigma_3\}$ and
$SS(X)=\{\Delta\cdot\sigma_1\sigma_3,
\;\Delta\cdot\sigma_1\cdot\sigma_1,\; \Delta\cdot\sigma_3\cdot
\sigma_3 \}$ (the factors in each left normal form are separated by
a dot). In general $SSS(X)$ is much smaller than $SS(X)$.

Starting by a given element $X$, one can find an element $\widetilde
X\in SSS(X)$ by a sequence of special conjugations, called {\em
cyclings} and {\em decyclings}. The conjugating elements involved in
a cycling or a decycling will play a crucial role later, so we start
by defining them.

\begin{definition} \label{D:iota_varphi}
{\rm Given $X\in G$ whose left normal form is $X=\Delta^p x_1\cdots
x_r$ ($r>0$), we define the {\em initial factor} of $X$ as
$\iota(X)=\tau^{-p}(x_1)$, and the {\em final factor} of $X$ as
$\varphi(X)=x_r$. If $r=0$ we define $\iota(\Delta^p)=1$ and
$\varphi(\Delta^p)=\Delta$.}
\end{definition}

 {\bf Remark:} Up to conjugation by $\Delta^p$, the simple element
$\iota(X)$ (resp. $\varphi(X)$) corresponds to the first (resp.
last) non-$\Delta$ factor in the left normal form of $X$. An
equivalent definition of $\iota(X)$ and $\varphi(X)$, which does not
involve the left normal form of $X$ (although it involves its
infimum and supremum), is the following.
\begin{itemize}

\item $\iota(X) =
X\Delta^{-p}\wedge\Delta$.

\item $\varphi(X) = (\Delta^{p+r-1}\wedge X)^{-1}X$.

\end{itemize}
This explains why $\iota(\Delta^p)$ and $\varphi(\Delta^p)$ are
defined in the above way.

\ms

The initial and final factors of $X$ and $X^{-1}$ are closely
related.

\begin{lemma}\label{L:initial & final factors of inverses}
For every $X\in G$ one has $\iota(X^{-1})=\partial(\varphi(X))$ and
$\varphi(X^{-1})=\partial^{-1}(\iota(X))$.
\end{lemma}

\begin{proof}
Let $\Delta^p x_1\cdots x_r $ be the left normal form of $X$,  and
suppose that $r>0$. We know that $\Delta^{-p-r} x_r'\cdots x_1'$ is
the left normal form of $X^{-1}$, where
$x_i'=\tau^{-p-i}(\partial(x_i))$. Hence one has $\iota(X^{-1})=
\tau^{p+r}(x_r')=
\tau^{p+r}(\tau^{-p-r}(\partial(x_r)))=\partial(x_r) =
\partial(\varphi(X))$.  Permuting $X$ and $X^{-1}$ in this formula
yields $\iota(X)=\partial(\varphi(X^{-1}))$, hence
$\varphi(X^{-1})=\partial^{-1}(\iota(X))$.

If $r=0$, that is if $X=\Delta^p$, then $\iota(X^{-1})= 1 =
\partial(\Delta) = \partial(\varphi(X))$, and
$\varphi(X^{-1}) = \Delta = \partial^{-1}(1) =
\partial^{-1}(\iota(X))$, so the result is also true in this case.
\end{proof}

\ms
 {\bf Remark:} The above result can be restated as follows: For every
$X\in G$, one has $\varphi(X)\iota(X^{-1}) =\Delta =
\varphi(X^{-1})\iota(X)$.

\ms

We can now define the very special conjugations called {\em
cyclings} and {\em decyclings}.

\begin{definition} {\rm Given $X\in G$, we call $\mathbf c(X)=X^{\iota(X)}$
the {\em cycling} of $X$ and we call $\mathbf
d(X)=X^{\varphi(X)^{-1}}$ the {\em decycling} of $X$. In other
words, if $\Delta^p x_1\cdots x_r$ is the left normal form of $X$
and $r>0$, then
$$
 \mathbf c(X)  = \Delta^p\: x_2\cdots x_r
 \:\tau^{-p}(x_1), \quad {\rm and} \quad  \mathbf d(X)  = x_r \:\Delta^p \:x_2\cdots
 x_{r-1}.
$$
In the case $\ell(X)=0$, we have $\mathbf c(X) = \mathbf d(X) = X$.
}
\end{definition}

Roughly speaking, for an element of positive canonical length, the
cycling of $X$ is computed by passing the first simple factor of $X$
to the end, while the decycling of $X$ is computed by passing the
last simple factor of $X$ to the front. However, the powers of
$\Delta$ are not taken into account, which is why one must use the
automorphism $\tau$. Notice that the above decompositions of
$\mathbf c(X)$ and $\mathbf d(X)$ are not, in general, left normal
forms. Hence, if one wants to perform iterated cyclings or
decyclings, one needs to compute the left normal form of the
resulting element at each iteration.

As we said above, cyclings and decyclings can be used to find an
element in $SSS(X)$, given $X$. The following result was shown for
braid groups, but the same proof is valid for every Garside group.
Choose any $X \in G$. Let $r=\ell(X)$ and let $m$ be the letter
length of $\Delta$ in the atoms of the given Garside structure.

\begin{theorem} {\rm \cite{El-M}, \cite{BKL-2001}}
Let $(G,P,\Delta)$ be a Garside structure of finite type. Choose $X
\in G$, and let $r=\ell(X)$. Let $m$ be the letter length of
$\Delta$.
\begin{enumerate}

  \item A sequence of at most $rm$ cyclings and decyclings applied to $X$
produces a representative $\widetilde X \in SSS(X)$.

  \item If $Y\in SSS(X)$ and $\alpha \in P$ is such that $Y^\alpha\in
SSS(X)$ then $Y^{\alpha\wedge \Delta} \in SSS(X)$.

\end{enumerate}
\end{theorem}

\medskip

Notice that $\alpha\wedge \Delta$ is always a simple element. Since
the set of simple elements is finite, one has the following:

\begin{corollary} {\em \cite{El-M}} Let $X\in G$ and ${\cal V}\subset SSS(X)$ be non-empty.
If ${\cal V}\neq SSS(X)$ then there exist $Y\in {\cal V}$ and a
simple element $s$ such that $Y^s\in SSS(X)\backslash {\cal V}$.
\end{corollary}

Since $SSS(X)$ is a finite set, the above corollary allows to
compute the whole $SSS(X)$. More precisely, if one knows a subset
${\cal V}\subset SSS(X)$ (at the beginning ${\cal V}=\{\widetilde
X\})$, one conjugates each element in ${\cal V}$ by {\em all} simple
elements (recall that $G$ is of finite type, that is, the set of
simple elements is finite). If one encounters a new element $Z$ with
the same canonical length as $\widetilde X$ (a new element in
$SSS(X)$), then consider ${\cal V}\cup \{Z\}$ and start again. If no
new element is found, this means that ${\cal V}=SSS(X)$, and we are
done.  One important remark is that this algorithm not only computes
the set $SSS(X)$, but it also provides conjugating elements joining
the elements in $SSS(X)$. Hence it solves both the CDP and the CSP
in Garside groups.

The computational cost of computing $SSS(X)$ depends mainly in two
ingredients: the size of $SSS(X)$ and the number of simple elements.
If we consider braid groups $B_n$ with the usual Garside structure,
for instance, all known upper bounds for the size of $SSS(X)$ are
exponential in $n$, although it is conjectured that for fixed $n$ a
polynomial bound in the canonical length of $X$
exists~\cite{Epstein}. Recall also that the number of simple
elements is $n!$, and one needs to conjugate {\em every} element in
$SSS(X)$ by {\em all} simple elements. Fortunately, this task can be
avoided thanks to the following result.

\begin{theorem}\label{T:convexity_SSS}{\rm \cite{F-GM}}
Let $X\in SSS(X)$. If $s,t\in G$ are such that $X^s\in SSS(X)$ and
$X^t\in SSS(X)$, then $X^{s\wedge t}\in SSS(X)$.
\end{theorem}

\begin{corollary}
Let $X\in G$ and $Y\in SSS(X)$. For every $u\in P$ there is a unique
$\preceq$-minimal element $\rho_Y(u)$ satisfying
$$
    u\preceq \rho_Y(u) \quad \mbox{and} \quad Y^{\rho_Y(u)}\in SSS(X).
$$
\end{corollary}

\begin{proof}
The gcd of $\{v\in P \:|\; u\preceq v,\; Y^v\in SSS(X)\}$ is the
element $\rho_Y(u)$ and has all the claimed properties.
\end{proof}

The set $\rho_Y(A)=\{\rho_Y(a)\;|\; a \mbox{ is an atom}\}$ contains
all nontrivial elements which are $\preceq$-minimal among those
conjugating $Y$ to an element in $SSS(X)$. We call the latter the
{\em minimal simple elements} for $Y$ with respect to $SSS(X)$.
Since one could have $\rho_Y(a)\prec \rho_Y(b)$ (strict) for two
distinct atoms $a$ and $b$, the set of minimal simple elements for
$Y$ is in general strictly contained in $\rho_Y(A)$.

\begin{corollary}
Let $X\in G$ and ${\cal V}\subset SSS(X)$ be non-empty. If ${\cal
V}\neq SSS(X)$ then there exist $Y\in {\cal V}$ and a minimal simple
element $\rho=\rho_Y(a)$ such that $Y^{\rho} \in SSS(X)\backslash
{\cal V}$.
\end{corollary}

Using the technique of minimal simple elements, super summit sets
can be computed as in~\cite{El-M}, but instead of conjugating each
element $Y\in SSS(X)$ by all simple elements, it suffices to
conjugate $Y$ by its minimal simple elements.  Notice that the
number of minimal simple elements for a given $Y\in SSS(X)$ is
bounded by the number of atoms. In the case of the braid group $B_n$
with the usual Garside structure, the number of atoms is $n-1$,
hence one just needs to perform $n-1$ conjugations instead of $n!$,
for each element in $SSS(X)$. Moreover, the minimal simple elements
for a given $Y\in SSS(X)$ can be computed very fast~\cite{F-GM}.

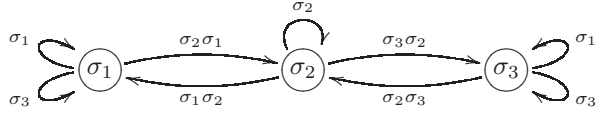
\begin{figure}[ht]
$\centerline{\xymatrix@C=2cm{
 *=<7mm,7mm>[o][F]{\sigma_1} \ar@(l,ul)[]^{\sigma_1} \ar@(l,dl)[]_{\sigma_3} \ar@/^/[r]^{\sigma_2\sigma_1} &
 *=<7mm,7mm>[o][F]{\sigma_2}  \ar@(ul,ur)[]^{\sigma_2}  \ar@/^/[l]^{\sigma_1\sigma_2} \ar@/^/[r]^{\sigma_3\sigma_2}  &
 *=<7mm,7mm>[o][F]{\sigma_3} \ar@(r,ur)[]_{\sigma_1} \ar@(r,dr)[]^{\sigma_3} \ar@/^/[l]^{\sigma_2\sigma_3}
 }}$
 \caption{Graph associated to $SSS(\sigma_1)$ in
$B_4$.}\label{figure:SSSsigma1}
\end{figure}

Notice that the algorithm just described computes not only the set
$SSS(X)$, but also the minimal simple elements that connect the
elements in $SSS(X)$ by conjugations. In other words, the algorithm
computes a directed graph whose vertices are the elements in
$SSS(X)$, and whose arrows are defined as follows: there is an arrow
labeled by $\rho$ starting at $Y$ and ending at $Z$ if $\rho$ is a
minimal simple element for $Y$ and $Y^\rho = Z$. In
Figure~\ref{figure:SSSsigma1} one can see the graph associated to
$\sigma_1\in B_4$. Notice that there are exactly 3 arrows starting
at every vertex (the number of atoms in $B_4$). In general, the
number of arrows starting at a given vertex can be smaller or equal,
but never bigger than the number of atoms.

Let us mention here a tool that will be used several times in this
paper, which is the {\em transport map} introduced
in~\cite{Gebhardt}. Let $X\in SSS(X)$ and let $\alpha$ be an element
such that $\alpha^{-1}X\alpha = Y \in SSS(X)$. We can write this by
$X\stackrel{\alpha}{\longrightarrow} Y$. We know from \cite{El-M}
that $\mathbf c(X)$ and $\mathbf c(Y)$ also belong to $SSS(X)$.
Notice that $X\stackrel{\iota(X)}{\longrightarrow} \mathbf c(X)$ and
$Y\stackrel{\iota(Y)}{\longrightarrow} \mathbf c(Y)$.
In~\cite{Gebhardt}, the {\it transport} $\alpha^{(1)}$ of $\alpha$
is defined as the element making the following diagram commutative
in the sense explained below:
$$
\begin{CD}
 X   @>\iota(X)>> \mathbf c(X) \\
@V\alpha VV   @VV\alpha^{(1)}V \\
 Y @>\iota(Y)>> \mathbf c(Y)
\end{CD}
$$
This means $\alpha^{(1)}= \iota(X)^{-1} \:\alpha \:\iota(Y)$.  The
nontrivial fact shown in~\cite{Gebhardt} is that if $\alpha$ is
simple, then $\alpha^{(1)}$ is simple, and if $\alpha$ is a minimal
simple element for $X$, then $\alpha^{(1)}$ is a minimal simple
element for $\mathbf c(X)$.

At this point, the size of the set of simple elements is no longer a
problem for the complexity of the algorithm, but there is still a
big problem to handle: The size of $SSS(X)$ is, in general, very
big. The most recent improvement, given in~\cite{Gebhardt}, is to
define a small subset of $SSS(X)$ satisfying all the good properties
described above, so that a similar algorithm can be used to compute
it. The definition of this new subset appeared after observing that
the {\em cycling} function maps $SSS(X)$ to itself. As $SSS(X)$ is
finite, iterated cycling of any representative of $SSS(X)$ must
eventually become periodic. Hence it is natural to define the
following:

\begin{definition} {\rm Given $X\in G$, define the {\em ultra summit set}
of $X$, $\; USS(X)$, to be the set of elements $Y\in SSS(X)$ such
that $\mathbf c^m(Y)=Y$, for some $m > 0$.}
\end{definition}

The ultra summit set $USS(X)$ thus consists of a (finite) set of
disjoint, closed orbits under cycling. For instance, in the braid
group $B_n$ one has
$USS(\sigma_1)=SSS(\sigma_1)=SS(\sigma_1)=\{\sigma_1,\ldots,\sigma_{n-1}\}$,
and each element corresponds to an orbit under cycling, since
$\mathbf c(\sigma_i)=\sigma_i$ for $i=1,\ldots,n-1$. A less trivial
example is given by the element
$$
X=\sigma_1\sigma_3\sigma_2\sigma_1 \cdot \sigma_1\sigma_2 \cdot
\sigma_2\sigma_1\sigma_3 \in B_4.
$$
In this example $USS(X)$ has 6 elements, while $SSS(X)$ has 22
elements. More precisely, the ultra summit set of $X$ consists of 2
closed orbits under cycling, $USS(X)=O_1\cup O_2$, each one
containing 3 elements:
$$
 O_1=\{
 \sigma_1\sigma_3\sigma_2\sigma_1 \cdot \sigma_1\sigma_2 \cdot \sigma_2\sigma_1\sigma_3,
 \quad  \sigma_1\sigma_2 \cdot \sigma_2\sigma_1\sigma_3 \cdot \sigma_1\sigma_3\sigma_2\sigma_1,
 \quad \sigma_2\sigma_1\sigma_3 \cdot \sigma_1\sigma_3\sigma_2\sigma_1 \cdot \sigma_1\sigma_2
 \},
$$
$$
 O_2=\{
  \sigma_3\sigma_1\sigma_2\sigma_3 \cdot \sigma_3\sigma_2 \cdot \sigma_2\sigma_3\sigma_1,
 \quad  \sigma_3\sigma_2 \cdot \sigma_2\sigma_3\sigma_1 \cdot \sigma_3\sigma_1\sigma_2\sigma_3,
 \quad \sigma_2\sigma_3\sigma_1 \cdot \sigma_3\sigma_1\sigma_2\sigma_3 \cdot \sigma_3\sigma_2
 \}.
$$
Notice that $O_2=\tau(O_1)$. Notice also that the cycling of every
element in $USS(X)$ gives another element which is already in left
normal form, hence iterated cyclings corresponds to cyclic
permutations of the factors in the left normal form. We will say
that elements satisfying this property are {\it rigid}. The precise
definition will be given in $\S$\ref{subsection:summary of results}.
We remark that the size of the ultra summit set of a generic braid
of canonical length $l$ is either $l$ or $2l$~\cite{Gebhardt}. This
means that, in the generic case, ultra summit sets consist of one or
two orbits (depending on whether $\tau(O_1)=O_1$ or not), containing
rigid braids.

The algorithm given in~\cite{Gebhardt} to solve the CDP/CSP in
Garside groups (of finite type) is analogous to the previous ones,
but this time one needs to compute $USS(X)$ instead of $SSS(X)$. In
order to do this, the following results, which are analogous to
those given for super summit sets, are used.

\begin{theorem}\label{T:convexity_USS} {\rm \cite{Gebhardt}}
Let $X\in USS(X)$. If $s,t\in G$ are such that $X^s\in USS(X)$ and
$X^t\in USS(X)$, then $X^{s\wedge t}\in USS(X)$.
\end{theorem}

\begin{corollary} {\rm \cite{Gebhardt}}
Let $X\in G$ and $Y\in USS(X)$. For every $u\in P$ there is a unique
$\preceq$-minimal element $c_Y(u)$ satisfying
$$
    u\preceq c_Y(u) \quad \mbox{and} \quad Y^{c_Y(u)}\in USS(X).
$$
\end{corollary}

\begin{definition} {\rm Given $X\in G$ and $Y\in USS(X)$, we say that a
simple element $s\neq 1$ is a {\em minimal simple element} for $Y$
with respect to $USS(X)$ if $Y^s=s^{-1}Ys \in USS(X)$, and no proper
prefix of $s$ satisfies this property.}
\end{definition}

Notice that the set of minimal simple elements for $Y$ with respect
to $USS(X)$ is contained in $c_Y(A) = \{c_Y(a)\; | \; a \mbox{ is an
atom}\}$, hence the number of minimal simple elements for $Y$ is
bounded by the number of atoms. For the rest of the paper, all
minimal simple elements will be considered with respect to ultra
summit sets (and not super summit sets).

\begin{corollary}{\rm \cite{Gebhardt}}
Let $X\in G$ and ${\cal V}\subset USS(X)$ be non-empty. If ${\cal
V}\neq USS(X)$ then there exist $Y\in {\cal V}$ and an atom $a$ such
that $c_Y(a)$ is a minimal simple element for $Y$, and $Y^{c_Y(a)}
\in USS(X)\backslash {\cal V}$.
\end{corollary}

In~\cite{Gebhardt} it is shown how to compute the minimal simple
elements corresponding to a given $Y\in USS(X)$, hence one can
compute the whole $USS(X)$ starting by a single element $\widetilde
X\in USS(X)$.

As above, the algorithm in~\cite{Gebhardt} not only computes
$USS(X)$ but also a graph which determines the conjugating elements.
This graph is defined as follows.

\begin{definition} {\rm Given $X\in G$, the directed graph $\Gamma_X$ is
defined by the following data:
\begin{enumerate}

 \item The set of vertices is $USS(X)$.

 \item For every $Y\in USS(X)$ and every minimal simple element $s$
 for $Y$ with respect to $USS(X)$, there is an arrow labeled by $s$
 going from $Y$ to $Y^s$.

\end{enumerate}
}
\end{definition}

We remark that one obtains an element $\widetilde X\in USS(X)$ by
iterated application of cycling to an element in $SSS(X)$, which we
know how to compute using cyclings and decyclings. The number of
times one needs to apply cycling, in order to go from an element in
$SSS(X)$ to an element in $USS(X)$ is not known in general.
Nevertheless, the theoretical complexity of the algorithm
in~\cite{Gebhardt} is not worse than the one of the algorithm
in~\cite{F-GM}, and is substantially better in practice, at least
for braid groups.

In fact, it follows from the work in~\cite{Garside,
El-M,BKL-2001,F-GM,Gebhardt} discussed above, that the complexity of
CDP/CSP for two elements $X$, $Y$ in a Garside group ($||X||\geq
||Y||)$ is $O(|USS(X)|\: p+ q)$, where $p$ is a polynomial in
$||X||$ and the number of atoms, and $q$ is related to the number of
times one must apply cycling to an element in $SSS(X)$ to transform
it into an element in $USS(X)$. We believe that the second term $q$
is negligible compared to $|USS(X)|\: p$, so our main interest is in
trying to bound the size of the ultra summit set of an element in a
Garside group.

In the particular case of braid groups, the size and structure of an
ultra summit set happen to depend heavily on the geometrical
properties of the braid, more precisely, on its Nielsen-Thurston
type. This is explained next.

\subsection{The Thurston-Nielsen trichotemy in the braid groups}
\label{subsection:the Thurston-Nielsen trichotemy in the braid
groups}

The braid group $B_n$ is isomorphic to the mapping class group
$\pi_0({\rm Diff}_+(D^2_n))$ of the disc with $n$ points removed.
Admissible diffeomorphisms preserve orientation, fix $\partial D^2$
pointwise and fix the n punctures or distinguished points setwise.
Admissible isotopies fix both pointwise.   As a mapping class group,
$B_n$ has structure which, at this time, has not been fully related
to its Garside structure, although some interesting relation between
the two structures can be found in~\cite{D-W}. We will use the
geometric structure in $\S$\ref{subsection:Consequences for
pseudo-Anosov braids} and also in~\cite{BGGM-II,BGGM-III}, so we
describe what we need here. The structure that we describe had its
origins in 3 very long papers of J. Nielsen \cite{Nielsen}, written
in the 1930's, but the grand sweep of the theory was not recognized
until much later, in the work of W. Thurston \cite{Thurston1988}. We
refer to it as the {\it Thurston-Nielsen trichotemy}.  There are
many ways to describe it. We choose one which is based upon the
action of $B_n$ on isotopy classes of simple closed curves (scc) on
$D^2_n$. The scc considered in $D^2_n$ are non-degenerate, which
means that they bound neither a single puncture nor all punctures
(otherwise they could be collapsed to a puncture or isotoped to the
boundary).

\begin{theorem} {\rm \cite{Thurston1988, Nielsen}}
Let $X\in B_n$. Then, after a suitable isotopy, $X$ belongs to
exactly one of the following pairwise disjoint classes: \be

\item $X$ is `periodic'. That is, some power of $X$ is a power of a Dehn twist on
$\partial D^2_n$ (this Dehn twist is precisely $\Delta^2$, with the
usual Garside structure).
\item $X$ is `pseudo-Anosov' or PA. That is, neither $X$ nor any power of $X$ fixes the isotopy class of any scc on $D^2_n$.  This case is the generic case.
\item $X$ is `reducible'. That is, there exists a family of scc on $D^2_n$ whose isotopy class is  fixed  by $X$, so that some power $X^m$ of $X$ fixes the isotopy class of each simple closed curve in the family.  Moreover, if the disc $D^2_n$ is split open along suitable representatives of the fixed curves, then the restriction of $X^m$ to the closure of each component of the split-open disc is either periodic or PA.
\ee
\end{theorem}

We note that there is a working algorithm, given in
\cite{BGN1993,BGN1995}, to determine whether a given braid is
periodic or reducible, and we used it  in basic ways when we
computed the millions of examples that suggested the different
structures of ultra summit sets, depending on the geometric type.
Note that if one can recognize whether a braid is periodic or
reducible, then if it is neither it must be PA.

\subsection{A project to solve the conjugacy problems in braid groups and a summary of our results}
\label{subsection:summary of results}

Making use of almost all of the ideas that we have just described,
we have developed a strategy for attacking the problem of the
complexity of the conjugacy decision and search problems (CDP/CSP)
in the braid groups. It uses the structure of centralizers of PA
braids, and the uniqueness of their roots~\cite{GM}, some particular
properties of periodic braids, and also the geometric decomposition
of a reducible braid along its invariant curves. Hence, our strategy
does not apply to an arbitrary Garside group, although many of the
results that we show (all results in~\cite{BGGM-II} and all results
in this paper, except Theorem~\ref{T:bound for pseudo-Anosov} and
those in $\S$\ref{subsection:Consequences for pseudo-Anosov braids})
are stated and hold in the general framework of Garside groups. The
results in~\cite{BGGM-III} and the remaining parts of our project
are conceived for braid groups, although we believe that they will
probably be generalized to other Garside groups, at least to
spherical type Artin-Tits groups.

As was noted in the previous sections, two elements $X,Y$ in a
Garside group are conjugate if and only if one element in $USS(X)$
is also in $USS(Y)$. This means that we must compute all of $USS(Y)$
in order to be able to test conjugacy. Thus we will need to
understand  the structure and size of the ultra summit set.
Unfortunately, however, $USS(X)$ can be quite complicated, partly
because cycling is not, in general, a cyclic permutation of the
factors in a left normal form, but also because it is not clear how
the distinct orbits in $USS(X)$ are related. The former problem is
avoided if $USS(X)$ is made of {\it rigid} elements:

In Section~\ref{section:rigidity} we will introduce and study {\it
rigid elements}.  Let $X=\Delta^px_1x_2\cdots x_r$ be in left normal
form. Assume $r>0$. Then $X$ is rigid if $\Delta^px_1x_2\cdots
x_r\tau^{-p}(x_1)$ is in normal form as written. We were lead to
study rigid elements when we realized, long ago, that it was often
very difficult to predict and understand the changes in normal form
of braids after cycling. If $X$ is rigid, the left normal form of
$\mathbf c(X)$ is precisely $\Delta^p x_2 \cdots x_r
\tau^{-p}(x_1)$, so cycling is simpler than in the general case, and
the combinatorics in $USS(X)$ are easier to understand.

In this paper we will see that obtaining a polynomial solution to
the CDP/CSP for certain elements in a Garside group, reduces to
obtaining such a solution for rigid elements. In the case of braid
groups, this happens for pseudo-Anosov (PA) braids. Since the
property of being PA is generic in $B_n$, this is an important step
in the case of braids.

Assuming that $X,Y \in B_n$, we consider the three cases separately:
$X,Y$ are PA, periodic  or reducible. We break our approach to $B_n$
into  the following 6 steps:

\be

\item [I] {\bf Determining if a braid is periodic, reducible or PA.}
We remark that it is very fast to decide whether a given braid is
periodic~\cite{GM}, so the main problem is to determine if a braid
is reducible, and find the reducing curves. This question is solved
in~\cite{BGN1993,BGN1995}, but the proposed algorithm computes
$SSS(X)$. In fact, one can replace $SSS(X)$ by $USS(X)$, but having
to compute $USS(X)$ means that the algorithm is not polynomial, in
general (In~\cite{BGGM-III} there are examples of USS's in $B_n$
whose size is exponential in $n$). This yields the following.

{\bf Open question 1:} Is there an algorithm to determine if a braid
in $B_n$ is reducible and to find its reducing curves, which is
polynomial in $n$ and $||X||$?

This problem was first studied in~\cite{BGN1993}. Some work in this
direction can be found in~\cite{Lee_reducible}.

\item  [II] {\bf PA braids: passing to powers.}
In Sections~\ref{section:cyclings and powers} and
\ref{section:rigidity} of this paper we show that if $X$ and $Y$ are
PA, there is some small power $m$ such that $USS(X^m)$ is made of
rigid braids, and it suffices to solve the CDP/CSP for $X^m$ and
$Y^m$.


In this regard we make two remarks: The first is that, in view of
the results in \cite{GM}, for every nonzero integer $m$, the braids
$X$ and $Y$ are conjugate if and only if $X^m$ and $Y^m$ are
conjugate. Furthermore, $PA$ braids have unique roots. Hence if $X$
and $Y$ are $PA$, and $Z$ conjugates $X^m$ to $Y^m$, then $Z$
conjugates $X$ to $Y$. Therefore nothing is lost in passing to
powers.

Our second remark is that we prove the non-emptiness of the {\it
stable ultra summit set} $SU(X)$ in a Garside group (compare with
\cite{Lee-Lee}, where the stable super summit set is introduced).
That is, for every $X\in G$ define $SU(X) = \{Y\in USS(X) \ |\  Y^k
\in USS(Y^k) \ {\rm for \ \ all} \ \ k\in \mathbb Z\}$.
Proposition~\ref{P:SU is not emptyset} of this paper proves that
$SU(X) \not= \emptyset$. However, we will not need to work in
$SU(X)$, it will suffice to control a bounded number of powers of
$X$, and we learn how to do that.

\item  [III] {\bf Understanding the USS graph.}
In \cite{BGGM-II} we uncover and study the structure of the ultra
summit set of an element in a Garside group. More precisely, we show
that the conjugations corresponding to minimal simple elements (the
arrows in the USS graph) are a very special kind of conjugation that
we call {\it partial cyclings}. This work is not restricted to
braids.  At the end of \cite{BGGM-II} we specialize our work to the
cases: (a) $X$ is a rigid element, and (b) $X$ is a periodic element
(in a Garside group, that is, a root of some power of $\Delta$).
This is a first step towards the solution of the following.

\item  [IV] {\bf Finding a polynomial bound for the size of
 USS(X), when X is rigid.}  At this writing this work is
incomplete. We have computed many many examples, using random
searches, and on the basis of the evidence found that in the generic
case $USS(X)$ has either 2 orbits, where one is the conjugate of the
other by $\Delta$, or 1 orbit which is conjugate to itself by
$\Delta$. However, there are exceptional cases where $USS(X)$ has
unexpected size. There is no indication whatsoever of uncontrolled
growth. Indeed, the combinatorial conditions that are uncovered in
\cite{BGGM-II} are so restrictive that exponential growth seems very
unlikely.  But since we do not have an affirmative answer, we state
the following.

{\bf Open question 2:} If $X$ is a rigid element in $B_n$, is the
size of $USS(X)$ bounded above by some polynomial in $n$ and
$||X||$?

We remark that, in $\S$\ref{The ultra summit set of a rigid element
is made of rigid elements} of this paper, we show that if $X\in G$
is a rigid element of canonical length greater than 1, then $USS(X)$
consists of rigid elements.

Finally, solving Open question 2 affirmatively would imply that the
algorithm in~\cite{Gebhardt} applied to rigid braids is polynomial
in $n$ and $||X||$, provided that the following is also true, at
least for conjugates of a rigid braid:

{\bf Open question 3:} Given $X\in B_n$ and $Y\in SSS(X)$, let $m$
be such that $\mathbf c^m (Y)\in USS(X)$. Is $m$ bounded above by a
polynomial in $n$ and $||X||$?

\item[V] {\bf Periodic braids.}
In \cite{BGGM-III} we settle the CSP for periodic braids in $B_n$,
in polynomial time with respect to $n$ and $||X||$. We remark that
the CDP for periodic braids was already known to be
polynomial~\cite{BGN1993,BGN1995}, but the usual algorithm to solve
the CSP is not polynomial in this case, so in~\cite{BGGM-III} we
find a new specific algorithm for periodic braids.

\item  [VI] {\bf Reducible braids.} Suppose that Open question 1, 2 and 3 above are solved. Note
that reducible braids are braids that are made up of braided tubes,
each containing braided tubes and so forth until one reaches an
irreducible braid, which is then either periodic or PA. Once that
reducing curves are known, and one knows how to solve the CDP/CSP
for irreducible braids, one can use techniques from~\cite{GM} to
solve the CDP/CSP for reducible ones in polynomial time.

 \ee

As a conclusion, the work in this paper and
in~\cite{BGGM-II,BGGM-III}, together with an affirmative answer to
Open questions 1, 2 and 3 above, would yield a polynomial algorithm
to solve the CDP/CSP in braid groups. Due to our increasing
understanding of the structure of ultra summit sets, we believe that
this final goal is within reach.

In this paper we will solve problem II above. Most of our results
(except those in $\S$\ref{subsection:Consequences for pseudo-Anosov
braids} and Theorem~\ref{T:bound for pseudo-Anosov}) hold in all
Garside groups.

In $\S$\ref{section:cyclings and powers} we determine the
relationship between the $m$-times iterated cycling $c^m(X)$ of
$X\in G$ and the $m^{th}$ power $X^m$ of $X$.  The main result is
Theorem~\ref{T:normal form of $C_m$ and $X^m$} of
$\S$\ref{subsection:Interplay between C_m and R_m}.  In
$\S$\ref{subsection:the stable ultra summit set} we introduce the
stable ultra summit set $SU(X)$ of $X\in G$ (cf \cite{Lee-Lee}) and
give a short proof that it is non-empty.  While we realized, after
we had completed the work in this paper, that we did not really need
$SU(X)$ in our work,  we include it for completeness, and because it
may be useful for others.

In $\S$\ref{section:rigidity} we study rigid elements in Garside
groups and prove some surprising results about them. In particular,
in Theorem~\ref{T:X rigid implies every Y in USS(X) rigid}  we prove
that if $X$ is rigid and $\ell(X)>1$ then every element in $USS(X)$
is also rigid. Theorem~\ref{T:characterization of pre-rigid}
characterizes exactly which elements in a Garside group have rigid
powers.  Using it, we prove in Theorem \ref{T:pA and USS implies
pre-rigid} that if $X$ is a pseudo-Anosov braid, then there exists
an $m$ such that $X^m$ is rigid.   In Theorem \ref{T:bound on
rigidity} we solve the problem that is described in II above in this
section, obtaining a polynomial bound for the power $m$.

\section{Cyclings and powers.}
\label{section:cyclings and powers}

Recall the definition of a {\em rigid} element in a Garside group
$G$ (Definition~\ref{D:rigid element}). Our goal in this paper is to
understand the conditions under which an element $X\in USS(X)\subset
G$ which is not necessarily rigid has a small power $X^m$ which is
rigid. This will be done by investigating the relationship between
iterated cyclings and iterated powers of $X$. However, the
connection between cycling and normal forms of powers is fairly
subtle. The problems that we will encounter and solve will be easier
to understand after we study an example.  They will probably have
been encountered by others who have worked with left normal forms
(see \cite{Adjan,El-M,Epstein}) in the braid group $B_n$, and
struggled to understand how they change after cycling.

{\bf An Example: } Let $X = \Delta^p x_1\cdots x_r \in USS(X)
\subset G, r > 0$. Since $\inf(X)=p$, $\ell(X) = r$, it is immediate
that $\inf(X^m) \geq mp$ and that $\ell(X^m) \leq mr$.  We can think
of the terms $\Delta^{mp}$ in the normal form of $X^m$ as the {\it
expected $\Delta's$}, and any additional ones as {\it unexpected
$\Delta's$}.  Similarly, $\sup(X) = p+r$, so that  $mp + mr$ is the
{\it expected supremum} of  $X^m$. There will be  an {\it unexpected
decrease} in $\sup(X^m)$ if and only if the actual value  of
$\sup(X^m)$ is less than $mp+mr$.   These two issues are closely
related, because by Corollary~\ref{C:inf, sup and length of the
inverse} $\sup(X) = -\inf(X^{-1})$ and $\inf(X) = -\sup(X^{-1})$, so
that if we arrive at an understanding of unexpected increases in the
infimum,  we will also have arrived at an understanding of
unexpected decreases in the supremum. Unfortunately, however, the
normal form of $X^m$ is not easily related to the normal form of
$X$, as is illustrated by the following example, taken from the
5-string braid group $B_5$.

Let $X=1 2 1 3 2 1 4 3 1 4 3\in B_n$, where the letter $i$ means the
elementary braid $\sigma_i$.  In this example $\inf(X) = 0$.  A
calculation shows that $X$ is in its ultra summit set, and there are
2 cycling orbits in $USS(X)$, each with 4 elements, with the second
being the conjugate of the first by $\Delta$.  In this simple case
$\ell(X) = 2$, that is, there are 2 simple factors in the left
normal form for $X = C_1\cdot R_1$, where dots are used to separate
the simple words in the left normal form.   Here is the first orbit:

$X = 1 2 1 3 2 1 4 3 \cdot 1 4 3 = C_1\cdot R_1$ \\
$c(X) = 1 2 1 3 2 4 3 2 1 \cdot 1 4 = C_2\cdot R_2$\\
$c^2(X) =  1 2 1 3 2 4 3 2 \cdot 2 1 4 = C_3\cdot R_3$\\
$c^3(X) =  1 2 1 3 4 3 \cdot 1 2 3 2 4 = C_4\cdot R_4$\\
$c^4(X) = X$

What about powers of $X$?  Calculating left normal forms, we find
that:

$X = 1 2 1 3 2 1 4 3 \cdot  1 4 3$\\
$X^2 =   \Delta \cdot  2 3 2 4 3 2 1 \cdot  1 4 \cdot  1 4 3$ \\
$X^3 =  \Delta^2 \cdot  1 2 3 2 4 \cdot  2 1 4 \cdot  1 4 \cdot  1 4 3$ \\
$X^4 =  \Delta^2 \cdot  1 2 1 3 2 1 4 3 \cdot  1 4 3 \cdot  1 2 3 2 4 \cdot  2 1 4 \cdot  1 4 \cdot  1 4 3$\\
$X^5 = \Delta^3 \cdot  2 3 2 4 3 2 1 \cdot  1 4 \cdot 1 4 3 \cdot 1 2 3 2 4 \cdot 2 1 4 \cdot 1 4 \cdot 1 4 3$\\

Since $\inf(X) = 0$,  the powers $X^2, X^3, X^4, X^5$ have 1,2,2,3
unexpected $\Delta's$.

A hint at how the normal forms of $X$ and $X^m$ might be related
comes from a more careful inspection of this example (and many many
other examples like it).  The initial factors of the elements in the
orbit of $X$ are   $C_1 = 1 2 1 3 2 1 4 3, \ \ C_2 =  1 2 1 3 2 4 3
2 1, \ \ C_3 = 1 2 1 3 2 4 3 2, \ \ C_4 = 1 2 1 3 4 3$.  These are
the `conjugating factors' that are used when we cycle, that is, if $
{\bf C}_m = C_1C_2\cdots C_m,$ then
$c^m(X) = X^{{\bf C}_m}. $  Here are the left normal forms for ${\bf C}_1, {\bf C}_2,\dots,{\bf C}_5$:\\
${\bf C}_1 = 1 2 1 3 2 1 4 3  = $ first simple factor in $X$\\
${\bf C}_2 = 1 2 1 3 2 1 4 3 1 2 1 3 2 4 3 2 1  = \Delta \cdot  2 3 2 4 3 2 1 = $ product of first 2 simple factors in $X^2$\\
$ {\bf C}_3 = 1 2 1 3 2 1 4 3 1 2 1 3 2 4 3 2 1 1 2 1 3 2 4 3 2  = \Delta^2 \cdot  1 2 3 2 4   =$ product of  first 3 simple factors in $X^3$\\
$ {\bf C}_4 =  \Delta^2 \cdot 1 2 1 3 2 1 4 3 \cdot  1 4 3  = $ product of  first 4 simple factors in $X^4$\\
${\bf C}_5 = \Delta^3 \cdot 2 3 2 4 3 2 1 \cdot  1 4 =$ product of first 5 simple factors in $X^5$\\

One of the main results in this paper states that, for every $X \in
USS(X)$, the product of the first $m$ factors in the left normal
form of $X^m \Delta^{-mp}$, where we include powers of $\Delta$ in
the count,  is precisely $\mathbf C_m$, the product of the
conjugating elements involved in the first $m$ cyclings of $X$. This
will allow us to determine which elements admit a rigid power and,
under some hypothesis, we find an upper bound for the smallest power
which is rigid. In the particular case of braid groups, these
results apply to pseudo-Anosov braids, since we will show in
$\S$\ref{subsection:Consequences for pseudo-Anosov braids} that
every pseudo-Anosov braid in its ultra summit set has a rigid power.

\subsection{Decomposition of powers of $X$.}

In this section we will decompose $X^m$ as a product of two
elements, each of which is determined by the iterated cyclings of
$X\in G$. Assume from now on that $X\in SSS(X)$ and $\ell(X)\geq 1$.
We will develop some basic properties of this decomposition for
elements of a Garside group.

We first need some notation. If the left normal form of $X$ is
$\Delta^p x_1\cdots x_r$, recall that $\iota(X)=\tau^{-p}(x_1)$ is
the {\em initial factor} of $X$, and that the cycling of $X$ is
defined by $\mathbf c(X)= X^{\iota(X)}$. If we apply iterated
cyclings to $X$, the conjugating elements will be denoted by $C_1,
C_2, \ldots$ That is, $C_i=\iota(\mathbf c^{i-1}(X))$ for $i\geq 1$.
Hence one has $\mathbf c^m(X) = X^{C_1\cdots C_m}$. The letter $C$
in the symbol $C_i$ comes from {\em conjugating element}, since one
conjugates $\mathbf c^{i-1}(X)$ by $C_i$ to obtain  $\mathbf
c^{i}(X)$.

The element $X$ can be decomposed as follows: $X= \Delta^p x_1\cdots
x_r = C_1 \Delta^p x_2\cdots x_r$. We denote $R_1= x_2\cdots x_r$,
so $X=C_1 \Delta^p R_1$.  For the iterated cyclings of $X$, we
denote $R_i$ in a similar way, that is, the element satisfying
$\mathbf c^{i-1}(X) = C_i \Delta^p R_i$. The letter $R$ comes from
{\em remainder}. Notice that every $C_i$ is a simple element, while
$R_i$ is simple only if $\ell(X)\leq 2$, and it is trivial if
$\ell(X)=1$. The important fact about these elements relies on how
they behave when they are multiplied in the right way.

\begin{definition} {\rm Let $X\in SSS(X)$ with $\inf(X)=p$ and $\ell(X)\geq
1$. For $i\geq 1$, let $C_i$ and $R_i$ be the elements defined
above. Then, for every $m\geq 1$, we define:
\begin{eqnarray}
\nonumber
\bullet \ \ \ \ {\bf C}_m & = & C_1 \cdots C_{m}\\
\nonumber \bullet  \ \ \ \  {\bf R}_m & = & \tau^{-p}(R_m)
\tau^{-2p}(R_{m-1})\cdots \tau^{-mp}(R_1).
\end{eqnarray}
}
\end{definition}

Notice that
$$
   \mathbf R_m \Delta^{pm} = (\Delta^p R_m) (\Delta^p R_{m-1})\cdots
   (\Delta^p R_1).
$$

Since later we will deal not only with $X$, but with successive
cyclings of $X$, we want to define the corresponding elements above,
for $\mathbf c^k (X)$. Hence we define $\mathbf C_{[k,m]}$ and
$\mathbf R_{[k,m]}$ to be the elements $\mathbf C_m$ and $\mathbf
R_m$ above, but defined with respect to $\mathbf c^k (X)$. This
yields the analogous definition with the indices shifted by $k$:

\begin{definition}{\rm Let $X\in SSS(X)$ with $\inf(X)=p$ and $\ell(X)\geq
1$. For $i\geq 1$, let $C_i$ and $R_i$ be the elements defined
above. Then, for every $m\geq 1$ and $k\geq 0$, we define:
\begin{eqnarray}
\nonumber
\bullet \ \ \ \ {\bf C}_{[k,m]} & = & C_{k+1} \cdots C_{k+m}\\
\nonumber \bullet  \ \ \ \  {\bf R}_{[k,m]} & = & \tau^{-p}(R_{k+m})
\tau^{-2p}(R_{{k+m-1}})\cdots \tau^{-mp}(R_{k+1}).
\end{eqnarray}
}
\end{definition}

Notice that
$$
   \mathbf R_{[k,m]} \Delta^{pm} = (\Delta^p R_{k+m}) (\Delta^p R_{k+m-1})\cdots
   (\Delta^p R_{k+1}).
$$

Clearly, ${\bf C}_m={\bf C}_{[0,m]}$ and ${\bf R}_m={\bf
R}_{[0,m]}$.

In the particular case in which $X\in USS(X)$, that is, $X\in
SSS(X)$ and $\mathbf c^t(X)=X$ for some positive integer $t$, we can
extend the above definition to negative values of $k$, as follows.
We know that $\mathbf c^m(X)\in USS(X)$ for every $m\geq 0$. If we
denote by $\mathcal O(X)$ the orbit of $X$ under cycling, we can
define $\mathbf c^{-m}(X)$ to be the element $Y\in \mathcal O(X)$
such that $\mathbf c^m(Y)=X$ (although cycling is not injective in
the whole $G$, it is a bijection in $USS(X)$, so we hope this
notation will not cause confusion).

Recall that we defined $C_i=\iota(\mathbf c^{i-1}(X))$, and $R_i$ in
such a way that $\mathbf c^{i-1}(X)= C_i \Delta^p R_i$, for every
$i\geq 1$. The same definitions can now be given for every $i\in
\mathbb Z$, as we have definitions for the {\em negative} cyclings
of $X$. Since $\mathcal O(X)$ is a finite set, the sequences
$\{C_i\}_{i\in \mathbb Z}$ and $\{R_i\}_{i\in \mathbb Z}$ are
periodic. Therefore, we have definitions for the elements $C_i$,
$R_i$, and also $\mathbf C_{[k,m]}$ and $\mathbf R_{[k,m]}$ for
every $i,k\in \mathbb Z$ and every $m\geq 1$.

Let us show a result that will be useful later.

\begin{lemma}\label{L: RC=CR}
Let $X\in SSS(X)$ with $\inf(X)=p$ and $\ell(X)\geq 1$. For every
$m,k\geq 1$ one has
$$
  \mathbf  R_{[k-1,m]} \Delta^{pm} \; C_k  = C_{k+m} \;\mathbf R_{[k,m]} \Delta^{pm}.
$$
In other words,
$$
(\Delta^p R_{k+m-1}) \cdots (\Delta^p R_k) C_k = C_{k+m} (\Delta^p
R_{k+m}) \cdots (\Delta^p R_{k+1}).
$$
Moreover, if $X\in USS(X)$, the same equality holds for every $k\in
\mathbb Z$.
\end{lemma}

\begin{proof}
We first show the result for $X\in SSS(X)$ and $k\geq 1$. If $m=1$,
the result is true since $(\Delta^p R_k) C_k = \mathbf c^k (X) =
C_{k+1} (\Delta^p R_{k+1})$ by definition. Suppose the result true
for $m-1$. Then one has
$$
(\Delta^p R_{k+m-1}) (\Delta^p R_{k+m-2}) \cdots (\Delta^p R_k) C_k
= (\Delta^p R_{k+m-1}) C_{k+m-1} (\Delta^p R_{k+m-1}) \cdots
(\Delta^p R_{k+1})
$$
$$
= \mathbf c^{k+m-1}(X) (\Delta^p R_{k+m-1}) \cdots (\Delta^p
R_{k+1}) =  C_{k+m} (\Delta^p R_{k+m}) (\Delta^p R_{k+m-1}) \cdots
(\Delta^p R_{k+1}),
$$
so the result is also true for $m$ and we are done.

If $X\in USS(X)$ and $k<0$, the same proof is valid.
\end{proof}

We will now see how the element $X^m$ can be decomposed in terms of
$\mathbf C_m$ and $\mathbf R_m$, together with some properties
concerning the normal form of these two factors.

\begin{lemma}\label{L:decomposition}
Let $X \in SSS(X)$, with $\ell(X)\geq 1$. Let $C_i$, $R_i$, $\mathbf
C_m$ and $\mathbf R_m$ be the elements defined above. Then:

1. The $m$-th power of $X$ has the decomposition
$$
X^m = {\bf C}_m{\bf R}_m\Delta^{mp}.
$$

2. In this decomposition, $\inf({\bf R}_m)=0$ and $\iota({\bf
R}_m)\preceq C_{m+1}$, for every $m\geq 1$.

3. In general $\inf({\mathbf C_m}) \geq 0$.  If $\ell(X)>1$, one has
$\sup({\bf C}_m)=m$, $\ell({\bf C}_m)>0$ and $\varphi({\bf
C}_m)\succeq \varphi(\mathbf c^{m}(X))$, for every $m\geq 1$.
\end{lemma}

{\bf Remark:} The left normal  form of ${\bf C}_m$ is not so easy to
understand, as we saw in the example that was given at the beginning
of Section~\ref{section:cyclings and powers}. Uncovering it, and
relating it to the left normal form of $X^m$, will be a major part
of our investigations.

{\it Proof:}

1. For $m=1$ the result is clear, since $X= C_1 \Delta^p R_1 =
\mathbf C_1 \mathbf R_1 \Delta^p$ by definition. Now suppose that
$$
X^{m-1}= \mathbf C_{m-1} \mathbf R_{m-1} \Delta^{(m-1)p}.
$$
Then one has
$$
  X^m = X^{m-1} X =  \mathbf C_{m-1} \mathbf R_{m-1} \Delta^{(m-1)p}
  \; C_1 (\Delta^p R_1).
$$
By Lemma~\ref{L: RC=CR} with $k=1$, it then follows that
$$
 X^m = \mathbf C_{m-1} C_m \mathbf R_{[1,m-1]} \Delta^{(m-1)p}
  (\Delta^p R_1)  = \mathbf C_{m} \mathbf R_{[1,m-1]} \Delta^{mp}
   R_1
$$
$$
 = \mathbf C_{m} \mathbf R_{[1,m-1]} \tau^{-mp}(R_1) \Delta^{mp}
 = \mathbf C_m \mathbf R_m \Delta^{mp},
$$
so the result is true for every $m\geq 1$.

\medskip

2. To prove that $\inf({\bf R}_m) = 0$ and that $\iota({\bf
R}_m)\preceq C_{m+1}$, for every $m\geq 1$, we notice that ${\bf
R}_m$ is positive by definition. Hence both statements will follow
if we can show that $\Delta \wedge {\bf R}_m \preceq C_{m+1}$.

If $m=1$ one has ${\bf R}_1=\tau^{-p}(x_2 \cdots x_r)$, hence
$\Delta \wedge {\bf R}_1 = \tau^{-p}(x_2)$. We also have $C_2=
\iota(\mathbf c(X)) = \iota(\Delta^p x_2\cdots x_r \tau^{-p}(x_1))$.
Since $X\in SSS(X)$, we have $\inf(x_2\cdots x_r \tau^{-p}(x_1))=0$,
so the first factor in its left normal form is equal to $x_2 s$ for
some simple element $s$. Hence $C_2 = \iota(\Delta^p x_2 \cdots x_r
\tau^{-p}(x_1)) = \tau^{-p}(x_2 s)$. Therefore $\Delta \wedge {\bf
R}_1 = \tau^{-p}(x_2) \preceq \tau^{-p}(x_2 s) = C_2$, and the
result is true for $m=1$.

Suppose that $\Delta \wedge {\bf R}_{m-1} \preceq C_{m}$ for some
$m$. By definition, ${\bf R}_{m} = \tau^{-p}(R_m {\bf R}_{m-1})$,
hence
$$
 \Delta \wedge {\bf R}_{m} = \Delta \wedge \tau^{-p}(R_m {\bf R}_{m-1}) =
 \tau^{-p}(\Delta \wedge (R_m {\bf R}_{m-1})).
$$

Notice that, since $\inf({\bf R}_{m-1})=0$ by the induction
hypothesis, the initial factor of $R_m {\bf R}_{m-1}$ depends only
on $R_m$ and on the initial factor of ${\bf R}_{m-1}$, that is,
$$
\Delta \wedge (R_m {\bf R}_{m-1}) = \Delta \wedge (R_m (\Delta
\wedge {\bf R}_{m-1})) \preceq \Delta \wedge (R_m C_m).
$$
But $\Delta^p R_m C_m = \mathbf c^{m}(X)$, and $\iota(\mathbf
c^{m}(X))=C_{m+1}$, hence $\Delta\wedge (R_mC_m) = \tau^p(C_{m+1})$.
Therefore
$$
\Delta \wedge (R_m {\bf R}_{m-1}) \preceq \tau^p(C_{m+1}),
$$
and then
$$
 \Delta \wedge {\bf R}_{m} = \tau^{-p}(\Delta \wedge (R_m {\bf R}_{m-1})) \preceq
 \tau^{-p} (\tau^p (C_{m+1})) = C_{m+1},
$$
as we wanted to show.

\medskip

3.  We prove (3) by induction on $m$. If $m=1$ then ${\bf C}_1=C_1$.
As $\ell(X)>0$, $C_1$ is a non-trivial simple element, whence
$\sup(C_1)=\ell(C_1)=1$.  As $X\in SSS(X)$, the number of canonical
factors cannot decrease when passing from $X$ to $\mathbf c(X)$.  In
particular, the factor $C_1$ moved to the end cannot be absorbed
completely which shows $\varphi(\mathbf C_1) = C_1 \succeq
\varphi(\mathbf c(X))$.

Suppose the result true for some $m\geq 1$ and let $F = \varphi({\bf
C}_m)= \varphi(C_1 \cdots C_{m})$. As above we see that
$\varphi(\mathbf c^{m}(X)) \; C_{m+1}$ cannot be simple, as $X$ (and
hence $\mathbf c^{m}(X)$) is super summit. Notice that we used
$\ell(X)>1$ here. Since $F \succeq \varphi(\mathbf c^m(X))$ by
induction, this implies that $F C_{m+1}$ also has supremum 2. It is
well known~\cite{Michel} that if one multiplies a left normal form
$y_1\cdots y_m$ by a simple element $C_{m+1}$, then the left normal
form of the product is computed by by applying $m$ local
transformations to pairs of consecutive factors, starting by the
rightmost pair. In this way one can compute the normal form of every
element, so we will refer to this as the {\it normal form
algorithm}. Since $\sup(\mathbf C_m)=m$ by induction hypothesis, and
$F C_{m+1}$ is not simple, the normal form algorithm implies that
$\sup({\bf C}_{m+1})=\sup(C_1 \cdots C_{m+1}) = m+1$. This algorithm
together with $F \succeq \varphi(\mathbf c^m(X))$ and $\ell(X)>1$
also implies that
$$
    \varphi({\bf C}_{m+1})=
    \varphi( C_1 \cdots C_{m} C_{m+1} ) = \varphi( F  C_{m+1} )
    \succeq \varphi( \varphi(\mathbf c^m(X)) \; C_{m+1} ) =
    \varphi(\mathbf c^{m+1}(X)).
$$
Finally, since $C_{m+1}\neq \Delta$ and $F  C_{m+1}$ has supremum 2,
it follows that $\varphi({\bf C}_{m+1})= \varphi( F  C_{m+1} ) \neq
\Delta$, so $\ell({\bf C}_{m+1})>0$.  \qed

Since the super summit set of an element is closed under cycling,
Lemma~\ref{L:decomposition} is still true if we apply it to every
iterated cycling of $X$. If furthermore $X\in USS(X)$, the same will
be true for every element in $\mathcal O(X)$. This yields the
following result.

\begin{lemma}\label{L:decomposition whole orbit}
Let $X \in SSS(X)$, with $\ell(X)\geq 1$. With the above notation,
one has, for every $k\geq 0$ and every $m\geq 1$:

1. The $m$-th power of $\mathbf c^k(X)$ has the decomposition
$$
(\mathbf c^k(X))^m = {\bf C}_{[k,m]}
    {\bf R}_{[k,m]} \Delta^{mp}.
$$

2. In this decomposition, $\inf({\bf R}_{[k,m]})=0$ and $\iota({\bf
R}_{[k,m]})\preceq C_{k+m+1}$, for every $m\geq 1$.

3. In general $\inf({\mathbf C_{[k,m]}}) \geq 0$.  If $\ell(X)>1$,
one has $\sup({\bf C}_{[k,m]})=m$, $\ell({\bf C}_{[k,m]})>0$ and
$\varphi({\bf C}_{[k,m]})\succeq \varphi(\mathbf c^{k+m}(X))$, for
every $m\geq 1$.

Moreover, if $X\in USS(X)$, the result holds for every $k\in \mathbb
Z$.
\end{lemma}

Notice that in ${\bf C}_{[k,m]}$ and ${\bf R}_{[k,m]}$, the first
index determines an element in the cycling orbit of $X$, and the
second index determines its power. One can also think of $m$ as
being the number of factors in the decompositions of ${\bf
C}_{[k,m]}$ and ${\bf R}_{[k,m]}$ given by the definitions. But this
is not necessarily the number of factors in their normal forms.

\subsection{Interplay between $\mathbf C_m$ and $\mathbf R_m$.}
\label{subsection:Interplay between C_m and R_m}

Having proved that $X^m = \mathbf C_m \mathbf R_m \Delta^{pm}$, we
will show that if $X$ belongs to its ultra summit set, this
decomposition is left weighted, that is, $\varphi(\mathbf C_m)
\iota(\mathbf R_m)$ is in left normal form as written. In other
words, since we know by Lemma~\ref{L:decomposition} that
$\sup(\mathbf C_m)=m$, we will show that the product of the first
$m$ factors, including $\Delta$'s, in the left normal form of
$X^m\Delta^{-pm}$ is precisely $\mathbf C_m$.

If $X\in USS(X)$, recall that by Lemma~\ref{L: RC=CR}, one has
$$
  \mathbf  R_{[k-1,m]} \Delta^{pm} \; C_k  = C_{k+m} \;\mathbf R_{[k,m]}
  \Delta^{pm},
$$
for every $k\in \mathbb Z$. We will actually see that the initial
factor of this element is precisely $C_{k+m}$, no matter how many
remainders we multiply on the right, that is, no matter how big is
$m$.

\begin{lemma}\label{L: RC=CR initial factor}
Let $X\in SSS(X)$ with $\inf(X)=p$ and $\ell(X)\geq 1$. For every
$m\geq 1$ and $k\geq 0$, one has
$$
(C_{k+m} \;\mathbf R_{[k,m]})\wedge \Delta = C_{k+m}.
$$
If furthermore $X\in USS(X)$, this is also true for every $k\in
\mathbb Z$.
\end{lemma}

\begin{proof}
If $m=1$, we need to show that $(C_{k+1} \mathbf R_{[k,1]})\wedge
\Delta = C_{k+1}$, but we have $\mathbf c^k(X)=C_{k+1} \Delta^p
R_{k+1} = C_{k+1} \mathbf R_{[k,1]} \Delta^p $, where this
decomposition is left weighted by definition. So the result is true
for $m=1$.

Suppose the result true for $m-1$. This means that $(C_{k+m-1}
\;\mathbf R_{[k,m-1]})\wedge \Delta = C_{k+m-1}$.  If we multiply on
the left by $R_{k+m-1}$ we get $(R_{k+m-1} C_{k+m-1} \;\mathbf
R_{[k,m-1]})\wedge (R_{k+m-1}\Delta) = R_{k+m-1} C_{k+m-1}$. Notice
that $\Delta\preceq R_{k+m-1}\Delta$, hence if we consider the
maximal simple prefix of each element in the above equality, we
obtain
$$
(R_{k+m-1} C_{k+m-1} \;\mathbf R_{[k,m-1]})\wedge \Delta =
(R_{k+m-1} C_{k+m-1}) \wedge \Delta.
$$
On the other hand, recall that $(\Delta^p R_{k+m-1}) C_{k+m-1} =
\mathbf c^{k+m-1}(X) = C_{k+m} (\Delta^p R_{k+m})$, hence $R_{k+m-1}
C_{k+m-1} = \tau^p(C_{k+m}) R_{k+m}$, where $\tau^p(C_{k+m})$ is the
maximal simple prefix of this element. Therefore, one has
$$
 (R_{k+m-1} C_{k+m-1} \mathbf R_{[k,m-1]}) \wedge \Delta=
 (R_{k+m-1} C_{k+m-1})\wedge \Delta  = \tau^p(C_{k+m}).
$$
Now notice that
$$
\Delta^p R_{k+m-1} C_{k+m-1} \mathbf R_{[k,m-1]} =
 C_{k+m} \Delta^p R_{k+m} \mathbf R_{[k,m-1]} = C_{k+m} \mathbf
 R_{[k,m]} \Delta^p.
 $$
This means that
$$
 R_{k+m-1} C_{k+m-1} \mathbf R_{[k,m-1]} = \tau^p(C_{k+m} \mathbf
 R_{[k,m]}),
$$
and we just showed that its maximal simple prefix is precisely
$\tau^p(C_{k+m})$.  Applying $\tau^{-p}$ to this element, one
obtains $(C_{k+m} \;\mathbf R_{[k,m]})\wedge \Delta = C_{k+m}$, as
we wanted to show.

The proof for $X\in USS(X)$ and $k<0$ is the same.
\end{proof}

It will possibly help to understand the situation if we extract a
particular case from the above result, assuming that $X\in USS(X)$:

\begin{lemma}
Let $X\in USS(X)$ with $\inf(X)=p$ and $\ell(X)\geq 1$. For every
$m\geq 1$, one has
$$
(C_1 \;\mathbf R_{[1-m,m]})\wedge \Delta = C_1.
$$
In other words, the biggest simple prefix of
$$
 C_1 \tau^{-p}(R_1) \tau^{-2p}(R_0) \tau^{-3p}(R_{-1})\cdots \tau^{-mp}(R_{2-m})
$$
is $C_1$, no matter how big is $m$.
\end{lemma}

We can now show that the decomposition $\mathbf C_m \mathbf R_m$ is
left weighted. This will actually be a particular case of the
following stronger result.

\begin{proposition}\label{P:CR left weighted}
Let $X\in USS(X)$ with $\ell(X)\geq 1$. With the above notation, for
every $k,l,m,n\in \mathbb Z$, such that $ m,n\geq 1$ and $k+m=l+n$,
the decomposition $ \mathbf C_{[k,m]} \mathbf R_{[l,n]}$ is left
weighted. That is, $\varphi(\mathbf C_{[k,m]}) \iota(\mathbf
R_{[l,n]}) \wedge \Delta = \varphi(\mathbf C_{[k,m]})$.
\end{proposition}

\begin{proof}
We proceed by induction on $m$. If $m=1$, one has $ \mathbf
C_{[k,1]} \mathbf R_{[l,n]} = C_{k+1} \mathbf R_{[l,n]}$. Since
$k+1=l+n$ by hypothesis, it follows from Lemma~\ref{L: RC=CR initial
factor} that the biggest simple prefix of this element is precisely
$C_{k+1}$. Hence $ \mathbf C_{[k,1]} \mathbf R_{[l,n]}$ is left
weighted, and the result is true for $m=1$.

Now consider $k,l,m,n$ as above, with $m>1$, and suppose the result
true for $m-1$. This implies that $\mathbf C_{[k,m-1]} \mathbf
R_{[l-1,n]}$ is left weighted, since $k,l-1,m-1,n$ satisfy the
required hypothesis. (Notice that if we had required $k,l\geq 0$, we
would not have been able to apply the induction hypothesis here,
since we could have had $l-1<0$; This is why we require $X\in
USS(X)$ and not only in $SSS(X)$.)

Let $\mathbf C_{[k,m-1]}=\Delta^q y_1\cdots y_s$ and $\mathbf
R_{[l-1,n]}=z_1\cdots z_t$ in left normal form. Then $\mathbf
C_{[k,m-1]} \mathbf R_{[l-1,n]} = \Delta^q y_1\cdots y_s z_1\cdots
z_t$ is in left normal form as written, by induction hypothesis. Now
multiply this element on the right by $\tau^{-pn}(C_l)$.  By the
normal form algorithm, since $\Delta^q y_1\cdots y_s z_1\cdots z_t$
is already in left normal form and $\tau^{-pn}(C_l)$ is a simple
element, then the left normal form of $\Delta^q y_1\cdots y_s
z_1\cdots z_t \tau^{-pn}(C_l)$ is computed by by applying $s+t$
local transformations to pairs of consecutive factors, starting by
the rightmost pair. When we apply the first $t$ transformations, we
obtain the left normal form of $z_1\cdots z_t \tau^{-pn}(C_l) =
\mathbf R_{[l-1,n]} \tau^{-pn}(C_l)$. By Lemma~\ref{L: RC=CR}, this
element is equal to $C_{l+n} \mathbf R_{[l,n]} = C_{k+m} \mathbf
R_{[l,n]}$.  Moreover, by Lemma~\ref{L: RC=CR initial factor},
$C_{k+m}$ is the biggest simple prefix of this element. Hence the
left normal form of $z_1\cdots z_t \tau^{-pn}(C_l)$ has the form
$C_{k+m} z_1'\cdots z_t'$, where $z_1'\cdots z_t'= \mathbf
R_{[l,n]}$. We then have
$$
\mathbf C_{[k,m-1]} \mathbf R_{[l-1,n]} \tau^{-pn}(C_l) = \mathbf
C_{[k,m-1]} C_{k+m} \mathbf R_{[l,n]} =  \Delta^q y_1\cdots y_s
C_{k+m} z_1'\cdots z_t',
$$
where the last $t+1$ factors in the latter decomposition are in left
normal form.

If we continue applying the normal form algorithm, we perform $s$
local transformations to the element $y_1\cdots y_s C_{k+m}$, which
is equal to $\mathbf C_{[k,m-1]} C_{k+m}={\bf C}_{[k,m]}$. Since the
resulting factorization of $\mathbf C_{[k,m]} \mathbf R_{[l,n]}$ is
in left normal form by construction, it follows that $\varphi({\bf
C}_{[k,m]})\iota(\mathbf R_{[l,n]})$ is left weighted, as we wanted
to show.
\end{proof}

This result implies one of the strongest relations between cyclings
and powers of an element $X$ in its ultra summit set:

\begin{theorem}\label{T:normal form of $C_m$ and $X^m$}
Let $X\in USS(X)$ with $\inf(X)=p$ and $\ell(X)>1$. For every $m\geq
1$, the product of the first $m$ factors (including $\Delta$'s) in
the left normal form of $X^m \Delta^{-mp}$ is equal to $\mathbf
C_m$. That is,
$$
    (X^m \Delta^{-mp}) \wedge \Delta^m =  \mathbf C_m.
$$
In particular, $\iota(X^m)=\iota(\mathbf C_m)$.
\end{theorem}

\begin{proof}
The first claim is a straightforward consequence of the previous
result and Lemma~\ref{L:decomposition}, since $X^m \Delta^{-mp} =
\mathbf C_m \mathbf R_m$, where the latter decomposition is left
weighted and $\sup(\mathbf C_{m})=m$.

The second claim follows from the fact that $\iota(Y)=\iota(Y
\Delta^{t})$ for every $Y\in G$ and every $t\in \mathbb Z$. Hence
$\iota(X^m) = \iota(X^m \Delta^{-mp}) = \iota(\mathbf C_m \mathbf
R_m)$. Since the latter decomposition is left weighted, and
$\ell(\mathbf C_m)>0$, it follows that $\iota(X^m) = \iota(\mathbf
C_m)$, as we wanted to show.
\end{proof}

\begin{corollary}\label{C:inf(X^m)}
For $X\in USS(X)$ and $m\geq 1$, one has
$$
   \inf(X^m)= m \inf(X) +\inf(\mathbf C_m).
$$
In particular, the unexpected $\Delta's$ in $X^m$ are determined
entirely by the $\mathbf C_m$ part of the normal form of $X^m$.
\end{corollary}

\begin{proof}
Since we know by Lemma~\ref{L:decomposition} that $\ell(\mathbf
C_m)>0$, it follows from Theorem~\ref{T:normal form of $C_m$ and
$X^m$} that all $\Delta$'s in the left normal form of $X^m$ can be
seen in $\Delta^{pm} \tau^{pm}(\mathbf C_m)$, where $p=\inf(X)$.
Hence $\inf(X^m)= pm+\inf(\mathbf C_m)$, and the result follows.
\end{proof}

We end this section with an immediate corollary concerning how the
infimum and supremum of an element, in its ultra summit set, behave
when one raises the element to some power. This is related to the
translation number of the element (see~\cite{Lee-Lee3}). Notice that
the following result is closely related to Proposition 3.6
in~\cite{Lee-Lee}.

\begin{corollary}
Let $X\in USS(X)$. For every $m\geq 1$, one has
$$
    \inf(X^m)+\inf(X) \leq \inf(X^{m+1}) \leq \inf(X^m)+\inf(X)+1.
$$
If $X^{-1}\in USS(X^{-1})$, then
$$
    \sup(X^m)+\sup(X)-1 \leq \sup(X^{m+1}) \leq \sup(X^m)+\sup(X).
$$
\end{corollary}

\begin{proof}
Let $p=\inf(X)$. By the previous corollary,
$\inf(X^m)=pm+\inf(\mathbf C_m)$ and
$\inf(X^{m+1})=p(m+1)+\inf(\mathbf C_{m+1})$. Hence, the first
inequality will be true if and only if
$$
   \inf(\mathbf C_m) \leq \inf(\mathbf C_{m+1}) \leq \inf(\mathbf
   C_m)+1.
$$
But $\mathbf C_{m+1} = \mathbf C_m  C_{m+1}$, where $C_{m+1}$ is a
simple element. The result then follows from the following well know
fact, which is a direct consequence of the normal form algorithm: If
an element in a Garside group is multiplied by a simple element,
then its infimum either is preserved or is increased by one.

The second inequality is equivalent to the first one, since
$\sup(X)=-\inf(X^{-1})$ by Corollary~\ref{C:inf, sup and length of
the inverse}.
\end{proof}

\subsection{The absolute initial and final factors}
\label{subsection:the absolute initial and final factors}

In this section we will define some simple factors related to an
element $X\in USS(X)$. They are defined in terms of the cycling
elements $C_i$, but they are closely related to powers of $X$, as we
will see. We called them the {\em absolute} initial and final
factors of $X$.

In general, the absolute initial and final factors are related to,
but do not coincide with, the initial and final factors of $X$.
Nevertheless we will see that, if $X$ has a rigid power $X^m$, the
absolute initial and final factors of $X$ coincide with the initial
and final factors of $X^m$.

Suppose that $X\in USS(X)$. We saw in Lemma~\ref{L:decomposition
whole orbit} that
$$
 \varphi({\bf C}_{[k,m]}) = \varphi(C_{k+1}\cdots C_{k+m}) \succeq \varphi(\mathbf
 c^{k+m}(X)),
$$
where this is true for every $m\geq 1$ and every $k\in \mathbb Z$.
This implies a very interesting fact: if we fix the number $k+m$ and
take different values of $k$, that is, if we start with $C_{k+m}$
and multiply it on the left by $C_{k+m-1}$, then by $C_{k+m-2}$,
etc., then the final factor of each of the resulting elements is a
left multiple of $\varphi(\mathbf c^{k+m}(X))$. For instance, if we
take $k+m=0$, we have
$$
   \varphi(C_{-m+1} C_{-m+2} \cdots C_{-1} C_0) \succeq \varphi(X),
$$
for every $m\geq 1$.

In the same way, by Lemma~\ref{L:decomposition whole orbit} we know
that
$$
\iota({\bf R}_{[k,m]}) = \iota(\tau^{-p}(R_{k+m})
\tau^{-2p}(R_{{k+m-1}})\cdots \tau^{-mp}(R_{k+1})) \preceq C_{k+m+1}
$$
for every $m\geq 1$ and every $k\in \mathbb Z$, where $p=\inf(X)$.
Hence, if we fix $k+m$, say $k+m=0$, we have (recall that
$C_1=\iota(X)$)
$$
   \iota(\tau^{-p}(R_{0})
\tau^{-2p}(R_{{-1}})\cdots \tau^{-mp}(R_{-m+1})) \preceq  \iota(X),
$$
for every $m\geq 1$. In the particular case in which $\inf(X)=p=0$,
this formula is even more similar to the above one, since one has:
$$
   \iota(R_0 R_{-1} \cdots R_{-m+1}) \preceq  \iota(X),
$$
for every $m\geq 1$.

But it is even more interesting to relate the values of
$$
  \varphi(C_{-m+1} C_{-m+2} \cdots C_{-1} C_0) \quad \mbox{and of}
  \quad \iota(\tau^{-p}(R_{0})
\tau^{-2p}(R_{{-1}})\cdots \tau^{-mp}(R_{-m+1})),
$$
respectively, for different values of $m$. It turns out that they
form ordered chains with respect to $\succeq$ and $\preceq$,
respectively, as shown in the following result.

\begin{proposition}\label{P:chains}
Let $X\in USS(X)$ with $\ell(X)>1$. For every $k,m\in \mathbb Z$
with $m>0$, one has
$$
    \varphi({\bf C}_{[k,m]}) \succeq \varphi({\bf C}_{[k-1,m+1]}) \hspace{1cm}
    \mbox{and} \hspace{1cm}    \iota({\bf R}_{[k,m]}) \preceq
    \iota({\bf R}_{[k-1,m+1]}).
$$
In other words, for every $k\in \mathbb Z$ there are chains:
$$
   \varphi({\bf C}_{[k,1]}) \succeq \varphi({\bf C}_{[k-1,2]}) \succeq \varphi({\bf C}_{[k-2,3]}) \succeq \cdots
$$
and
$$
   \iota({\bf R}_{[k,1]}) \preceq \iota({\bf R}_{[k-1,2]}) \preceq \iota({\bf R}_{[k-2,3]}) \preceq \cdots
$$
\end{proposition}

\begin{proof}
We know that $\sup({\bf C}_{[k,m]})=m$ and $\sup({\bf
C}_{[k-1,m+1]})=m+1$. Moreover, ${\bf C}_{[k-1,m+1]}= C_k C_{k+1}
\cdots C_{k+m} = C_k {\bf C}_{[k,m]}$. Hence, if we write ${\bf
C}_{[k,m]}=\Delta^q c_1\cdots c_s$ in left normal form (where
$q+s=m$ and $s>0$), then $\varphi({\bf C}_{[k,m]})= c_s \succeq
\varphi(C_k \Delta^q c_1\cdots c_s) =\varphi({\bf C}_{[k-1,m+1]})$.

On the other hand, let $p=\inf(X)$. One has
$$
 {\bf R}_{[k-1,m+1]}= \tau^{-p}(R_{k+m})
\tau^{-2p}(R_{k+m-1}) \cdots \tau^{-mp}(R_{k+1})
 \tau^{-(m+1)p}(R_{k}) = {\bf R}_{[k,m]} \cdot \tau^{-(m+1)p}(R_{k}),
$$
that is, ${\bf R}_{[k,m]}\preceq {\bf R}_{[k-1,m+1]}$. Since we know
that $\inf({\bf R}_{[k,m]})= \inf({\bf R}_{[k-1,m+1]}) = 0$, it
follows that $\iota({\bf R}_{[k,m]})\preceq \iota({\bf
R}_{[k-1,m+1]})$.
\end{proof}

Since the chains given by the above proposition consist of proper
simple elements, we know that the chains must stabilize. But we will
furthermore show that they stabilize fast. More precisely, they
stabilize exactly at the first repetition. Moreover, the
corresponding chains for all elements in ${\cal O}(X)$ stabilize at
the same time. This is proved by the next 4 lemmas and the
proposition that follows them.

\bigskip

\begin{lemma}\label{L:suffix equality transported}
Let $X\in USS(X)$ with $\ell(X)>1$ and $m>0$. If $\varphi({\bf
C}_{[k,m]})= \varphi({\bf C}_{[k-1,m+1]})$ for some $k \in \mathbb
Z$, then $\varphi({\bf C}_{[i,m]})= \varphi({\bf C}_{[i-1,m+1]})$
for every $i\in \mathbb Z$.
\end{lemma}

\begin{proof}
Since $X$ belongs to a closed orbit under cycling, the sequences
$\{{\bf C}_{[i,m]}\}_{i\in \mathbb Z}$ and $\{{\bf
C}_{[i,m+1]}\}_{i\in \mathbb Z}$ are periodic, hence it suffices to
show the case $i=k+1$.

The property $\varphi({\bf C}_{[k,m]})= \varphi({\bf
C}_{[k-1,m+1]})$, that is, $\varphi(C_{k+1}\cdots C_{k+m}) =
\varphi(C_{k}\cdots C_{k+m})$, can be reinterpreted as follows.
Since $\sup(C_{k+1}\cdots C_{k+m})=m$ and $\sup(C_{k}\cdots
C_{k+m})=m+1$, their final factors coincide if and only if the first
$m-1$ factors of the first element, multiplied on the left by $C_k$,
coincide with the first $m$ factors of the second element. In other
words,
$$
 (C_k \cdots C_{k+m})\wedge (C_k\Delta^{m-1}) = (C_k\cdots C_{k+m})\wedge
\Delta^m.
$$
We can now apply Gebhardt's transport~\cite{Gebhardt} to the whole
equality. We know that the transport of $\Delta^m$ is $\Delta^m$.
Notice that the transport of $C_k$ (based at $\mathbf c^{k-1}(X)$)
is $C_{k+1}$. This implies, by recurrence, that the transport of
$C_{k}\cdots C_{k+m}$ is $C_{k+1}\cdots C_{k+m+1}$, and also that
the transport of $C_{k}\Delta^{m-1}$ is $C_{k+1}\Delta^{m-1}$. Since
the transport preserves greatest common divisors, the transport of
the above equality yields
$$
 (C_{k+1} \cdots C_{k+m+1})\wedge (C_{k+1}\Delta^{m-1}) = (C_{k+1}\cdots C_{k+m+1})\wedge
\Delta^m,
$$
that is, $\varphi({\bf C}_{[k+1,m]})= \varphi({\bf C}_{[k,m+1]})$,
and the result is shown.
\end{proof}

\bigskip

\begin{lemma}\label{L:suffix equality stabilizes}
Let $X\in USS(X)$ with $\ell(X)>1$. If $\varphi({\bf C}_{[k,m]})=
\varphi({\bf C}_{[k-1,m+1]})$ for some $k,m \in \mathbb Z$ with
$m>0$, then $\varphi({\bf C}_{[i,j]})= \varphi({\bf C}_{[i-1,j+1]})$
for every $i\in \mathbb Z$ and every $j\geq m$.
\end{lemma}

\begin{proof}
We know by Lemma~\ref{L:suffix equality transported} that
$\varphi({\bf C}_{[i,m]})= \varphi({\bf C}_{[i-1,m+1]})$ for every
$i\in \mathbb Z$. We just need to be able to increase the second
subindex.  But if $\varphi(C_{i+1} \cdots C_{i+m})= \varphi(C_{i}
\cdots C_{i+m}) \neq \Delta$, and we multiply both elements on the
right by $C_{i+m+1}$, since we know that no unexpected decrease of
supremum will happen ($\sup(\mathbf C_{[i,m+1]})=m+1$), it follows
that
$$
\varphi(C_{i+1} \cdots C_{i+m} C_{i+m+1}) = \varphi(\varphi(C_{i+1}
\cdots C_{i+m}) C_{i+m+1})
$$
$$
= \varphi(\varphi(C_{i} \cdots C_{i+m}) C_{i+m+1}) = \varphi(C_{i}
\cdots C_{i+m} C_{i+m+1}).
$$
Hence $\varphi({\bf C}_{[i,m+1]})=\varphi({\bf C}_{[i-1,m+2]})$ for
every $i\in \mathbb Z$. By induction on $m$, it follows that
$\varphi({\bf C}_{[i,j]})=\varphi({\bf C}_{[i-1,j+1]})$ for every
$j\geq m$, as we wanted to show.
\end{proof}
The analogous results can now be shown for the chain involving
prefixes of ${\bf R}_{[k,m]}$.

\bigskip
\begin{lemma}\label{L:prefix equality transported}
Let $X\in USS(X)$ with $\ell(X)>1$ and $m\geq 1$. If $\iota({\bf
R}_{[k,m]})= \iota({\bf R}_{[k-1,m+1]})$ for some $k \in Z$, then
$\iota({\bf R}_{[i,m]})= \iota({\bf R}_{[i-1,m+1]})$ for every $i\in
\mathbb Z$.
\end{lemma}

\begin{proof}
As above, since $X$ belongs to a closed orbit under cycling, it
suffices to show the case $i=k+1$.  We want to reinterpret the
equality $\iota({\bf R}_{[k,m]})= \iota({\bf R}_{[k-1,m+1]})$. If we
recall that $(\mathbf c^{k}(X))^{m} = {\bf C}_{[k,m]} {\bf
R}_{[k,m]} \Delta^{mp}$ where $p=\inf(X)$, and that $\inf(\mathbf
R_{[k,m]})=0$, then we see that
$$
 (\mathbf c^{k}(X))^{m}\wedge \left({\bf C}_{[k,m]} \Delta^{mp+1}\right) = {\bf C}_{[k,m]}
 \: \iota({\bf R}_{[k,m]}) \:\Delta^{mp} = C_{k+1}\cdots C_{k+m} \: \iota({\bf R}_{[k,m]}) \: \Delta^{mp}.
$$
In the same way, we obtain
$$
 (\mathbf c^{k-1}(X))^{m+1}\wedge \left({\bf C}_{[k-1,m+1]} \Delta^{(m+1)p+1}\right) = {\bf C}_{[k-1,m+1]}
 \: \iota({\bf R}_{[k-1,m+1]}) \:\Delta^{(m+1)p}
$$
$$
 = C_{k}\cdots C_{k+m} \: \iota({\bf R}_{[k-1,m+1]}) \: \Delta^{mp+p}.
$$
Therefore the equality $\iota({\bf R}_{[k,m]})= \iota({\bf
R}_{[k-1,m+1]})$ can be rewritten as follows:
$$
   C_k\: \left[(\mathbf c^{k}(X))^{m} \wedge ({\bf C}_{[k,m]} \: \Delta^{mp+1})\right] \Delta^p
   = (\mathbf c^{k-1}(X))^{m+1} \wedge \left({\bf C}_{[k-1,m+1]} \Delta^{(m+1)p+1}\right).
$$
If we apply Gebhardt's transport to the whole equality, it follows
that
$$
   C_{k+1}\: \left[(\mathbf c^{k+1}(X))^{m} \wedge ({\bf C}_{[k+1,m]} \: \Delta^{mp+1})\right] \Delta^p
   = (\mathbf c^{k}(X))^{m+1} \wedge \left({\bf C}_{[k,m+1]}
   \Delta^{(m+1)p+1}\right),
$$
hence $\iota({\bf R}_{[k+1,m]})= \iota({\bf R}_{[k,m+1]})$, and the
result is shown.
\end{proof}

\bigskip

\begin{lemma}\label{L:prefix equality stabilizes}
Let $X\in USS(X)$ with $\ell(X)>1$. If $\iota({\bf R}_{[k,m]})=
\iota({\bf R}_{[k-1,m+1]})$ for some $k,m \in Z$ with $m>0$, then
$\iota({\bf R}_{[i,j]})= \iota({\bf R}_{[i-1,j+1]})$ for every $i\in
\mathbb Z$ and every $j\geq m$.
\end{lemma}

\begin{proof}
We know by Lemma~\ref{L:prefix equality transported} that
$\iota({\bf R}_{[i,m]})= \iota({\bf R}_{[i-1,m+1]})$ for every $i\in
\mathbb Z$. We just need to be able to increase the second subindex.
But $\iota({\bf R}_{[i,m]})= \iota({\bf R}_{[i-1,m+1]})$ is
equivalent to $\iota((\Delta^p R_{i+m}) (\Delta^p R_{i+m-1}) \cdots
(\Delta^p R_{i+1})) = \iota((\Delta^p R_{i+m}) (\Delta^p R_{i+m-1})
\cdots (\Delta^p   R_{i}))$, where $p=\inf(X)$. If we multiply both
elements on the left by $\Delta^p R_{i+m+1}$, since we know that
there is no unexpected increase of infimum ($\inf(\mathbf
R_{[i,m+1]})=0$), it follows that
$$
\iota((\Delta^p R_{i+m+1})(\Delta^p R_{i+m})\cdots (\Delta^p
R_{i+1})) = \iota((\Delta^p R_{i+m+1}) \iota((\Delta^p R_{i+m})
\cdots (\Delta^p R_{i+1})))
$$
$$
= \iota((\Delta^p R_{i+m+1}) \iota((\Delta^p R_{i+m}) \cdots
(\Delta^p R_{i}))) = \iota((\Delta^p R_{i+m}) (\Delta^p R_{i+m-1})
\cdots (\Delta^p   R_{i})).
$$
Hence $\iota({\bf R}_{[i,m+1]})=\iota({\bf R}_{[i-1,m+2]})$ for
every $i\in \mathbb Z$. By induction on $m$, it follows that
$\iota({\bf R}_{[i,j]})=\iota({\bf R}_{[i-1,j+1]})$ for every $j\geq
m$, as we wanted to show.
\end{proof}

\bigskip

\begin{proposition}\label{P:suffix chain stabilizes}
Let $X\in USS(X)$ with $\ell(X)> 1$. Given $k\in \mathbb Z$, the
chain
$$
   \varphi({\bf C}_{[k,1]}) \succeq \varphi({\bf C}_{[k-1,2]}) \succeq \varphi({\bf C}_{[k-2,3]}) \succeq \cdots
$$
stabilizes whenever $\varphi({\bf C}_{[k-j+1,j]})=\varphi({\bf
C}_{[k-j,j+1]})$, and this happens for some $j<||\Delta||$.
Moreover, for all $i\in \mathbb Z$, the analogous chains starting at
$\varphi({\bf C}_{[i,1]})$ stabilize at the same value of $j$. Also,
the chain
$$
   \iota({\bf R}_{[k,1]}) \preceq \iota({\bf R}_{[k-1,2]}) \preceq \iota({\bf R}_{[k-2,3]}) \preceq \cdots
$$
stabilizes whenever $\iota({\bf R}_{[k-j+1,j]})=\iota({\bf
R}_{[k-j,j+1]})$, and this happens for some $j<||\Delta||$.
Moreover, for all $i\in \mathbb Z$, the analogous chains starting at
$\iota({\bf R}_{[i,1]})$ stabilize at the same value of $j$.
\end{proposition}

\begin{proof}
By Lemma~\ref{L:suffix equality stabilizes}, all chains stabilize
whenever $\varphi({\bf C}_{[k-j+1,j]})=\varphi({\bf C}_{[k-j,j+1]})$
for some $j$. Up to that point, the chains must be made of strict
inequalities. But the maximal length of such a chain (formed by
nontrivial simple elements) is bounded by the length of $\Delta$.
The proof that the second sequence stabilizes is  identical.
\end{proof}

\begin{definition}{\rm Given $X\in USS(X)$ with $\ell(X)> 1$, we define
the {\em absolute final factor} $F(X)$ of $X$ as the factor in which
the above descending chain stabilizes, for $k=-1$, that is:
$$
    F(X) = \varphi(\mathbf C_{[-m, m]}),
$$
for $m\geq ||\Delta||-1$. In other words,
$
  F(X) = \varphi(C_{-m+1}C_{-m+2}\cdots C_{-1} C_0),
$ for $m$ big enough. }
\end{definition}

\begin{definition}{\rm Given $X\in USS(X)$ with $\ell(X)>1$, we define the
{\em absolute initial factor} $I(X)$ of $X$ as the factor in which
the above ascending chain stabilizes, for $k=-1$, that is:
$$
    I(X) = \iota(\mathbf R_{[-m, m]}),
$$
for $m\geq ||\Delta||-1$. In other words,
$
  I(X) = \iota(\tau^{-p}(R_0) \tau^{-2p}(R_{-1}) \cdots \tau^{-mp}(R_{-m+1})),
$
or alternatively
$
  I(X) = \iota((\Delta^p R_0) (\Delta^p R_{-1}) \cdots (\Delta^p  R_{-m+1})),
$ for $m$ big enough, where $p=\inf(X)$.}
\end{definition}

\begin{proposition}\label{P:FI left weighted}
Given $X\in USS(X)$ with $\ell(X)>1$, the decomposition $F(X)I(X)$
is left weighted as written.
\end{proposition}

\begin{proof}
This is an immediate consequence of Proposition~\ref{P:CR left
weighted}, since
$$
 F(X)I(X) = \varphi(\mathbf C_{[-m,m]})
\iota(\mathbf R_{[-m,m]})
$$
for $m$ big enough.
\end{proof}

We have seen at the beginning of this section that $F(X)\succeq
\varphi(X)$ and $I(X)\preceq \iota(X)$. But we will see now that the
absolute factors are also related to the initial and final factors
of {\em powers} of $X$.

\begin{proposition}
Let $X\in USS(X)$ with $\inf(X)=p$ and $\ell(X)>1$. For every $m\geq
1$ such that $X^m\in SSS(X^m)$, one has:
\begin{itemize}
 \item $F(X)\succeq \varphi(X^m)$.

 \item $I(X)\preceq \iota(X^m)$.
\end{itemize}
\end{proposition}

\begin{proof}
The case $m=1$ is a straightforward consequence of
Lemma~\ref{L:decomposition whole orbit}, for $k=-m$.

Suppose that $m>1$. We know from Theorem~\ref{T:normal form of $C_m$
and $X^m$} that the left normal form of
$$
(\mathbf c^{-m}(X))^m= (C_{-m+1}\cdots C_0)X^m (C_0^{-1}\cdots
C_{-m+1}^{-1})
$$
is equal to
$$
\Delta^{pm+q} y_1\cdots y_s z_1\cdots z_t,
$$
where $\Delta^q y_1\cdots y_s = \tau^{pm}(C_{-m+1}\cdots C_0)$ and
$q+s=m$. If we conjugate this element by $C_{-m+1}\cdots C_0$ we
obtain
$$
 X^m = \Delta^{pm} z_1\cdots z_t \Delta^q \tau^{-pm}(y_1\cdots y_s).
$$

But if $X^m\in SSS(X^m)$, we also have $Y^m\in SSS(X^m)$ for every
$Y$ in the cycling orbit of $X$: Indeed, $Y^m=(X^{\mathbf C_t})^m =
(X^m)^{\mathbf C_t}$ for some $t$, where $\mathbf
C_t=(X^m\Delta^{-mp})\wedge \Delta^t$ by Theorem~\ref{T:normal form
of $C_m$ and $X^m$}. Since $X^m\Delta^{-mp}$ and $\Delta^t$
conjugate $X^m$ to elements in their super summit sets (namely
$\tau^{-mp}(X^m)$ and $\tau^t(X^m)$, respectively), it follows by
Theorem~\ref{T:convexity_SSS} that $Y^m= (X^m)^{\mathbf C_t}\in
SSS(X^m)$.  In particular $(\mathbf c^{-m}(X))^m \in SSS(X^m)$,
hence $\ell(X^m)= \ell((\mathbf c^{-m}(X))^m) = s+t$. Since the
above decomposition of $X^m$ has precisely $s+t$ non-$\Delta$
factors, and the final one is $ \tau^{-pm}(y_s)$, it follows that
$\tau^{-pm}(y_s) \succeq \varphi(X^m)$. That is,
$\varphi(C_{-m+1}\cdots C_0)\succeq \varphi(X^m)$.

Notice that we can apply the same reasoning to every element in the
cycling orbit of $X$, in particular to $\mathbf c^{-m}(X)$. It
follows that $\varphi(C_{-2m+1}\cdots C_{-m}) \succeq
\varphi((\mathbf c^{-m}(X))^m) = z_t$. Hence, since multiplying
$C_i$'s never decreases the supremum, one has
$$
\varphi(C_{-2m+1}\cdots C_{-m} C_{-m+1} \cdots C_0) =
\varphi(\varphi(C_{-2m+1}\cdots C_{-m}) C_{-m+1} \cdots C_0)
$$
$$
\succeq \varphi(z_t \Delta^q \tau^{-pm}(y_1\cdots y_s)) =
\varphi(X^m).
$$
Applying the same reasoning again, one obtains by induction on $k$
that $\varphi(C_{-km+1}\cdots C_0) \succeq \varphi(X^m)$ for every
$k\geq 1$. When $k$ is big enough so that $km\geq ||\Delta||-1$,
this implies $F(X)\succeq \varphi(X^m)$, as we wanted to show.

The relation $I(X)\preceq \iota(X^m)$ is shown in a similar way.
Since $\tau^{-pm}(z_1\cdots z_t) = \mathbf R_{[-m,m]}$, from the
above decomposition of $X^m$ is follows that $\iota(\mathbf
R_{[-m,m]})= \tau^{-pm}(z_1) \preceq \iota(X^m)$. Applying the same
reasoning to $\mathbf c^{-m}(X)$, is follows that $\iota(\mathbf
R_{[-2m,m]}) \preceq \iota((\mathbf c^{-m}(X))^m)=
\tau^{-pm-q}(y_1)$. Hence
$$
   \iota(\mathbf R_{[-2m,2m]})= \iota(\mathbf R_{[-m,m]}\:\tau^{pm}(\mathbf
   R_{[-2m,m]})) = \iota(\mathbf R_{[-m,m]}\:\iota(\tau^{pm}(\mathbf
   R_{[-2m,m]})))
$$
$$
 \preceq \iota (\mathbf R_{[-m,m]}\: \tau^{-q}(y_1)) = \iota(\Delta^{pm} z_1\cdots z_t \Delta^q
 \tau^{-mp}(y_1)) = \iota(X^m).
$$
Iterating the same reasoning one shows that $\iota(\mathbf
R_{[-km,km]})\preceq \iota(X^m)$ for every $k\geq 1$, and when $k$
is big enough this yields $I(X)\preceq \iota(X^m)$.
\end{proof}

\subsection{The stable ultra summit set.} \label{subsection:the stable
ultra summit set}

We have studied, up to now, how powers and cyclings of $X$ are
related under the hypothesis, in most cases, that $X\in USS(X)$. But
this fact does not imply that $X^m\in USS(X^m)$ for every $m\in
\mathbb Z$, not even for every $m\in \mathbb N$. If we want to
extract more information from the powers of $X$, it would be
desirable that all these powers belonged to their ultra summit sets.

\begin{definition}{\rm Given $X$ in a Garside group $G$, the {\em stable
ultra summit set} of $X$ is defined as
$$
 SU(X)=\{Y\in USS(X); \quad Y^m\in USS(X^m), \; m\in \mathbb Z\}.
$$
}
\end{definition}

The first obvious question is whether $SU(X) \not= \emptyset$.

\begin{proposition}
\label{P:SU is not emptyset} For every $X\in G$, the set $SU(X)$ is
nonempty.
\end{proposition}

\begin{proof}
We can clearly assume that $\ell(Y)>0$ for every $Y\in USS(X)$.
Given an element $Z=\Delta^p z_1\cdots z_r$, its initial factor
$\iota(Z)$ can be described as $Z \Delta^{-p} \wedge \Delta$,
because $Z \Delta^{-p}$ is equal to $\tau^{-p}(z_1\cdots z_r)$,
hence the initial factor of $Z$ is equal to the first factor of
$Z\Delta^{-p}$, which is computed by considering its gcd with
$\Delta$.

Recall that $\iota(Z)$ is also the conjugating element for cycling.
On the other hand, the conjugating element for decycling is
$\Delta^p z_1\cdots z_{r-1}$, which can be described as $Z\wedge
\Delta^{p+r-1}$  (even if $p$ is negative).

Therefore, if we want to perform a cycling or a decycling to $Z$, we
must conjugate it by $Z \Delta^{-p} \wedge \Delta$ or by $Z \wedge
\Delta^{p+r-1}$, respectively.

Now consider an element $X$.  For every $Y \in USS(X)$, define
$S(Y)$ to be the set of integers $k$ such that $Y^k$ belongs to its
$USS$. Let $V\in USS(X)$ be such that $S(V)$ is maximal. Such a $V$
exists because $USS(X)$ is finite. We will see that
$S(V)=\mathbb{Z}$, hence $V\in SU(X)$.

Suppose that $S(V) \neq \mathbb{Z}$. This means that some power of
$V$, say $V^k$, does not belong to its USS. We would then like to
apply cyclings and decyclings to $V^k$ to bring it into $USS(X^k)$.
Let us conjugate $V$ by $\iota(V^k)$, to obtain some $W$. In this
way, {\bf all} powers of $V$ will be conjugate by $\iota(V^k)$. In
particular, $W^k$ will be the cycling of $V^k$. Moreover, if some
other power of $V$, say $V^s$, belongs to its USS, then $W^s$ also
belongs to its USS. Indeed, let $m=\inf(V^k)$. It is clear that
$V^k\Delta^{-m}$ conjugates $V^s$ to an element in its USS, namely
$\tau^{-m}(V^s)$. In the same way, $\Delta$ conjugates $V^s$ to an
element in its USS, $\tau(V^s)$. Therefore, by
Theorem~\ref{T:convexity_USS}, $\iota(V^k)= (V^k \Delta^{-m})\wedge
\Delta$ conjugates $V^s$ to an element ($W^s$) in its USS.

Therefore, the set $S(W)$ contains $S(V)$, and the power $W^k$ is
the cycling of $V^k$. The same can be done for decycling, since the
conjugating element for decycling $V^k$ is $V^k \wedge \Delta^t$
(for some $t$). Hence, by suitable conjugations of $V$, we can apply
iterated cyclings and decyclings to $V^k$, until we obtain a
conjugate $Z$ of $V$ such that $S(Z)$ contains $S(V)$, and $Z^k$
belongs to its USS. But then $S(Z)$ strictly contains $S(V)$, which
contradicts the maximality of $S(V)$.  This shows that $SU(X)$ is
non empty. \end{proof}

\ms

{\bf Remark:} Although we had a different name for the {\em stable
ultra summit set}, we chose the latter when we learnt about the
paper~\cite{Lee-Lee}, in which the {\em stable super summit set} was
defined in a similar way as above, but considering $Y^m\in SSS(Y^m)$
for every $m\in \mathbb N$. We remark that Proposition~\ref{P:SU is
not emptyset} was made public by the authors at a meeting on braid
groups held in Luminy, in June 2005, some months before the
appearance of~\cite{Lee-Lee}. Notice also that the proof of
Proposition~\ref{P:SU is not emptyset} can be applied to show the
non-emptyness of the stable super summit set, using
Theorem~\ref{T:convexity_SSS} instead of
Theorem~\ref{T:convexity_USS}. Moreover, one can extend the set of
exponents to the whole $\mathbb Z$, in both cases.  The proof of the
nonemptyness of the stable super summit set in~\cite{Lee-Lee} is
much more involved, and was found independently from ours.

\ms

Now notice that the proof of Proposition~\ref{P:SU is not emptyset}
yields an algorithm to compute the set
$$
\{Y\in USS(X); \quad Y^m\in USS(X^m), \;   m \in [A,B] \}
$$
for every pair of integers $A<B$, that is, we can assume that $X^m$
belongs to its ultra summit set for all integers between $A$ and
$B$. A priori, no matter how big is the interval $[A,B]$, this does
not say that $X\in SU(X)$, since one could have $X^t\not\in
USS(X^t)$ for some $t\not\in [A,B]$. But for our purposes we will
only need that $X^m\in USS(X^m)$ for $m\in [1,\: ||\Delta||\:]$,
hence the proof of Proposition~\ref{P:SU is not emptyset} allows us
to assume this hypothesis.

{\bf Remark:} Very recently we learnt that in~\cite{Lee-Lee3}, a
finite time algorithm to compute the stable {\em super} summit
set~\cite{Lee-Lee} is given. It is possible that similar methods can
be used to compute $SU(X)$ in finite time, but as we said above, we
will not need that for our purposes in this paper.

Let us then assume that $X^m\in USS(X^m)$ for $m=1,\ldots,
||\Delta||$. Recall that we have defined some factors $C_i$, $R_i$,
$\mathbf C_i$ and $\mathbf R_i$ related to $X$, for every $i\in
\mathbb Z$. We can thus define the same elements related to each
$X^m$, but we need some notation to make the distinction between
them, for different values of $m$. The notation
$C_i(X^m)=\iota(\mathbf c^{i-1}(X^m))$ would not cause confusion,
but it would be too awkward for the formulae below, so we will
simplify it by denoting:
$$
   C_i^{(m)} = C_i(X^m) = \iota(\mathbf c^{i-1} (X^m)).
$$
Thus $C_i^{(m)}$ has the same definition as $C_i$, but related to
$X^m$ instead of $X$. Later on, we will study the sequence
$$
     C_1, C_1^{(2)}, C_1^{(3)}, \ldots
$$
that is,
$$
   \iota(X), \iota(X^2), \iota(X^3),\ldots
$$
Notice that, a priori, there does not have to be a relation between
them, due to the unexpected increases of supremum that one
encounters when taking powers.

Similarly to $C_i^{(m)}$, one defines $R_i^{(m)}$, $\mathbf
C_i^{(m)}$ and $\mathbf R_i^{(m)}$ in the same way as $R_i$,
$\mathbf C_i$ and $\mathbf R_i$, but related to $X^m$ instead of
$X$. The relation between these elements for different powers of $X$
will be crucial in the sequel.

\section{Rigidity}
\label{section:rigidity}

\subsection{Rigidity of an element and behavior under cyclings and powers.}

In this section we will define a notion of rigidity for elements in
a Garside group, and we will study how rigidity is affected when
applying some cyclings or taking some powers. The idea of studying
rigidity came from the study of elements whose left normal form
changes only in the obvious way under cyclings, decyclings and
powers, so their ultra summit sets are easier to study. We call them
{\it rigid} elements:

\begin{definition}[\bf rigid element] \label{D:rigid element}
{\rm Let $X = \Delta^p x_1\cdots x_r$ be in left normal form, with
$r>0$. Then $X$ is {\it rigid} if the element $\Delta^p\: x_1\cdots
x_r
 \:\tau^{-p}(x_1)\, $ is in left normal form as written.}
\end{definition}

Notice that if $X$ is rigid, then the cycling of $X$, that is,
$\mathbf c(X) = \Delta^p\: x_2\cdots x_r
 \:\tau^{-p}(x_1)\, $ is in left normal form as written. Actually,
this latter property is equivalent to $X$ being rigid if $r>1$. But
we prefer the definition above, otherwise every element of canonical
length 1 would be rigid.

The following are equivalent definitions of rigid elements.

\begin{proposition}\label{P:definitions of rigidity}
Given $X=\Delta^p x_1\cdots x_r \in G$ with $r>0$, the following
conditions are equivalent.
\begin{enumerate}

\item $X$ is rigid.

\item $\varphi(X)\iota(X)$ is left weighted as written.

\item $\iota(X)\wedge \iota(X^{-1})=1$.

\end{enumerate}
\end{proposition}

\begin{proof}
By definition $X$ is rigid if $\Delta^p x_1\cdots x_r
\tau^{-p}(x_1)$ is in left normal form as written. Since $x_1\cdots
x_r$ is already in left normal form, this is equivalent to the left
weightedness of $x_r \tau^{-p}(x_1)=\varphi(X)\iota(X)$ so
conditions 1 and 2 are equivalent.  But condition 2 means
$\tau^{-p}(x_1) \wedge
\partial(x_r) =1$. We know that $\iota(X)=\tau^{-p}(x_1)$ and also,
by Lemma~\ref{L:initial & final factors of inverses},
$\iota(X^{-1})=\partial (x_r)$. Hence conditions 2 and 3 are also
equivalent.
\end{proof}

In general, we define the {\it rigidity} of an element, in such a
way that rigid elements have rigidity 1.

\begin{definition} {\rm Given $X=\Delta^p x_1\cdots x_r$ in left normal
form, with $r>0$, we define the {\it rigidity} of $X$ as
$$
    \mathcal R(X) = k/r,
$$
where $k$ is the biggest integer in $\{0,\ldots,r\}$ such that the
first $k$ factors in the left normal form of $x_1\cdots x_r
\tau^{-p}(x_1)$ are precisely $x_1\cdots x_k$.  If $r=0$, we define
$\mathcal R(X)=0$.}
\end{definition}

The rigidity of an element tells us how many (non-$\Delta$) factors
of the left normal form of $X$ are preserved when considering $X^2$.
Notice that $X$ is rigid if and only if $\mathcal R(X)=1$.

{\bf Examples:} \be

\item If $X=\Delta \cdot 12\cdot 21 \cdot 12 \in B_3$,
then $\mathcal R(X)=3/3=1$, since $\iota(X)=\tau^{-1}(12)=21$, and
$12\cdot 21 \cdot 12 \cdot 21$ is in left normal form as written.
Hence $X$ is rigid.

\item If $X= 13\cdot 13 \cdot 1 \in B_4$, then $\mathcal R(X)=2/3$,
since the left normal form of $(13\cdot 13 \cdot 1) \: 13 $ is $13
\cdot 13 \cdot 13 \cdot 1$, hence $k=2$ and $r=3$. This means that
two thirds of the left normal form of $X$ are preserved when
considering $X^2 = 13\cdot 13 \cdot 13 \cdot 13 \cdot 1 \cdot 1$.

\item If $X= 12132143 \cdot 143\in B_5$ (this is the example at
the beginning of Section~\ref{section:cyclings and powers}), then
$\mathcal R(X)=0$, since the left normal form of $(12132143 \cdot
143) \cdot 12132143 $ is $\Delta \cdot 2324321 \cdot 14$. Hence,
nothing from the left normal form of $X$ is preserved when computing
its square $X^2=\Delta \cdot 2324321 \cdot 14 \cdot 143$. In this
case we say that $X$ has no rigidity, or that it is 0-rigid. This
is, of course, the most difficult case if one tries to relate
cyclings and powers of $X$.

\ee

Let us see some characterizations of rigidity, and then how rigidity
behaves under cyclings or powers of an element.

\begin{lemma}\label{L:k-rigidity}
Let $X\in G$ with $\ell(X)=r>0$ and $\inf(X)=p$. Then $\mathcal
R(X)=k/r$  if and only if $k$ is the biggest integer such that
$$
(X^2\Delta^{-2p}) \wedge \Delta^k = (X \Delta^{-p}) \wedge \Delta^k.
$$
In particular, $\mathcal R(X)>0$ if and only if $\inf(X^2)= 2p$ and
$\iota(X^2)=\iota(X)$.
\end{lemma}

\begin{proof}
First notice that $k=0$ always satisfies the above condition, since
$(X^2\Delta^{-2p}) \wedge 1 = 1 = (X \Delta^{-p}) \wedge 1$. Also,
no $k>r$ can satisfy the condition, since $(X \Delta^{-p})$ only has
$r$ factors, and this would imply that $X^2\Delta^{-2p} = X
\Delta^{-p}$, which is not possible if $r>0$. Hence the biggest
integer $k$ satisfying the condition must belong to
$\{0,\ldots,r\}$.

The rigidity of $X=\Delta^p x_1\cdots x_r$ is at least $k/r$ if the
first $k$ factors in the left normal form of $x_1\cdots x_r
\tau^{-p}(x_1)$ are $x_1\cdots x_k$. This is the case if and only if
the biggest simple prefix of $x_k \cdots x_r \tau^{-p}(x_1)$ is
$x_k$, which in turn is the case if and only if the biggest simple
prefix of $x_k\cdots x_r \tau^{-p}(x_1\cdots x_r)$ is $x_k$. Since
$x_1\cdots x_k$ is in left normal form, the above condition holds if
and only if the left normal form of $X^2 = \Delta^{2p}
\tau^{p}(x_1\cdots x_r) (x_1\cdots x_r) $ has the form $\Delta^{2p}
\tau^p(x_1)\cdots \tau^{p}(x_k) z_1 \cdots z_t$ for some simple
elements $z_1,\ldots,z_t$. This happens if and only if
$(X^2\Delta^{-2p}) \wedge \Delta^k  = \tau^{-p}(x_1\cdots x_k) =
(X\Delta^{-p})\wedge \Delta^k$. This shows that $\mathcal R(X)=k/r$
if and only if $k$ is the biggest integer satisfying the latter
condition.

Now $\mathbf R(X)>0$ if $k$ is at least $1$, where for $k=1$ the
above condition reads $X^2\Delta^{-2p}\wedge \Delta =
\tau^{-p}(x_1)= \iota(X)$, which is equivalent to $\inf(X^2)=2p$ and
$\iota(X^2)=\iota(X)$. (Notice that if $\inf(X^2)>2p$ then
$X^2\Delta^{-2p}\wedge \Delta =\Delta \neq \iota(X^2)$.)
\end{proof}

One can also check the rigidity of an element by looking at its
inverse.

\begin{lemma}\label{L:rigidity and inverses}
Let $X\in G$  with $\ell(X)=r>0$. Then $\mathcal R(X)=k/r>0$ if and
only if $\sup(X^{-2})=2\: \sup(X^{-1})$ and the final $k$ factors in
the left normal forms of $X^{-1}$ and $X^{-2}$ coincide. In
particular, $X$ is rigid if and only if $X^{-1}$ is rigid. And also
$\mathcal R(X)>0$ if and only if $\sup(X^{-2})=2\:\sup(X^{-1})$ and
$\varphi(X^{-2})=\varphi(X^{-1})$.
\end{lemma}

\begin{proof}
The rigid case can be shown independently. We know by
Proposition~\ref{P:definitions of rigidity} that $X$ is rigid if and
only if $\iota(X)\wedge \iota(X^{-1})=1$, and this condition is
invariant under taking inverses. Hence $X$ is rigid if and only if
so is $X^{-1}$.

On the other hand, let $\Delta^p x_1\cdots x_r $ be the left normal
form of $X$. By Lemma~\ref{L:k-rigidity}, $\mathcal R(X)=k/r>0$ if
and only if the first $k$ factors in the left normal forms of
$X^2\Delta^{-2p}$ and $X\Delta^{-p}$ coincide. This means that the
left normal form of $X^2$ is $\Delta^{q} y_1\cdots y_t$, where
$q=2p$ and $y_i=\tau^p(x_i)$ for $i=1,\ldots,k$. By
Theorem~\ref{T:left normal form of the inverse}, the left normal
form of $X^{-1}$ is $\Delta^{-p-r} x_r'\cdots x_1'$ where
$x_i'=\tau^{-p-i}(\partial(x_i))$, and the left normal form of
$X^{-2}$ is $\Delta^{-q-t} y_t'\cdots y_1'$, where
$y_i'=\tau^{-q-i}(y_i)$. Then $q=2p$ means
$\sup(X^{-2})=2\sup(X^{-1})$, and   $y_i=\tau^{p}(x_i)$ means
$y_i'=\tau^{-2p-i}(\partial (y_i)) = \tau^{-2p-i}(\partial
(\tau^{p}(x_i)))=\tau^{-2p-i}(\tau^p(\partial(x_i)))=
\tau^{-p-i}(\partial(x_i))= x_i'$, hence the result follows.
\end{proof}

In the case of nonzero rigidity, we will be able to state some
common property of all powers of $X$.

\begin{corollary}\label{C:rigidity>0 and powers}
Let $X\in G$  with $\ell(X)>0$. If $\mathcal R(X)>0$, then
$\iota(X)=\iota(X^m)$  (whence $\varphi(X^{-1})=\varphi(X^{-m})$),
and also $\inf(X^m)=m\inf(X)$ (whence $\sup(X^{-m})=m\sup(X^{-1})$)
for every $m\geq 1$.
\end{corollary}

\begin{proof}
We just need to show the equalities $\iota(X)=\iota(X^m)$ and
$\inf(X^m)=m\inf(X)$, since in that case, by Corollary~\ref{C:inf,
sup and length of the inverse} one has $\varphi(X^{-1})=
\partial^{-1} (\iota(X)) =
\partial^{-1}(\iota(X^m))= \varphi(X^{-m})$ and
$\sup(X^{-m})=-\inf(X^m)= -m\inf(X)= m\sup(X^{-1})$.

The result is trivially true for $m=1$. Let $\Delta^p x_1\cdots x_r$
be the left normal form of $X$, and suppose that
$\iota(X^m)=\iota(X)=\tau^{-p}(x_1)$ and $inf(X^m)=m\inf(X)=mp$ for
some $m\geq 1$. Write then $X^m= \Delta^{mp} \tau^{(mp-p}(x_1)
y_2\cdots y_s$ in left normal form.

Notice that $\Delta\not\prec x_1\cdots x_r \tau^{-p}(x_1)$ as
$\mathcal R(X)>0$, hence there is no unexpected appearance of
$\Delta$ in the product $X\: X^m$. This implies $\inf(X^{m+1})=
\inf(X)+\inf(X^m)=p+mp = (m+1)\inf(X)$. Also $\iota(X^{m+1}) =
\iota(\Delta^p x_1\cdots x_r \: \tau^{-p}(x_1) \tau^{-mp}(y_2 \cdots
y_s) \Delta^{mp}) = \iota(\Delta^p x_1\cdots x_r \tau^{-p}(x_1)) =
\tau^{-p}(x_1) =\iota(X)$, and the result follows.
\end{proof}

Let us see that rigidity cannot decrease by cyclings.

\begin{proposition}\label{P:rigidity and cyclings}
Let $X\in SSS(X)$ with $\ell(X)>0$. Then $\mathcal R(X)\leq \mathcal
R(\mathbf c^t(X))$ for all $t\geq 1$. Furthermore, if $X\in USS(X)$,
equality holds.
\end{proposition}

\begin{proof}
Let $p=\inf(X)$. By definition of rigidity, $\mathcal R(X)=k/r$
means that $k$ is the biggest integer such that $X\iota(X)\wedge
\Delta^{p+k} = X \wedge \Delta^{p+k}$. If we apply the transport map
defined in~\cite{Gebhardt} to this equality, we get
$$
   \mathbf c(X) \iota(\mathbf c(X)) \wedge \Delta^{p+k} = \mathbf
   c(X)\wedge \Delta^{p+k},
$$
which means that $\mathcal R(\mathbf c(X))\geq k/r$ (notice that
$\ell(\mathbf c(X))=r$ since $X\in SSS(X)$). Applying the same
reasoning to every cycling of $X$, one has $\mathcal R(\mathbf
c^{t-1}(X))\leq \mathcal R(\mathbf c^{t}(X))$ for every $t\geq 1$,
so the result follows.

If $X\in USS(X)$, one cannot have $\mathcal R(X)< \mathcal R(\mathbf
c^t(X))$ for some $t$, since some further cycling of $\mathbf
c^t(X)$ would be equal to $X$, yielding the contradiction $\mathcal
R(X)< \mathcal R(X)$.
\end{proof}

If an element has some rigidity, its conjugating elements for
cycling $C_i$ satisfy the following useful property.

\begin{lemma}\label{L:1-rigid implies C_m rigid}
Let $X\in SSS(X)$ with $\ell(X)>0$. If $\mathcal R(X)>0$, the left
normal form of $\mathbf C_m$ is precisely $C_1 C_2\cdots C_m$ for
every $m\geq 1$.
\end{lemma}

\begin{proof}
Let $\Delta^p x_1\cdots x_r$ be the left normal form of $X$, and let
$\Delta^p y_1\cdots y_r$ be the left normal form of $\mathbf c(X)=
\Delta^p x_2\cdots x_r \tau^{-p}(x_1)$. Since $\mathcal R(X)>0$, the
left normal form of $x_1 x_2\cdots x_r \tau^{-p}(x_1)$ is
$x_1y_1\cdots y_r$. Since $C_1=\tau^{-p}(x_1)$, $C_2=\tau^{-p}(y_1)$
and $x_1y_1$ is left weighted, it follows that $C_1C_2$ is left
weighted. Applying the same reasoning to $\mathbf c^{i-1}(X)$ for
every $i\geq 2$, it follows that $C_i C_{i+1}$ is left weighted as
written, hence $C_1\cdots C_m$ is the left normal form of $\mathbf
C_m$.
\end{proof}

It is easy to see that if an element $X$ is rigid, then every power
of $X$ is rigid. We can generalize this to every element $X\in
USS(X)$, showing that the rigidity of $X$ can never decrease by
taking powers.

\begin{proposition}\label{P:rigidity and powers}
Let $X\in USS(X)$ with $\ell(X)\geq 1$.  Then $\mathcal R(X)\leq
\mathcal R(X^m)$ for every $m>1$.
\end{proposition}

\begin{proof}
Let $\ell(X)=r$ and $\mathcal R(X)=k/r$. If $k=0$ the result is
trivial, so we can assume that $k>0$. In the case $r=1$ this would
mean that $X$ is rigid, hence every power of $X$ is rigid and the
result would also be true. Therefore we will also assume that
$\ell(X)=r>1$. Let $\Delta^p x_1\cdots x_r$ be the left normal form
of $X$. Recall that by Corollary~\ref{C:rigidity>0 and powers}
$\inf(X^t)=tp$ for every $t\geq 1$. We will show that for $m\geq 1$
one has
$$
 (X^m\Delta^{-pm})\wedge \Delta^{mk}= \mathbf C_{mk}.
$$
Recall from Theorem~\ref{T:normal form of $C_m$ and $X^m$} that $
 (X^{mk}\Delta^{-pmk})\wedge \Delta^{mk}= \mathbf C_{mk}.
$ Since one has $X^m\Delta^{-pm} \preceq X^{mk}\Delta^{-pmk}$ for
every $k>0$, (notice that the infimum of both elements is 0) it
follows that
$$
(X^m\Delta^{-pm})\wedge \Delta^{mk} \preceq
(X^{mk}\Delta^{-pmk})\wedge \Delta^{mk}= \mathbf C_{mk},
$$
hence we only need to show that $\mathbf C_{mk} \preceq X^m
\Delta^{-pm}$.

We will first show that $C_i=\tau^{-p}(x_i)$ for $i=1,\ldots,k$.
This will be done by proving that, for $i=0,\ldots,k-1$ the first
$k-i$ non-$\Delta$ factors in the left normal form of $\mathbf
c^i(X)$ are $x_{i+1}\cdots x_k$. Indeed, this is trivially true for
$i=0$. If we assume the claim true for some $i$, $0\leq i <k-1$, we
have $\mathbf c^{i}(X)=\Delta^p x_{i+1}\cdots x_k y_{k+1}\cdots
y_r$. By Proposition~\ref{P:rigidity and cyclings} we know that
$\mathcal R(\mathbf c^i(X))\geq k/r$, hence the first $k-1$
non-$\Delta$ factors in the left normal form of $\mathbf c^{i+1}(X)
= \Delta^p x_{i+2}\cdots x_k y_1\cdots y_r \tau^{-p}(x_{i+1})$ are
precisely $x_{i+2} \cdots x_k y_1\cdots y_i$. In particular, the
first $k-(i+1)$ non-$\Delta$ factors are $x_{i+2}\cdots x_k$, thus
the claim is shown.  This implies that
$$
   \mathbf C_k = C_1\cdots C_k = \tau^{-p}(x_1\cdots x_k),
$$
hence $\mathbf C_k \preceq X\Delta^{-p}$.

Now suppose that $\mathbf C_{mk} \preceq X^m \Delta^{-pm}$ for some
$m\geq 1$. If we apply $k$ times the transport defined
in~\cite{Gebhardt}, which preserves $\preceq$, we obtain
$$
   \mathbf C_{[k,mk]}\preceq (\mathbf c^k(X))^m \Delta^{-pm}.
$$
But since $\mathbf C_k=\tau^{-p}(x_1\cdots x_k)$, it follows that $
\mathbf c^{k}(X)=X^{\mathbf C_k} = \Delta^p x_{k+1}\cdots x_r
\tau^{-p}(x_1\cdots x_k). $ Hence
$$
  X^{m+1} = (\Delta^p x_1\cdots x_r)^{m+1} =
  \Delta^p x_1\cdots x_k \;\left(x_{k+1} \cdots x_r \Delta^p x_1\cdots
  x_k\right)^m \; x_{k+1}\cdots x_r.
$$
$$
=\Delta^p x_1\cdots x_k \;(\Delta^p \tau^{p}(x_{k+1} \cdots x_r)
x_1\cdots x_k)^m \;x_{k+1}\cdots x_r
$$
$$
  = \tau^{-p}(x_1\cdots x_k) \; (\mathbf c^k(X))^m \;
  \tau^{-p}(x_{k+1}\cdots  x_r) \Delta^p.
$$
Since $\tau^{-p}(x_1\cdots x_k)=\mathbf C_k$, and $\mathbf
C_{[k,mk]}\preceq (\mathbf c^k(X))^m \Delta^{-pm}$, it follows that
$$
   \mathbf C_{(m+1)k} = \mathbf C_k \mathbf C_{[k,mk]} \preceq
   X^{m+1} \Delta^{-(m+1)p},
$$
as we wanted to show. Hence $
 (X^m\Delta^{-pm})\wedge \Delta^{mk}= \mathbf C_{mk}.
$ for every $m\geq 1$.

Now recall from Lemma~\ref{L:1-rigid implies C_m rigid} that, since
$k>0$, the left normal form of $\mathbf C_{mk}$ is $C_1\cdots
C_{mk}$. Hence, for every $m\geq 1$ one has
$$
   (X^{2m}\Delta^{-2pm})\wedge \Delta^{mk} = \mathbf C_{2mk}\wedge
   \Delta^{mk} = \mathbf C_{mk} =(X^m \Delta^{-pm}) \wedge
   \Delta^{mk}.
$$
By Lemma~\ref{L:k-rigidity}, and since $\ell(X^m)\leq mr$, this
implies that $\mathcal R(X^m) \geq \frac{mk}{mr} = \frac{k}{r} =
\mathcal R(X)$, as we wanted to show.
\end{proof}

{\bf Remark:} The main difference between rigidity of cyclings and
rigidity of powers is that, while iterated cycling of $X\in SSS(X)$
yields a non-decreasing sequence
$$
  \mathcal R(X) \leq \mathcal R(\mathbf c(X)) \leq \mathcal R(\mathbf
  c^2(X)) \leq \cdots,
$$
this does not happen for powers of $X$, even if $X\in USS(X)$. For
instance, if $X= 12132143 \cdot 143\in B_5$ is the example at the
beginning of Section~\ref{section:cyclings and powers}, one has
$$
 \mathcal R(X)=0, \quad  \mathcal R(X^2)=0, \quad \mathcal R(X^3)=1,
\quad \mathcal R(X^4)=0, \quad  \mathcal R(X^5)=0, \quad \mathcal
 R(X^6)=1, \ldots
$$
Notice that this is not in contradiction with
Proposition~\ref{P:rigidity and powers}, where the rigidity of $X^m$
is compared with that of $X$, not with the rigidity of the
intermediate powers.

\ms

The above results imply that elements having some rigidity behave
nicely with respect to powers and cyclings, in the following sense:

\begin{corollary}\label{C:1-rigid implies cyc & pow commute}
If $X\in USS(X)$ with $\ell(X)\geq 1$ and $\mathcal R(X)>0$, then
one has $(\mathbf c^t(X))^m = \mathbf c^t(X^m)$ for every $t,m\geq
1$.
\end{corollary}

\begin{proof}
By Lemma~\ref{L:1-rigid implies C_m rigid}, the left normal form of
$\mathbf C_m$ is $C_1\cdots C_m$. If $\ell(X)>1$ this implies, by
Theorem~\ref{T:normal form of $C_m$ and $X^m$}, that
$C_1=\iota(C_1\cdots C_m)=\iota(X^m)$. If $\ell(X)=1$ then
$X^m=C_1\cdots C_m \Delta^{mp}$ where $p=\inf(X)$, so also in this
case we have $\iota(X^m)=C_1$. Hence $\mathbf c(X^m) = (X^m)^{C_1} =
(X^{C_1})^m =(\mathbf c(X))^m$, and the result is true for $t=1$. If
the result is true for some $t-1$, it suffices to apply the previous
case to $\mathbf c^{t-1}(X)$, which has some rigidity by
Proposition~\ref{P:rigidity and cyclings}, to obtain
$$
    \mathbf c^t(X^m) = \mathbf c (\mathbf c^{t-1}(X^m)) = \mathbf c
    ((\mathbf c^{t-1}(X))^m) = (\mathbf c (\mathbf c^{t-1}(X)))^m
    = (\mathbf c^t(X))^m.
$$
\end{proof}

\begin{corollary}\label{C:1-rigid in USS has powers in closed orbit}
Let $X\in USS(X)$ with $\ell(X)\geq 1$. If $\mathcal R(X)>0$ then
$X^m$ belongs to a closed orbit under cycling, for every $m\geq 1$.
\end{corollary}

\begin{proof}
Let $N$ be the orbit length of $X$. By Corollary~\ref{C:1-rigid
implies cyc & pow commute}, $\mathbf c^N(X^m) = (\mathbf c^N(X))^m =
X^m$, so the result follows.
\end{proof}

\ms

 {\bf Remark:} The above result does not imply that $X^m\in USS(X)$,
 since it could happen that $X^m\not\in SSS(X^m)$. But the fact that
 $X^m$ belongs to a closed orbit under cycling will be enough for
 our purposes.

\subsection{The ultra summit set of a rigid element is made of rigid
elements} \label{The ultra summit set of a rigid element is made of
rigid elements}

By the above discussion on rigidity, we know that if $X\in USS(X)$
is rigid, then the whole orbit of $X$ under cycling is made of rigid
braids. But what about the other orbits in $USS(X)$?  In this
subsection we will show that, if $\ell(X)>1$, all orbits in $USS(X)$
are made of rigid braids. Hence $USS(X)$ is just the set of rigid
conjugates of $X$.  We start with three small results.

\begin{lemma}\label{L:USS closed under decycling}
Given $X\in USS(X)$ then $\mathbf d(X)\in USS(X)$.
\end{lemma}

\begin{proof}
Let $\Delta^p x_1\cdots x_r$ be the left normal form of $X$. Notice
that $X^X=X\in USS(X)$ and that $X^{\Delta^{p+r-1}} =
\tau^{p+r-1}(X) \in USS(X)$. Then, by Theorem~\ref{T:convexity_USS},
$\mathbf d(X) = X^{(\Delta^p x_1\cdots x_{r-1})} = X^{X\wedge
\Delta^{p+r-1}} \in USS(X)$.
\end{proof}

\ms

\begin{lemma}\label{L:orbits are rigid}
If $X\in G$ is rigid, then $\mathbf c^i(X)$ and $\mathbf d^i(X)$ are
rigid for every $i\geq 1$. Moreover $X\in USS(X)$, and if $e\geq 1$
is such that $\Delta^e$ is central, then $\mathbf c^m(X)=X$ for some
$m\leq \ell(X)\: e$.
\end{lemma}

\begin{proof}
Let $\Delta^p x_1\cdots x_r$ be the left normal form of $X$. By
Proposition~\ref{P:rigidity and cyclings}, all iterated cyclings of
$X$ are rigid. Then one can easily show by recurrence that if
$i=kr+j$ with $0\leq j \leq r-1$, the left normal form of $\mathbf
c^i(X)$ is equal to $\Delta^p \tau^{kp}(x_{j+1})\cdots
\tau^{kp}(x_r) \tau^{(k+1)p}(x_1)\cdots \tau^{(k+1)p}(x_j)$. Hence,
if $e$ is such that $\tau^e=1$, one has $\mathbf c^{er}(X)= \Delta^p
\tau^{ep} (x_1)\cdots \tau^{ep}(x_r) = \Delta^p x_1 \cdots x_r=X$,
so $X$ belongs to a closed orbit under cycling, and the orbit length
is a divisor of $re=\ell(X)\: e$.

By Lemma~\ref{L:rigidity and inverses} $X$ is rigid if and only if
so is $X^{-1}$. This means that $\mathbf c^i(X^{-1})$ is rigid for
every $i\geq 1$. But we know by~\cite{El-M} that $\mathbf
c^i(X^{-1})=(\mathbf d^i(X))^{-1}$, hence $\mathbf d^i(X)$ is also
rigid for every $i\geq 1$.

Furthermore, by the above arguments $X^{-1}$ belongs to a closed
orbit under cycling, thus $X$ belongs to a closed orbit under
decycling. But an element belonging to closed orbits under cycling
and decycling belongs to its ultra summit set (since one can always
reach the ultra summit set by iterated cycling and decycling), so it
follows that $X\in USS(X)$.
\end{proof}

For elements which belong to their ultra summit set, the converse of
Lemma~\ref{L:orbits are rigid} is also true.

\begin{lemma}\label{L:nonrigid_c_d}
If $Y\in USS(X)$ is not rigid, then neither $\mathbf c^i(Y)$ nor
$\mathbf d^i(Y)$ are rigid for any $i\geq 1$.
\end{lemma}

\begin{proof}
It clearly suffices to show the result for $i=1$. First, $\mathbf
c(Y)$ has the same rigidity as $Y$ by Proposition~\ref{P:rigidity
and cyclings}, hence it cannot be rigid.

Now suppose that $\mathbf d(Y)$ is rigid. This is clearly not
possible if $\ell(Y)=1$, so we can suppose that $\ell(Y)>1$. If
$\Delta^p y_1\cdots y_r$ is the left normal form of $Y$, then
$\mathbf d(Y)= \Delta^p \tau^p(y_r)y_1\cdots y_{r-1}$, although this
decomposition is not the left normal form of $\mathbf d(Y)$.
Nevertheless, since $y_1\cdots y_{r-1}$ is in left normal form, we
know, by the left normal form algorithm (see for
instance~\cite{Epstein}) that there is a decomposition $y_i=a_ib_i$
for $i=1,\ldots,r-1$ such that the left normal form of $\mathbf
d(Y)$ is precisely $\Delta^p (\tau^p(y_r)a_1)(b_1a_2)\cdots
(b_{r-2}a_{r-1})(b_{r-1})$.

Since we are assuming that $\mathbf d(Y)$ is rigid, we have that
$(b_{r-1})(y_r\tau^{-p}(a_1))$ is left weighted as written. Notice
that this implies that $y_{r-1} (y_r \tau^{-p}(a_1))$ is left
weighted as written.  Hence the left normal form of $y_1\cdots y_r
\tau^{-p}(y_1)$ is precisely $y_1\cdots y_{r-1} (y_r \tau^{-p}(a_1))
\tau^{-p}(b_1)$. In other words, $\mathcal R(Y)= \frac{r-1}{r}$. By
Proposition~\ref{P:rigidity and cyclings}, all iterated cyclings of
$Y$ have rigidity $\frac{r-1}{r}$. This implies, in particular, that
$$
\mathbf c^{r-1}(Y) = Y^{\tau^{-p}(y_1\cdots y_{r-1})} = \Delta^p y_r
\tau^{-p}(y_1\cdots y_{r-1}) = \tau^{-p}(\mathbf d(Y)),
$$
but this latter element is supposed to be rigid. A contradiction.
Hence no iterated decycling of $Y$ can be rigid.
\end{proof}

We can finally prove the main result concerning the elements of the
ultra summit set of a rigid element.

\begin{theorem}
\label{T:X rigid implies every Y in USS(X) rigid} Let $X$ be rigid
and $\ell(X)>1$.  Then every element in $USS(X)$ is rigid.
\end{theorem}

\begin{proof}
Suppose that there exists an element in $USS(X)$ which is not rigid.
Since every two elements in $USS(X)$ are connected by a sequence of
conjugations by simple elements, there must be a non-rigid element
in $USS(X)$ which is the conjugate of a rigid one by a simple
element. Hence we can assume without loss of generality that $X^s=Y$
for some non-rigid element $Y\in USS(X)$ and some simple element
$s$. We will also assume that $s$ is a maximal element (with respect
to $\preceq$) in the set of all simple elements conjugating $X$ to a
non-rigid element in $USS(X)$. We will get a contradiction by
showing that $X$ is conjugate to $\mathbf c(\mathbf d(Y))$ by a
simple element which is a proper right multiple of $s$.

Let $\Delta^p x_1\cdots x_r$ be the left normal form of $X$ and let
$\Delta^p y_1\cdots y_r$ be the left normal form of $Y$. Since
$X^s=Y$, it is known by~\cite{Gebhardt} that there exist simple
elements $s_0,\ldots,s_r$ such that $s_0=\tau^p(s)$, $\;s_r=s$ and
$y_i=s_{i-1}^{-1}x_is_i$ for $i=1,\ldots,r$. That is, the left
normal form of $Y$ is $\Delta^p y_1\cdots y_r = \Delta^p (s_0^{-1}
x_1s_1)(s_1^{-1}x_2s_2)\cdots (s_{r-1}^{-1}x_rs_r)$.

Now consider $\mathbf d(X)$ and $\mathbf d(Y)$. By Lemma~\ref{L:USS
closed under decycling}, these two elements belong to $USS(X)$.
Since $X$ is rigid, the left normal form of $\mathbf d(X)$ is
$\Delta^p \tau^p(x_r)x_1\cdots x_{r-1}$. However, $Y$ is not rigid,
so the left normal form of $\mathbf d(Y)$ is not $\Delta^p
\tau^p(y_r) y_1\cdots y_{r-1}$, since $\tau^p(y_r)y_1$ is not left
weighted as written. (Here we use the fact that $r=\ell(X)>1$.)
Hence, $\iota(\mathbf d(Y))=y_r t$ for some nontrivial simple
element $t$.

The elements $\mathbf d(X)$ and $\mathbf d(Y)$ are also connected
through a conjugation by a simple element. Namely, $\mathbf
d(X)^{s_{r-1}}= \Delta^p \tau^p(s_{r-1})^{-1} \tau^p(x_r) x_1\cdots
x_{r-1}s_{r-1} = \Delta^p \tau^p(y_r) y_1\cdots y_{r-1} = \mathbf
d(Y)$. Hence, by~\cite{Gebhardt} again, there exist simple elements
$t_0,\ldots,t_r$ such that $t_0=\tau^p(s_{r-1})$, $\; t_r=s_{r-1}$
and the left normal form of $\mathbf d(Y)$ is
$\Delta^p(t_0^{-1}\tau^p(x_r)t_1)(t_1^{-1}x_1t_2) \cdots
(t_{r-1}x_{r-1}t_r)$. Therefore $\iota(\mathbf
d(Y))=\tau^{-p}(t_0)^{-1} x_r \tau^{-p}(t_1)=s_{r-1}^{-1}x_r
\tau^{-p}(t_1)$. Since we saw that $\iota(\mathbf d(Y))=y_r t =
(s_{r-1}^{-1} x_r s)t$, it follows that $st=\tau^{-p}(t_1)$, which
is a simple element. If we denote $u=\tau^{-p}(t_1)$, we just showed
that $s\prec u$ (strict) and that $\iota(\mathbf d(Y))=s_{r-1}^{-1}
x_r u$.

Finally, notice that $X^u = ((X^{x_r^{-1}})^{s_{r-1}})^{s_{r-1}^{-1}
x_r u} = (\mathbf d(X)^{s_{r-1}})^{s_{r-1}^{-1} x_r u} = \mathbf
d(Y)^{\iota(\mathbf d(Y))} = \mathbf c(\mathbf d(Y))$. But since $Y$
is not rigid and belongs to $USS(X)$, Lemma~\ref{L:nonrigid_c_d}
tells us that $\mathbf d(Y)$ is not rigid. Since $\mathbf d(Y)$ also
belongs to $USS(X)$ by Lemma~\ref{L:USS closed under decycling}, it
follows again by Lemma~\ref{L:nonrigid_c_d} that $\mathbf c(\mathbf
d(Y))$ is not rigid, and belongs to $USS(X)$. But $s\precneqq u$, so
this contradicts the maximality of $s$, and we are done.

\end{proof}

\begin{corollary}\label{C:USS_of_rigid}
If $X$ is rigid and $\ell(X)>1$, then $USS(X)$ is the set of rigid
conjugates of $X$.
\end{corollary}

\begin{proof}
Let $Y$ be a rigid conjugate of $X$. Since $Y$ is rigid, it belongs
to its ultra summit set and since it is conjugate to $X$, its ultra
summit set is precisely $USS(X)$. Conversely, every element in
$USS(X)$ is rigid by the above result.
\end{proof}

\begin{corollary}
\label{C:USS_inverse} If $X$ is rigid and $\ell(X)>1$, then
$USS(X^{-1})$ is the set of inverses of the elements in $USS(X)$.
\end{corollary}

\begin{proof}
This is a direct consequence of Corollary~\ref{C:USS_of_rigid} and
the fact that $Y\in G$ is rigid if and only if $Y^{-1}$ is rigid
(Lemma~\ref{L:rigidity and inverses}).
\end{proof}

{\bf Remark:}  If $\ell(X)=1$, then $USS(X)$ may contain rigid and
non-rigid elements. For instance, the simple element $1 2 3 2 1 4 3
5  \in B_6$ is rigid (since $1 2 3 2 1 4 3 5  \cdot 1 2 3 2 1 4 3 5
$ is left weighted), but it is conjugate (by $23$) to the simple
element $1 2 1 3 4 3 2 5 $, which is not rigid (the left normal form
of $1 2 1 3 4 3 2 5 \; 1 2 1 3 4 3 2 5 $ is $1 2 1 3 4 3 2 5 1 4
\cdot 2 1 3 2 4 5$). Clearly both elements belong to the ultra
summit set, since they are simple.

\subsection{Elements having a rigid power}

In this section we will characterize elements $X\in G$ having a
rigid power $X^m$ for some integer $m\neq 0$. Notice that such an
element cannot be periodic. Otherwise, since rigidity is preserved
by powers, some rigid power of $X$ (which, by definition, has
positive canonical length) would also be a power of $\Delta$ (which
has zero canonical length), and this is not possible. If the element
$X$ belongs to its ultra summit set, we can say something more.

\begin{proposition}\label{P:pre-rigid implies C_M is power of X}
Let $X\in USS(X)$ with $\ell(X)\geq 1$. If $X$ has a rigid power,
then there exists some $M>0$ such that $\mathbf C_M = \Delta^{k}
X^t$ for some integers $k,t$, where $t>0$ and $\Delta^k$ is central.
\end{proposition}

\begin{proof}
Let $p=\inf(X)$. If $\ell(X)=1$ then $X^M= \mathbf C_M \Delta^{-pM}$
for every $M$, so we just need to take $M$ big enough so that
$\Delta^M$ is central, and we are done. Hence we can assume that
$\ell(X)>1$.

Let $e>0$ be such that $\Delta^e$ is central, let $m>0$ be such that
$X^m$ is rigid, and let $N$ be the orbit length of $X$ under
cycling. Consider $T=emN$. By Lemma~\ref{L:decomposition} one has
$X^{T}=\mathbf C_{T} \mathbf R_{T}\Delta^{pT}$, and by
Proposition~\ref{P:CR left weighted}, $\varphi(\mathbf C_T)
\iota(\mathbf R_T)$ is left weighted. Since $T$ is a multiple of
$e$, $\Delta^{pT}$ is central. Since $T$ is a multiple of $m$, it
follows that $X^T$ is rigid, hence $\varphi(X^T)\iota(X^T) =
\varphi(\mathbf R_T) \iota(\mathbf C_T)$ is also left weighted
(notice that the equality holds since $\Delta^{pT}$ is central).
Finally, since $T$ is a multiple of $N$, it follows that $\mathbf
C_T$ commutes with $X$, thus it commutes with $X^T$. Then one has:
$$
   X^T = (X^T)^{\mathbf C_T} = \mathbf R_T \Delta^{pT} \mathbf C_T =
   \Delta^{pT} \mathbf R_T \mathbf C_T.
$$
Moreover,  since $\varphi(\mathbf R_T) \iota(\mathbf C_T)$ is left
weighted, $\varphi(X^T) = \varphi(\mathbf C_T)$. Hence
$\varphi(\mathbf C_T)\iota(\mathbf C_T)=\varphi(X^T)\iota(X^T)$ is
left weighted by the rigidity of $X^T$, so it follows that $\mathbf
C_T$ is also rigid. In particular, $\inf((\mathbf
C_T)^k)=k\inf(\mathbf C_T)$ for every $k>0$. Since $T$ is a multiple
of $N$, one has $(\mathbf C_T)^k= \mathbf C_{Tk}$.  Hence, by
considering a suitable multiple of $T$, we can assume that
$\inf(\mathbf C_T)$ is a multiple of $e$, that is, $\mathbf
C_T=\Delta^{eq} y_1\cdots y_r$, and $\mathbf R_T=z_1\cdots z_s$,
where $y_r z_1$ and $z_sy_1$ are left weighted. Then one has
$$
   \Delta^{e(pmN+q)} y_1\cdots y_r z_1\cdots z_s = X^T =
   \Delta^{e(pmN+q)} z_1\cdots z_s y_1\cdots y_r,
$$
where both decompositions of $X^T$ are in left normal form. In other
words, the left normal form of $X^T$ is invariant under some cyclic
permutations of its factors. This is only possible if there is some
rigid element $Y=y_1 \cdots y_a$ (where $a=\gcd (r,s)$), such that
$Y^{i}=y_1\cdots y_r$ and $Y^{j}=y_1\cdots y_r z_1\cdots z_s$ for
some $i,j>0$.  But then $\mathbf C_{Tj} = (\mathbf C_T)^{j} =
(\Delta^{e(pmN+q)} y_1\cdots y_r )^j = \Delta^{ej(pmN+q)} Y^{ij}$.
Since $Y^{ij}= (y_1\cdots y_r z_1\cdots z_s)^i$, it follows that
$Y^{ij}=\Delta^{ek'} X^{Ti}  $ for some $k'\in \mathbb Z$. Denoting
$M=Tj$,  $t=Ti$ and $k=e(k'+jpmN+jq)$, one finally obtains $\mathbf
C_{M}=\Delta^{k} X^{t}$, as we wanted to show.
\end{proof}

We will now show that the converse of Proposition~\ref{P:pre-rigid
implies C_M is power of X} is also true for elements of canonical
length greater than 1, by the following two results.

\begin{proposition}\label{P:some power is 1-rigid}
Let $X\in USS(X)$  with $\ell(X)>1$. Suppose that $\mathbf C_M =
\Delta^{k} X^t$ for some integers $M,k,t$, where $M,t>0$ and
$\Delta^k$ is central.  Then $\mathcal R(X^{-m})>0$ for some $m>0$.
\end{proposition}

\begin{proof} Notice that $\mathbf C_M$ commutes with $X$, hence $M$
is a multiple of the orbit length of $X$, and then $(\mathbf C_M)^r=
\mathbf C_{Mr}$ for every $r\geq 1$. This implies $\mathbf
C_{2M}=(\mathbf C_M)^2 = \Delta^{2k} (X^t)^2$, where $\Delta^{2k}$
is also central. Hence, replacing $M$ by a multiple if necessary, we
can assume that $M\geq ||\Delta||$.

By Lemma~\ref{L:decomposition}, $\sup(\mathbf C_M)=M$ and
$\sup(\mathbf C_{2M})=2M$. Hence $\sup((X^t)^2)=-2k+2M = 2 (-k+M) =
2 \sup(X^t)$.  At the same time, since $M$ is a multiple of the
orbit length of $X$ and $M\geq ||\Delta||$, one has $\varphi(\mathbf
C_M)=F(\mathbf c^{M}(X))= F(X)$ and also $\varphi(\mathbf
C_{2M})=F(X)$. Therefore $\varphi((X^t)^2)=\varphi(\mathbf C_{2M})=
F(X) = \varphi(\mathbf C_M) = \varphi(X^t)$. By
Lemma~\ref{L:rigidity and inverses}, this means that $\mathcal
R(X^{-t})>0$, so we take $m=t$ and we are done.
\end{proof}

\begin{proposition}\label{P:some power is rigid}
Let $X\in USS(X)$  with $\ell(X)>1$. Suppose that $\mathbf C_M =
\Delta^{k} X^t$ for some integers $M,k,t$, where $M,t>0$ and
$\Delta^k$ is central.  Then $X^T$ is rigid for some $T>0$.
\end{proposition}

\begin{proof}
We know by Proposition~\ref{P:some power is 1-rigid} that $\mathcal
R(X^{-m})>0$ for some $m>0$. We also know that $\mathbf C_{Mr}=
\Delta^{kr} X^{tr}$ for every $r\geq 1$. Hence, replacing $M$ by a
multiple, if necessary, we can assume that both $M$ and $t$ are
multiples of $m$. Since $\mathcal R(X^{-m})>0$, and $M$ and $t$ are
multiples of $m$, Corollary~\ref{C:rigidity>0 and powers} implies
that $\varphi(X^M)=\varphi(X^{m})=\varphi(X^t)= \varphi(\mathbf
C_M)$. Notice that $M-t\neq 0$, otherwise $\mathbf R_M$ would be a
power of $\Delta$, while $\ell(\mathbf R_M)>0$ by
Lemma~\ref{L:decomposition}. Hence $|M-t|$ is also a nontrivial
multiple of $m$, so we have $\varphi(\mathbf
C_M)=\varphi(X^{|M-t|})$.

On the other hand, by Lemma~\ref{L:decomposition} one has $X^M =
\mathbf C_M \mathbf R_M \Delta^{pM} = \Delta^k X^t \mathbf R_M
\Delta^{pM}$, where $p=\inf(X)$. Since $\Delta^k$ is central, this
means $\mathbf R_M = X^{-t} \Delta^{-k} X^M \Delta^{-pM} = X^{M-t}
\Delta^{-k-Mp}$. Recall that $\ell(\mathbf R_M)>0$, hence
$\iota(\mathbf R_M)=\iota(X^{M-t})$.  Moreover $M-t>0$, otherwise we
would have $\varphi(\mathbf C_M)=\varphi(X^{t-M})$ and
$\iota(\mathbf R_M)=\iota(X^{M-t})=\partial(\varphi(X^{t-M}))$,
whence $\varphi(\mathbf C_M)\iota(\mathbf R_M)=\Delta$ and this
contradicts Proposition~\ref{P:CR left weighted}, which states that
$\varphi(\mathbf C_M)\iota(\mathbf R_M)$ is left weighted.

Therefore $M-t>0$, and then $\varphi(\mathbf C_M)=\varphi(X^{M-t})$.
Therefore $\varphi(X^{M-t})\iota(X^{M-t})= \varphi(\mathbf C_M)
\iota(\mathbf R_M)$ is left weighted, hence $X^{M-t}$ is rigid and
we just take $T=M-t$.
\end{proof}

We have then shown the following result.

\begin{theorem}\label{T:characterization of pre-rigid}
Let $X\in USS(X)$ with $\ell(X)>1$. Then $X$ has a rigid power if
and only if $\mathbf C_M = \Delta^{k} X^t$ for some integers
$M,k,t$, where $M,t>0$ and $\Delta^k$ is central.

Moreover, in this case $M-t>0$ and $M$ can be chosen so that
$\mathbf R_M = \Delta^s X^{M-t}$ where $\Delta^s$ is central, and
$X^M$, $X^t$ and $X^{M-t}$ are all rigid.
\end{theorem}

\begin{proof}
The first claim is shown in Propositions~\ref{P:pre-rigid implies
C_M is power of X} and \ref{P:some power is rigid}.  In the proof of
Proposition~\ref{P:some power is rigid} it is also shown that in
this case $M-t>0$. Replacing $M$ by some suitable multiples, we
replace $X^M$, $X^t$ and $X^{M-t}$ by powers, hence we can choose
$M$ in such a way that these three elements are rigid.
\end{proof}

A very interesting consequence of this result is the following

\begin{theorem}\label{T:pre-rigid for whole USS}
Let $X \in USS(X)$ with $\ell(X)\geq 1$. If $X$ has a rigid power,
then all elements in $USS(X)$ have rigid powers.
\end{theorem}

\begin{proof}
By Theorem~\ref{T:characterization of pre-rigid}, $\mathbf C_M=
\Delta^k X^t$ for some integers $M,k,t$, where $M,t>0$ and
$\Delta^k$ is central. We can also assume that $M$ is a multiple of
the orbit length of $X$.

Let $Y\in USS(X)$. For $i\geq 1$, let $C_i'$, $R_i'$, $\mathbf C_i'$
and $\mathbf R_i'$ denote the elements analogous to $C_i$, $R_i$,
$\mathbf C_i$ and $\mathbf R_i$, defined for $Y$ instead of $X$.

Let $\alpha$ be a positive element such that $\alpha^{-1} X \alpha
=Y$. In~\cite{Gebhardt}, the $M$-th transport of $\alpha$ is defined
as the element $\beta$ such that $\alpha^{-1} \mathbf C_M \beta =
\mathbf C_M'$. It is shown in~\cite{Gebhardt} that some iterated
transport of $\alpha$ will be equal to $\alpha$. Hence, replacing
$M$ by a multiple if necessary, we can assume that $\beta=\alpha$.
But then $\mathbf C_M' = \alpha^{-1} \mathbf C_M \alpha =
\alpha^{-1} \Delta^k X^t \alpha = \Delta^k \alpha^{-1} X^t \alpha =
\Delta^k Y^t$. By Theorem~\ref{T:characterization of pre-rigid} this
means that $Y$ also has a rigid power.
\end{proof}

\subsection{Consequences for pseudo-Anosov braids.}
\label{subsection:Consequences for pseudo-Anosov braids}

The results from the previous subsection have a very important
consequence in the case of braid groups. The structure of
centralizers of pseudo-Anosov braids is well known, and this allows
to show that pseudo-Anosov braids in their ultra summit set have
rigid powers.

\begin{theorem}\label{T:pA and USS implies pre-rigid}
Let $X\in B_n$ be a pseudo-Anosov braid. If $X\in USS(X)$ and
$\ell(X)>1$, then $X$ has a rigid power.
\end{theorem}

\begin{proof}
Let $N$ be the orbit length of $X$ under cycling. Then $\mathbf C_N$
commutes with $X$. It is known~\cite{Ivanov,GM-W} that if $X$ is
pseudo-Anosov, every element in the centralizer of $X$ has a common
power with $X$, up to multiplication by a central power of $\Delta$.
Hence $(\mathbf C_N)^r=\Delta^k X^t$ for some integers $r$, $k$,
$t$, where $\Delta^k$ is central. Moreover we can assume that $r>0$,
otherwise we consider the inverse of the above equation. Since $N$
is the orbit length of $X$, $(\mathbf C_N)^r = C_{Nr}$, hence taking
$M=Nr$ one has $\mathbf C_M=\Delta^k X^t$ for some positive $M$ and
some integers $k,t$ such that $\Delta^k$ is central. By
Theorem~\ref{T:characterization of pre-rigid}, we only need to show
that $t>0$.

Suppose that $t<0$. Replacing $M$ by a multiple if necessary, we can
assume that $M\geq ||\Delta||$, hence $\varphi(\mathbf
C_M)=F(X)=\varphi(\mathbf C_{Mr})= \varphi(X^{tr})$ for every $r\geq
1$. We can also assume that $M$ is a (positive) multiple of $t$,
hence $t-M$ will be a (negative) multiple of $t$, and then
$\varphi(\mathbf C_M)=\varphi(X^{t-M})$.

On the other hand, by Lemma~\ref{L:decomposition}, we know that
$X^M= \mathbf C_M \mathbf R_M \Delta^{Mp}$, where $p=\inf(X)$. Hence
$\mathbf R_M = X^{M-t} \Delta^{-k-Mp}$, and then $\iota(\mathbf
R_M)=\iota(X^{M-t})=\partial(\varphi(X^{t-M}))$. This would
contradict Proposition~\ref{P:CR left weighted}, which states that
$\varphi(\mathbf C_M)\iota(\mathbf R_M)$ is left weighted. Therefore
$t>0$, and Theorem~\ref{T:characterization of pre-rigid} implies
that $X$ has a rigid power.
\end{proof}

\begin{corollary}\label{C:pA has a rigid power}
Every pseudo-Anosov braid has a rigid power, up to conjugacy.
\end{corollary}

\begin{proof}
Let $Y$ be a pseudo-Anosov braid, and let $X\in SU(Y)$. That is, $X$
is conjugate to $Y$ and all powers of $X$ belong to their ultra
summit set.

We know that powers of $\Delta$ are not pseudo-Anosov but periodic,
hence $\ell(X)\geq 1$. We will first show that we have $\ell(X^m)>1$
for some $m>1$. Indeed, if $\ell(X^m)=1$ for all $m>1$, since the
set of simple elements is finite we would have $X^a=\Delta^u s$ and
$X^b=\Delta^v s$ for the same simple element $s$ and $a\neq b$. But
then $X^{b-a}= \Delta^{v-u}$, which is not possible since a
pseudo-Anosov braid cannot be periodic.

Since the property of being pseudo-Anosov is preserved by powers,
$X^m$ is a pseudo-Anosov braid such that $\ell(X^m)>1$. Moreover,
$X^m\in USS(X^m)$, as $X\in SU(Y)$. Hence we can apply
Theorem~\ref{T:pA and USS implies pre-rigid} to $X^m$ and it follows
that some power of $X^m$, thus some power of $X$, is rigid. Since
$X$ is conjugate to $Y$, the result follows.
\end{proof}

We remark that generic elements of $B_n$ are pseudo-Anosov. This
means that most elements in $B_n$ have rigid powers, up to
conjugacy.

\subsection{A bound for the rigid power of an element.}
\label{subsection:A bound for the rigid power of an element}

In this section we will show that, if $X\in USS(X)$ has a rigid
power, and several powers of $X$ belong to their ultra summit sets,
then $X^m$ is rigid for some {\it small} $m$, namely $m<
||\Delta||^3$. Moreover, if $\ell(X)>1$ we can take
$m<||\Delta||^2$. In the particular case of braid groups, using the
Artin structure one has $||\Delta||=n(n-1)/2$, and using the
Birman-Ko-Lee structure $||\Delta||=n-1$. Hence in both cases the
bound is polynomial on the number of strands, and does not depend on
the length of the braid.

We first need to show two results concerning elements having rigid
powers and absolute final factors.

\begin{proposition}\label{P:F(X)=F(X^t) if X is PA}
Let $X\in USS(X)$ with $\ell(X)>1$. Suppose that $X$ has a rigid
power, and that $X^t\in USS(X^t)$ for some $t>1$. Then
$F(X)=F(X^t)$.
\end{proposition}

\begin{proof}
By Theorem~\ref{T:characterization of pre-rigid}, one has $\mathbf
C_M=\Delta^k X^r$ for some $M,k,r$ such that $M,r>0$ and $\Delta^k$
is central. Replacing $M$ by a multiple if necessary, so that $M$ is
a multiple of the orbit length of $X$ under cycling, one has
$\mathbf C_{[-M,M]} = \mathbf C_M =\Delta^k X^r$.

In the same way, $X^t$ also has a rigid power. Moreover, since
$\ell(\mathbf C_t)\geq 1$, $\ell(\mathbf R_t)\geq 1$ and $\mathbf
C_t \mathbf R_t$ is left weighted, one has $\ell(X^t)>1$. Hence we
can apply Theorem~\ref{T:characterization of pre-rigid} to $X^t$ and
we obtain $\mathbf C_{[-M',M']}^{(t)} = \mathbf C_{M'}^{(t)} =
\Delta^{k'} X^{tr'}$ for some $M',r'>0$ and some $k'$ such that
$\Delta^{k'}$ is central.

Replacing $M$ and $M'$ above by some suitable multiples, we can
assume that $r=tr'$, and also that $M,M'\geq ||\Delta||$. Hence
$F(X^t)=\varphi(\mathbf C_{[-M',M']}) = \varphi(\Delta^{k'} X^{tr'})
= \varphi(\Delta^k X^r) = \varphi(\mathbf C_{[-M,M]}) = F(X)$.
\end{proof}

\begin{proposition}\label{P:1-rigid PA and absolute factors}
Let $X\in USS(X)$ with $\ell(X)>1$, and suppose that $X$ has a rigid
power. If $\mathcal R(X)>0$, then $\iota(X)=I(X)$. If $\mathcal
R(X^{-1})>0$, then $\varphi(X)=F(X)$.
\end{proposition}

\begin{proof}
By Theorem~\ref{T:characterization of pre-rigid}, there exists some
$M>0$ such that $\mathbf C_M = \Delta^k X^t$ and $\mathbf
R_M=\Delta^s X^{M-t}$, where $t>0$, $M-t>0$, and both $\Delta^k$ and
$\Delta^s$ are central.

Suppose that $\mathcal R(X^{-1})>0$. By Corollary~\ref{C:rigidity>0
and powers} this implies that $\varphi(X^m)=\varphi(X)$ for every
$m\geq 1$. We can assume that $M$ is a multiple of the orbit length
of $X$, hence $\mathbf C_{[-M,M]}=\mathbf C_M$. If one chooses $M$
big enough (replacing it by a multiple if necessary), one has
$F(X)=\varphi(\mathbf C_{[-M,M]})= \varphi(\mathbf C_M)=
\varphi(X^t) = \varphi(X)$.

Now suppose that $\mathcal R(X)>0$. Then $\iota(X^m)=\iota(X)$ for
every $m>0$, by Corollary~\ref{C:rigidity>0 and powers}. In the same
way as above, since $M$ is a multiple of the orbit length of $X$ and
$M-t>0$, replacing $M$ (and thus $t$) by a multiple if necessary one
has $I(X)=\iota(\mathbf R_{[-M,M]})= \iota(\mathbf R_{M}) =
\iota(\Delta^{s} X^{M-t}) = \iota(X^{M-t}) = \iota(X)$.
\end{proof}

In order to obtain the claimed bound on rigid powers, we need to
investigate how the left normal form of $\mathbf C_{[k,m]}$ is
modified when we multiply it on the left by $C_k$. We actually show
the following, more general result.

\begin{proposition}\label{P:normal forms of products of C_m}
Let $X\in USS(X)$ with $\ell(X)>1$, and suppose that $X$ has a rigid
power. Let $t>1$ be such that $X^t\in USS(X^t)$. Consider $C_1\cdots
C_m=\Delta^k y_1\cdots y_s$ and $C_{m+1}\cdots C_{m+t}= \Delta^q z_1
\cdots z_r$ in left normal form. Then, the final $t-1$ factors in
the left normal form of $C_1\cdots C_{m+t}$ are precisely $z_2\cdots
z_r$.
\end{proposition}

\begin{proof}
We will need to use the factors $C_i$ and $R_i$ corresponding to the
element $Y=(\mathbf c^m(X))^t$. In order to avoid an excessive use
of indices, we will denote them by $C_i'$ and $R_i'$. That is,
$C_1'= \iota(Y)$, and the other elements $C_i'$ and $R_i'$ are
defined in the same way as the corresponding elements for $X$.

We know by Theorem~\ref{T:normal form of $C_m$ and $X^m$} that
$C_{m+1}\cdots C_{m+t}$ is equal to the product of the first $t$
factors (including $\Delta$'s) in the left normal form of $(\mathbf
c^m(X))^t \Delta^{-pt}$ (where $p=\inf(X)$). Hence $C_1'=
\iota((\mathbf c^m(X))^t)= \tau^{-q}(z_1)$, and $\tau^{-q}(z_2)
\preceq R_1'$.

Now one has
$$
    C_1\cdots C_m C_{m+1} \cdots C_{m+t}
$$
$$
 = C_1\cdots C_m \Delta^q z_1\cdots z_r
$$
$$
 = C_1\cdots C_m \tau^{-q}(z_1) \tau^{-q}(z_2) \Delta^q z_3\cdots z_r.
$$
By Proposition~\ref{P:chains} and the definition of absolute final
factors, one has $\varphi(C_1\cdots C_m)\succeq F(\mathbf c^m(X))$.
Recall from the proof of Proposition~\ref{P:SU is not emptyset} that
if $X^t\in USS(X^t)$ then $(\mathbf c^m(X))^t$ also belongs to its
ultra summit set, since the action of cycling or decycling any power
of $X$ (in particular $X$) preserves the set of powers of $X$
belonging to their ultra summit set. Moreover, since $X$ has a rigid
power, Theorem~\ref{T:pre-rigid for whole USS} implies that $\mathbf
c^m(X)$ also has a rigid power. We can then apply
Proposition~\ref{P:F(X)=F(X^t) if X is PA} to obtain $F(\mathbf
c^m(X))=F((\mathbf c^m(X))^t)=F(Y)$. Hence $\varphi(C_1\cdots
C_m)\succeq F(Y)$.

This yields the following:
$$
   \varphi(C_1\cdots C_m \tau^{-q}(z_1)) = \varphi(C_1\cdots C_m
   C_1') \succeq \varphi( F(Y) C_1' ) = \varphi(C_{-a+1}' \cdots  C_{-1}' C_0' C_1' )
$$
for $a$ big enough. But we know by Lemma~\ref{L:decomposition whole
orbit} that the decomposition
$$
   \varphi(C_{-a+1}' \cdots  C_{-1}' C_0' C_1' ) \; R_1'
$$
is left weighted. Hence
$$
  \varphi(C_1\cdots C_m \tau^{-q}(z_1))\; \tau^{-q}(z_2)
$$
is also left weighted, and the factors $z_2\cdots z_r$ are not
modified when computing the left normal form of $C_1\cdots C_{m+t}$.
\end{proof}

The fact that the left normal forms of $\mathbf C_{[-m,m]}$ are not
modified too much when one increases $m$, implies a strong property
on the initial factors of powers of $X$: they are comparable by
$\preceq$. First we need the following technical result.

\begin{lemma}\label{L:i(As) and i(A)}
Let $A=\Delta^p x_1\cdots x_r \in G$, and let $s$ be a simple
element. Then either $\iota(As)\preceq \iota(A)$ or $\iota(A)\preceq
\iota(As)$.
\end{lemma}

\begin{proof}
Suppose there is no unexpected $\Delta$ when multiplying $A$ by $s$,
that is, $\inf(As)=p$. Then $\iota(\Delta^p x_1\cdots x_r s)
=\tau^{-p}(x_1 t)$ for some (possibly trivial) simple element $t$.
Hence $\iota(A)=\tau^{-p}(x_1)\preceq \tau^{-p}(x_1t) =\iota(As)$.

Now suppose there is an unexpected $\Delta$, that is,
$\inf(As)=p+1$. Let $\alpha=(x_2\cdots x_r s)\wedge \Delta$. Then
$\iota(As)=\iota(\Delta^p x_1\cdots x_r s)=\iota(\Delta^p x_1
\alpha)$. Moreover, $x_1\alpha = \Delta \beta$ for some simple
element $\beta$, that is, $\alpha= \partial(x_1) \beta$. Since
$\alpha$ is simple, it follows that $\beta \preceq \partial^2(x_1) =
\tau(x_1)$. Therefore
$$
  \iota(As)=\iota(\Delta^p x_1\alpha)=\iota(\Delta^{p+1}\beta)
  \preceq \iota(\Delta^{p+1} \tau(x_1) ) = \tau^{-p}(x_1) =
  \iota(A).
$$
\end{proof}

\begin{proposition}
Let $X\in USS(X)$ with $\ell(X)>1$. Suppose that $X$ has a rigid
power, and let $t\geq 1$ such that $X^t\in USS(X^t)$. Then for every
$m\geq 1$, the simple elements $\iota(C_1\cdots C_m)$ and
$\iota(C_1\cdots C_{m+t})$ are comparable. That is, either
$$
 \iota(C_1\cdots C_m)\preceq \iota(C_1\cdots C_{m+t})
$$
or
$$
 \iota(C_1\cdots C_{m+t})\preceq \iota(C_1\cdots C_{m}).
$$
\end{proposition}

\begin{proof}
Write $C_1\cdots C_m = \Delta^k y_1\cdots y_s$ and $C_{m+1}\cdots
C_{m+t}= \Delta^q z_1\cdots z_r$. We know by
Proposition~\ref{P:normal forms of products of C_m} that
$$
  \varphi(y_1\cdots y_s \tau^{-q}(z_1)) \: \tau^{-q}(z_2)
$$
is left weighted. Hence
$$
  \iota(C_1\cdots C_{m+t}) = \iota(\Delta^k y_1\cdots y_s
  \tau^{-q}(z_1)),
$$
where $\tau^{-q}(z_1)$ is simple. By Lemma~\ref{L:i(As) and i(A)},
this implies that $\iota(C_1\cdots C_{m+t})$ is comparable to
$\iota(\Delta^k y_1\cdots y_s)=\iota(C_1\cdots C_{m})$, as we wanted
to show.
\end{proof}

We can finally state the result concerning the initial factors of
powers of $X$.

\begin{corollary}
Let $X\in USS(X)$ with $\ell(X)>1$, and suppose that $X$ has a rigid
power. If $X^m\in USS(X^m)$ for $m=1,\ldots,||\Delta||$, then the
set
$$
\{\iota(X), \iota(X^2), \cdots \iota(X^{||\Delta||})\}
$$
is totally ordered by $\preceq$ (although the total order given by
$\preceq$ does not necessarily coincide with the above enumeration).
\end{corollary}

\begin{proof}
We just need to recall from Theorem~\ref{T:normal form of $C_m$ and
$X^m$} that $\iota(X^m)=\iota(C_1\cdots C_m)$ for every $m\geq 1$,
and use the above result. Since every two elements are comparable by
$\preceq$, the set is totally ordered.
\end{proof}

\begin{corollary}\label{C:repetition of initial factors}
With the above conditions, there exist some integers $a,b$ with
$1\leq a<b\leq ||\Delta||$ such that $\iota(X^a)=\iota(X^b)$.
\end{corollary}

\begin{proof}
The length of a strict chain of simple elements $1\prec s_1 \prec
s_2 \prec \cdots \prec s_r \prec \Delta$ is bounded by $||\Delta||$.
Since the elements in $\{\iota(X), \iota(X^2), \cdots
\iota(X^{||\Delta||})\}$ are totally ordered, the lack of a repeated
pair would provide a chain of bigger length, which is not possible.
\end{proof}

It is important to notice that, in the sequence $\iota(X),
\iota(X^2), \iota(X^3),\ldots$, when one encounters the first
repetition, the sequence becomes periodic. And the period is the
distance between the two repeated elements. This is given by the
following result.

\begin{proposition}\label{P:chain of initial factors becomes
periodic} With the above conditions, if $\iota(X^a)=\iota(X^b)$ then
$\iota(X^{a+k})=\iota(X^{b+k})$ for every $k\geq 0$.
\end{proposition}

\begin{proof}
By hypothesis $\iota(C_1\cdots C_a)= \iota(X^a) =\iota(X^b) =
\iota(C_1\cdots C_b)$. Applying to this equality the transport
defined in~\cite{Gebhardt}, one obtains $\iota(C_2\cdots C_{a+1}) =
\iota(C_2\cdots C_{b+1})$. Let $\iota= \iota(C_2\cdots C_{a+1}) =
\iota(C_2\cdots C_{b+1})$. By Proposition~\ref{P:normal forms of
products of C_m}, if we multiply $C_2\cdots C_{a+1}$ or $C_2\cdots
C_{b+1}$ on the left by $C_1$, only their initial factors (which in
both cases are equal to $\iota$) are modified. Moreover, since
$\sup(C_1\cdots C_{a+1})=a+1$, it follows that $\sup(C_1 \iota)=2$.
Hence the initial factor of $C_1\cdots C_{a+1}$ is equal either to
$\Delta\wedge (C_1\iota)$ (if the infimum does not increase) or to
$(C_1\iota) \Delta^{-1}$ (if the infimum increases). In any case,
$\iota(C_1\cdots C_{a+1})= \iota(C_1 \iota)$. In the same way,
$\iota(C_1\cdots C_{b+1})= \iota(C_1 \iota)$, hence
$\iota(X^{a+1})=\iota(X^{b+1})$. Induction on $k$ finishes the
proof.
\end{proof}

The above results can be used to bound the smallest power of $X$
having some rigidity.

\begin{proposition}\label{P:bound on 1-rigidity}
Let $X\in USS(X)$ with $\ell(X)>1$, and suppose that $X$ has a rigid
power. If $X^t\in USS(X^t)$ for $t=1,\ldots, ||\Delta||$, then
$\mathcal R(X^m)>0$ for some positive $m<||\Delta||$.
\end{proposition}

\begin{proof}
We know by Corollary~\ref{C:repetition of initial factors} that
$\iota(X^a) = \iota(X^b) $ for some $1\leq a < b \leq ||\Delta||$,
and by Proposition~\ref{P:chain of initial factors becomes periodic}
that the sequence $\iota(X^a),\iota(X^{a+1}),\iota(X^{a+2}),\ldots$
is periodic of period $d=b-a$.

Since the interval $[a,b]$ has length $d$, there exists a unique
$m$, $a\leq m < b$, which is a multiple of $d$. Then
$\iota(X^m)=\iota(X^{m+m}) = \iota((X^m)^2)$.  Hence, by
Lemma~\ref{L:k-rigidity}, we will have $\mathcal R(X^m)>0$ if we
show that $\inf(X^{2m})= 2\: \inf(X^m)$.

Suppose that $X^m=\Delta^p y_1\cdots y_r$. Since $m<b\leq
||\Delta||$ one has $X^m\in USS(X^m)$, so $\Delta \not\preceq
y_2\cdots y_r \tau^{-p}(y_1)$, which implies $\Delta \not \preceq
y_2\cdots y_r \tau^{-p}(y_1\cdots y_r)$ since $\tau^{-p}(y_1)$ is
the biggest simple prefix of $\tau^{-p}(y_1\cdots y_r)$. Then
$$
 y_1\Delta \not
\preceq y_1 y_2\cdots y_r \tau^{-p}(y_1\cdots y_r).
$$
But if $\inf(X^{2m})=2p+1$ then, since $X^{2m} = (X^m)^2 = \Delta^p
y_1\cdots y_r \Delta^p y_1\cdots y_r$, we would have $\Delta\preceq
\tau^p(y_1\cdots y_r) y_1\cdots y_r$, and since
$\iota(X^{2m})=\iota(X^m)=\tau^{-p}(y_1)$, it would follow that
$\Delta \tau^{p+1}(y_1)\preceq \tau^p(y_1\cdots y_r) y_1\cdots y_r$.
Applying $\tau^{-p}$ to this inequality, we would obtain $\Delta
\tau(y_1) = y_1\Delta \preceq y_1\cdots y_r \tau^{-p}(y_1\cdots
y_r)$, a contradiction. Therefore $\inf(X^{2m})=2p = \inf(X^m)$, and
since $\iota(X^{2m})=\iota(X^m)$ it follows from
Lemma~\ref{L:k-rigidity} that $\mathcal R(X^m)>0$.
\end{proof}

\begin{theorem}\label{T:bound on rigidity}
Let $X\in USS(X)$ with $\ell(X)>1$, and suppose that $X$ has a rigid
power. If $X^t\in USS(X^t)$ for every $t$ such that $-||\Delta||
\leq t \leq ||\Delta||$, then there is some $m<||\Delta||^2$ such
that $X^m$ is rigid.
\end{theorem}

\begin{proof}
By Proposition~\ref{P:bound on 1-rigidity}, $\mathcal R(X^p)>0$ for
some $0<p<||\Delta||$. Applying Proposition~\ref{P:bound on
1-rigidity} to $X^{-1}$, one obtains that $\mathcal R(X^{-q})>0$ for
some $0<q<||\Delta||$.

Let $m={\rm lcm} (p,q)\leq pq < ||\Delta||^2$. Since $X^m$ is a
power of $X^p$, it follows from Corollary~\ref{C:1-rigid in USS has
powers in closed orbit} that $X^m$ it belongs to a closed orbit
under cycling, hence it has maximal infimum in its conjugacy
class~\cite{El-M}. In the same way, since $X^{-m}$ is a power of
$X^{-q}$, it follows from Corollary~\ref{C:1-rigid in USS has powers
in closed orbit} that $X^{-m}$ belongs to a closed orbit under
cycling, thus $X^m$ belongs to a closed orbit under decycling, and
hence it has minimal supremum in its conjugacy class~\cite{El-M}.
Therefore $X^m\in SSS(X^m)$, and since it belongs to a closed orbit
under cycling, $X^m\in USS(X^m)$.

Moreover, by Proposition~\ref{P:rigidity and powers}, $\mathcal
R(X^m)>0$ and $\mathcal R(X^{-m})>0$, since they are powers of $X^p$
and $X^{-q}$. This implies by Proposition~\ref{P:1-rigid PA and
absolute factors} that $\iota(X^m)=I(X^m)$ and
$\varphi(X^m)=F(X^m)$.  Since $F(X^m)I(X^m)$ is left weighted by
Proposition~\ref{P:FI left weighted}, it follows that
$\varphi(X^m)\iota(X^m)$ is left weighted as written, hence $X^m$ is
rigid, as we wanted to show.
\end{proof}

{\bf Remark:} The proof of the above result is based on the fact
that if an element $X$ is such that $\mathcal R(X)>0$ and $\mathcal
R(X^{-1})>0$, and if $X$ has a rigid power, then $X$ is already
rigid. The hypothesis of $X$ having a rigid power is necessary,
since we could have $X$ and $X^{-1}$ with some rigidity without $X$
being rigid, even if $X\in SU(X)$. For instance, if we consider the
reducible braid $X=\sigma_1 \sigma_3 \cdot \sigma_3 \in B_4$, we
have $\mathcal R(X)=\mathcal R(X^{-1})=1/2$, but neither $X$ nor any
power of $X$ is rigid, since $X^m = (\sigma_1 \sigma_3 )^m \cdot
\sigma_3^m$. In this case $X\in SU(X)$ and $\mathcal R(X^m)=1/2$ for
every $m\neq 0$.

\ms

We can also find a bound for the smallest rigid power in the case
$\ell(X)=1$, thanks to the following result.

\begin{lemma}\label{L:powers of length 1}
Let $X\in G$ with $\ell(X)=1$. If $X^t\in USS(X^t)$ for $t=1,\ldots,
||\Delta||$, then $\ell(X^m)\neq 1$ for some $m\leq ||\Delta||$.
\end{lemma}

\begin{proof}
Suppose that $\ell(X^t)=1$ for $t=1,\ldots, ||\Delta||$. We will
show that the set $\{\iota(X),\iota(X^2),\ldots,
\iota(X^{||\Delta||})\}$ is totally ordered by showing that any two
elements in that set are comparable. Indeed, given $s,t\in
\{1,\ldots, ||\Delta||\}$ with $s<t$, we have $\ell(X^{t-s})=1$,
hence $\iota(X^t) = \iota(X^s X^{t-s}) = \iota(X^s \iota(X^{t-s}))$.
Since $\iota(X^{t-s})$ is simple, Lemma~\ref{L:i(As) and i(A)}
implies that either $\iota(X^t)\preceq \iota(X^s)$ or
$\iota(X^s)\preceq \iota(X^t)$. Therefore,
$\{\iota(X),\iota(X^2),\ldots, \iota(X^{||\Delta||})\}$ is a totally
ordered set of proper simple elements, thus $\iota(X^a)=\iota(X^b)$
for some $1\leq a < b \leq ||\Delta||$. But since
$\ell(X^a)=\ell(X^b)=1$, this means that $X^{b-a}$ is a power of
$\Delta$, a contradiction. Therefore, $\ell(X^m)\neq 1$ for some
$m\leq ||\Delta||$.
\end{proof}

We can finally remove the hypothesis $\ell(X)>1$ in order to bound
the rigid power of an element.

\begin{theorem}\label{T:bound on rigidity for any length}
Let $X\in USS(X)$, and suppose that $X$ has a rigid power. If
$X^t\in USS(X^t)$ for every $t$ such that $-||\Delta||^2 \leq t \leq
||\Delta||^2$, then there is some $m<||\Delta||^3$ such that $X^m$
is rigid.
\end{theorem}

\begin{proof}
If $\ell(X)>1$ the result follows from Theorem~\ref{T:bound on
rigidity}. If $\ell(X)=1$, Lemma~\ref{L:powers of length 1} implies
that $\ell(X^r)\neq 1$ for some $r\leq ||\Delta||$. We cannot have
$\ell(X^r)=0$, otherwise $X$ would be periodic and would not have a
rigid power. Hence $\ell(X^r)>1$, and $(X^r)^t \in USS((X^r)^t)$ for
$-||\Delta|| \leq t \leq ||\Delta||$. The hypothesis of
Theorem~\ref{T:bound on rigidity} are then satisfied by $X^r$, hence
$(X^r)^s$ is rigid for some $s< ||\Delta||^2$. Therefore $X^{rs}$ is
rigid with $rs < ||\Delta||^3$.
\end{proof}

In the case of braid groups, the above result implies the following.

\begin{theorem}\label{T:bound for pseudo-Anosov}
If $X\in B_n$ is a pseudo-Anosov braid, then $USS(X^m)$ consists of
rigid braids, for $m<||\Delta||^3$. Moreover, if the canonical
length of the elements in $USS(X)$ is greater than 1, then
$m<||\Delta||^2$.
\end{theorem}

\begin{proof}
By Corollary~\ref{C:pA has a rigid power}, some conjugate $Y$ of $X$
has a rigid power. Moreover, one can choose $Y\in SU(X)$, hence by
Theorem~\ref{T:bound on rigidity for any length} $Y^m$ is rigid for
$m<||\Delta||^3$ (and $m<||\Delta||^2$ if $\ell(Y)>1$). In the proof
of Theorem~\ref{T:bound on rigidity for any length} we see that we
can assume $\ell(Y^m)>1$, hence it follows from Theorem~\ref{T:X
rigid implies every Y in USS(X) rigid} that $USS(Y^m)$ consists of
rigid elements.
\end{proof}

To summarize the consequences for pseudo-Anosov braids, we can solve
the CDP/CSP problem for two pseudo-Anosov elements $X,Y\in B_n$
using rigid braids. We just need to compute an element in $USS(X^t)$
for each $t=1,2,\ldots$ until we find one of them, say $\widetilde
X^m\in USS(X^m)$ which is rigid and has canonical length greater
than one. By the above result, $m<||\Delta||^3$. Then all elements
in $USS(X^m)$ will be rigid, so the computation of $USS(X^m)$ is
easier than in the general case, as will be seen in~\cite{BGGM-II},
and we will possibly be able to bound the size of $USS(X^m)$.
Moreover, since pseudo-Anosov braids have unique roots, if one
solves the CDP/CSP for $X^m$ and $Y^m$, finding some conjugating
element $Z$, then $Z$ is also a conjugating element for $X$ and $Y$,
so this solves the CDP/CSP for $X$ and $Y$.

%
%
%

\footnotesize
\begin{tabular}{lll}
 {\bf Joan S. Birman}        &  {\bf Volker Gebhardt}  &  {\bf Juan Gonz\'alez-Meneses} \\
 Department of Mathematics, & School of Computing and Mathematics, & Departamento de \'Algebra, \\
 Barnard College andColumbia University, & University of Western Sydney, & Universidad de Sevilla,\\
 2990 Broadway, & Locked Bag 1797, &  Apdo. 1160, \\
 New York, New York 10027, USA. & Penrith South DC NSW 1797, Australia, & 41080 Sevilla, Spain.\\
 {\tt jb@math.columbia.edu} & {\tt v.gebhardt@uws.edu.au} & {\tt meneses@us.es}
\end{tabular}

\end{document}